\documentclass[sigconf,screen,nonacm,acmsmall]{acmart}
\pdfoutput=1

\let\savedsection\section

\usepackage{natbib}
\usepackage{lscape}  
\usepackage{subcaption}
\usepackage{scalerel}
\usepackage{hyphenat}
\usepackage{mathpartir}
\usepackage{braket}
\usepackage{cancel}
\usepackage{bookmark}
\usepackage{graphicx}
\usepackage{mathtools}
\usepackage[T1]{fontenc}
\usepackage{libertine}
\usepackage{transparent}
\usepackage{framed}

\newenvironment{recycle}
  {\leavevmode\color{gray}}
  {}

\usepackage{tikz}
\usepackage{tikz-cd}
\usepackage{quiver}

\tikzset{every picture/.style={line width=0.85pt}}

\usepackage[leftmargin=1em,rightmargin=1ex,vskip=1ex]{quoting}
\usepackage[utf8]{inputenc}
\usepackage{cmll}

\newcommand\scalemath[3]{\scalebox{#1}[#2]{\mbox{\ensuremath{\displaystyle #3}}}}

\newcommand{\leftarrowtip}{\ensuremath{\tikz\draw[line width=0.5pt,->] (10pt,0) -- (0,0);}}
\newcommand{\leftarrowtailnotip}{\ensuremath{\tikz\draw[line width=0.5pt,-<] (0,0) -- (10pt,0);}}

\newcommand{\unicodeStar}{\ensuremath{\star}}
\DeclareUnicodeCharacter{2605}{\unicodeStar}

\DeclareUnicodeCharacter{21D2}{\ensuremath{\Rightarrow}}
\DeclareUnicodeCharacter{2218}{\ensuremath{\circ}}
\DeclareUnicodeCharacter{2022}{\ensuremath{\bullet}}
\DeclareUnicodeCharacter{2219}{\ensuremath{\bullet}}
\DeclareUnicodeCharacter{2026}{\ensuremath{\dots}}
\DeclareUnicodeCharacter{2208}{\ensuremath{\in}}
\DeclareUnicodeCharacter{2192}{\ensuremath{\to}}
\DeclareUnicodeCharacter{2190}{\ensuremath{\leftarrowtip}}
\DeclareUnicodeCharacter{2919}{\ensuremath{\leftarrowtailnotip}}

\DeclareUnicodeCharacter{00D7}{\ensuremath{\times}}
\DeclareUnicodeCharacter{00B7}{\ensuremath{\cdot}}
\DeclareUnicodeCharacter{222B}{\ensuremath{\int}}
\DeclareUnicodeCharacter{22A4}{\ensuremath{\top}}
\DeclareUnicodeCharacter{22A5}{\ensuremath{\bot}}
\DeclareUnicodeCharacter{2264}{\ensuremath{\leq}}

\newcommand{\unicodecolon}{\ensuremath{\colon}}
\DeclareUnicodeCharacter{FE55}{\unicodecolon}
\newcommand{\unicodeleftpar}{\ensuremath{\left(}}
\DeclareUnicodeCharacter{27EE}{\unicodeleftpar}
\newcommand{\unicoderightpar}{\ensuremath{\right)}}
\DeclareUnicodeCharacter{27EF}{\unicoderightpar}
\DeclareUnicodeCharacter{2260}{\neq}
\DeclareUnicodeCharacter{22A9}{\Vdash}
\DeclareUnicodeCharacter{2237}{\proportion}
\DeclareUnicodeCharacter{2124}{\mathbb{Z}}
\DeclareUnicodeCharacter{27E8}{\langle}
\DeclareUnicodeCharacter{27E9}{\rangle}
\DeclareUnicodeCharacter{21A6}{\mapsto}
\DeclareUnicodeCharacter{22A2}{\vdash}
\DeclareUnicodeCharacter{2090}{\ensuremath{{}_a}}
\DeclareUnicodeCharacter{A71B}{{}^\uparrow}
\DeclareUnicodeCharacter{A71C}{{}^\downarrow}
\DeclareUnicodeCharacter{27E6}{\llbracket}
\DeclareUnicodeCharacter{27E7}{\rrbracket}

\DeclareUnicodeCharacter{2460}{\textcircled{\small{1}}}
\DeclareUnicodeCharacter{2461}{\textcircled{\small{2}}}
\DeclareUnicodeCharacter{2462}{\textcircled{\small{3}}}
\DeclareUnicodeCharacter{2463}{\textcircled{\small{4}}}
\DeclareUnicodeCharacter{2464}{\textcircled{\small{5}}}
\DeclareUnicodeCharacter{2465}{\textcircled{\small{6}}}
\DeclareUnicodeCharacter{2466}{\textcircled{\small{7}}}
\DeclareUnicodeCharacter{24EA}{\textcircled{\small{0}}}

\newcommand{\unicoderightcircle}{\ensuremath{\RIGHTcircle}}
\DeclareUnicodeCharacter{25D1}{\unicoderightcircle}
\newcommand{\unicodeleftcircle}{\ensuremath{\LEFTcircle}}
\DeclareUnicodeCharacter{25D0}{\unicodeleftcircle}
\DeclareUnicodeCharacter{229B}{\circledast}

\newcommand{\unicodebbA}{\ensuremath{\mathbb{A}}}
\DeclareUnicodeCharacter{1D538}{\unicodebbA}
\newcommand{\unicodebbB}{\ensuremath{\mathbb{B}}}
\DeclareUnicodeCharacter{1D539}{\unicodebbB}
\newcommand{\unicodebbC}{\ensuremath{\mathbb{C}}}
\DeclareUnicodeCharacter{2102}{\unicodebbC}
\DeclareUnicodeCharacter{1D53B}{\ensuremath{\mathbb{D}}}
\DeclareUnicodeCharacter{2115}{\ensuremath{\mathbb{N}}}
\DeclareUnicodeCharacter{211D}{\ensuremath{\mathbb{R}}}
\DeclareUnicodeCharacter{1D543}{\ensuremath{\mathbb{L}}}
\newcommand\UnicodeBlackboardP{\ensuremath{\mathbf{P}}} \DeclareUnicodeCharacter{2119}{\UnicodeBlackboardP}
\DeclareUnicodeCharacter{211A}{\ensuremath{\mathbb{Q}}}
\DeclareUnicodeCharacter{1D544}{\ensuremath{\mathbb{M}}}
\DeclareUnicodeCharacter{1D54C}{\ensuremath{\mathbb{U}}}
\DeclareUnicodeCharacter{1D54D}{\ensuremath{\mathbf{V}}}
\DeclareUnicodeCharacter{1D54E}{\ensuremath{\mathbb{W}}}
\DeclareUnicodeCharacter{1D542}{\ensuremath{\mathbb{K}}}
\DeclareUnicodeCharacter{1D546}{\ensuremath{\mathbb{O}}}
\DeclareUnicodeCharacter{1D540}{\ensuremath{\mathbb{I}}}
\DeclareUnicodeCharacter{1D54A}{\ensuremath{\mathbb{S}}}
\DeclareUnicodeCharacter{1D53C}{\ensuremath{\mathbb{E}}}

\DeclareUnicodeCharacter{1D405}{\ensuremath{\mathbf{F}}}
\DeclareUnicodeCharacter{1D406}{\ensuremath{\mathbf{G}}}
\DeclareUnicodeCharacter{1D413}{\ensuremath{\mathbf{T}}}
\DeclareUnicodeCharacter{1D5E4}{\ensuremath{\mathbf{Q}}}

\newcommand{\unicodecalS}{\ensuremath{\mathcal{S}}}
\newcommand{\unicodecalT}{\ensuremath{\mathcal{T}}}
\newcommand{\unicodecalC}{\ensuremath{\mathcal{C}}}

\newcommand{\unicodecalX}{\ensuremath{\mathcal{X}}}
\newcommand{\unicodecalN}{\ensuremath{\mathcal{N}}}
\newcommand{\unicodecalE}{\ensuremath{\mathcal{E}}}
\DeclareUnicodeCharacter{1D4D4}{\unicodecalE}
\DeclareUnicodeCharacter{1D4D2}{\unicodecalC}

\DeclareUnicodeCharacter{1D4DE}{\mathcal{O}}
\DeclareUnicodeCharacter{1D4AA}{\mathcal{O}}
\DeclareUnicodeCharacter{210B}{\mathcal{H}}
\DeclareUnicodeCharacter{1D4E2}{\unicodecalS}
\DeclareUnicodeCharacter{1D4E3}{\unicodecalT}
\DeclareUnicodeCharacter{1D4E7}{\unicodecalX}
\DeclareUnicodeCharacter{1D4D0}{\ensuremath{\mathcal{A}}}
\DeclareUnicodeCharacter{1D4D1}{\ensuremath{\mathcal{B}}}
\DeclareUnicodeCharacter{1D4D6}{\ensuremath{\mathcal{G}}}
\DeclareUnicodeCharacter{1D4D7}{\ensuremath{\mathcal{H}}}
\DeclareUnicodeCharacter{1D4DB}{\ensuremath{\mathcal{L}}}

\DeclareUnicodeCharacter{1D4DD}{\unicodecalN}
\DeclareUnicodeCharacter{1D4B1}{\ensuremath{\mathcal{V}}}
\DeclareUnicodeCharacter{1D4E5}{\ensuremath{\mathcal{V}}}
\DeclareUnicodeCharacter{1D4E6}{\ensuremath{\mathcal{W}}}
\DeclareUnicodeCharacter{1D4E4}{\ensuremath{\mathcal{U}}}

\DeclareUnicodeCharacter{03B1}{\alpha}
\DeclareUnicodeCharacter{03B2}{\beta}
\DeclareUnicodeCharacter{03BC}{\mu}
\DeclareUnicodeCharacter{03B4}{\delta}
\DeclareUnicodeCharacter{03B5}{\varepsilon}
\DeclareUnicodeCharacter{03B7}{\eta}
\DeclareUnicodeCharacter{03BB}{\lambda}
\DeclareUnicodeCharacter{03C1}{\rho}
\DeclareUnicodeCharacter{03C8}{\psi}
\DeclareUnicodeCharacter{03C4}{\tau}
\DeclareUnicodeCharacter{03A8}{\Psi}
\DeclareUnicodeCharacter{03C3}{\sigma}
\DeclareUnicodeCharacter{03C6}{\varphi}
\DeclareUnicodeCharacter{03A6}{\Phi}
\DeclareUnicodeCharacter{03A3}{\Sigma}
\DeclareUnicodeCharacter{03D5}{\phi}
\DeclareUnicodeCharacter{03B8}{\theta}
\DeclareUnicodeCharacter{03C0}{\ensuremath{\pi}}
\DeclareUnicodeCharacter{0393}{\Gamma}
\DeclareUnicodeCharacter{0394}{\Delta}
\DeclareUnicodeCharacter{03BA}{\kappa}
\DeclareUnicodeCharacter{03BD}{\nu}
\DeclareUnicodeCharacter{03A9}{\Omega}
\DeclareUnicodeCharacter{03C9}{\omega}
\DeclareUnicodeCharacter{03D2}{\Upsilon}
\DeclareUnicodeCharacter{03C5}{\upsilon}
\DeclareUnicodeCharacter{22A0}{\boxtimes}
\DeclareUnicodeCharacter{1D412}{\ensuremath{\mathbf{S}}}

\DeclareFontFamily{U}{mathb}{\hyphenchar\font45}
\DeclareSymbolFont{mathb}{U}{mathb}{m}{n}
\DeclareMathSymbol{\sqbullet}{\mathbin}{mathb}{"0D}
\DeclareUnicodeCharacter{25A0}{\sqbullet}
\DeclareUnicodeCharacter{25AA}{\mathbin{\sqbullet}}

\DeclareUnicodeCharacter{2205}{\varnothing}
\DeclareUnicodeCharacter{2254}{\coloneq}
\DeclareUnicodeCharacter{25AB}{\blacksquare}
\DeclareUnicodeCharacter{03B3}{\gamma}

\DeclareUnicodeCharacter{112A9}{\mathbin{\pmb{;}}}

\newcommand{\hirayo}{\scaleobj{0.9}{\text{\usefont{U}{min}{m}{n}\symbol{'210}}}}
\DeclareUnicodeCharacter{3088}{\hirayo}
\DeclareFontFamily{U}{min}{}
\DeclareFontShape{U}{min}{m}{n}{<-> udmj30}{}

\newcommand\UnicodeWhiteRightPointingSmallTriangle{\triangleright}
\DeclareUnicodeCharacter{25B9}{\mathbin{\UnicodeWhiteRightPointingSmallTriangle}}
\newcommand\UnicodeWhiteDownPointingSmallTriangle{\triangledown}
\DeclareUnicodeCharacter{25BF}{\mathbin{\UnicodeWhiteDownPointingSmallTriangle}}
\newcommand\UnicodeWhiteUpPointingSmallTriangle{\scalemath{1}{-1}{{}^{\triangledown}}}
\DeclareUnicodeCharacter{25B5}{\mathbin{\UnicodeWhiteUpPointingSmallTriangle}}

\DeclareUnicodeCharacter{2080}{\ensuremath{{}_0}}
\DeclareUnicodeCharacter{2081}{\ensuremath{{}_1}}
\DeclareUnicodeCharacter{2082}{\ensuremath{{}_2}}
\DeclareUnicodeCharacter{2083}{\ensuremath{{}_3}}

\DeclareUnicodeCharacter{1D62}{\ensuremath{{}_i}}
\DeclareUnicodeCharacter{2C7C}{\ensuremath{{}_j}}
\DeclareUnicodeCharacter{02B3}{\ensuremath{{}^r}}
\DeclareUnicodeCharacter{02E1}{\ensuremath{{}^\ell}}
\DeclareUnicodeCharacter{1D48}{\ensuremath{{}^d}}
\DeclareUnicodeCharacter{1D50}{\ensuremath{{}^m}}
\DeclareUnicodeCharacter{1D58}{\ensuremath{{}^u}}
\DeclareUnicodeCharacter{209A}{\ensuremath{{}_p}}
\DeclareUnicodeCharacter{2096}{\ensuremath{{}_k}}

\DeclareUnicodeCharacter{2245}{\ensuremath{\cong}}
\DeclareUnicodeCharacter{2286}{\subseteq}

\DeclareUnicodeCharacter{22C5}{\cdot}
\DeclareUnicodeCharacter{25C3}{\ensuremath{\triangleleft}}
\DeclareUnicodeCharacter{25B9}{\ensuremath{\triangleright}}

\DeclareUnicodeCharacter{2217}{\ast}
\DeclareUnicodeCharacter{039E}{\Xi}
\DeclareUnicodeCharacter{2295}{\oplus}
\DeclareUnicodeCharacter{2297}{\otimes}
\DeclareUnicodeCharacter{214B}{\parr}
\DeclareUnicodeCharacter{2298}{\oslash}
\DeclareUnicodeCharacter{25C0}{\mathbin{\blacktriangleleft}}
\DeclareUnicodeCharacter{25C1}{\mathbin{\vartriangleleft}}
\DeclareUnicodeCharacter{22B3}{\mathbin{\triangleright}}
\DeclareUnicodeCharacter{22B2}{\mathbin{\triangleleft}}
\DeclareUnicodeCharacter{FF5C}{\mid}
\DeclareUnicodeCharacter{227A}{\mathbin{\prec}}
\DeclareUnicodeCharacter{227B}{\mathbin{\succ}}
\DeclareUnicodeCharacter{22A3}{\mathbin{\dashv}}
\DeclareUnicodeCharacter{219D}{\ensuremath{\leadsto}}
\DeclareUnicodeCharacter{2191}{\ensuremath{\uparrow}}
\DeclareUnicodeCharacter{1361}{\colon}

\DeclareUnicodeCharacter{29D1}{\mathrel{\multimapdotbothB}}
\DeclareUnicodeCharacter{29D2}{\mathrel{\multimapdotbothA}}
\DeclareUnicodeCharacter{22C4}{\mathbin{\diamond}}
\DeclareUnicodeCharacter{226B}{\mathrel{\gg}}
\DeclareUnicodeCharacter{25A1}{\mathbin{\square}}
\DeclareUnicodeCharacter{266F}{\sharp}

\DeclareUnicodeCharacter{2113}{\ell}
\DeclareUnicodeCharacter{2261}{\equiv}
\DeclareUnicodeCharacter{2099}{_n}
\DeclareUnicodeCharacter{2098}{_m}
\DeclareUnicodeCharacter{1D5AD}{\ensuremath{\mathsf{N}}}
\DeclareUnicodeCharacter{1D5E1}{\ensuremath{\mathbf{N}}}

\DeclareUnicodeCharacter{1D4DF}{\mathcal{P}}
\DeclareUnicodeCharacter{1D415}{\mathbf{V}}
\DeclareUnicodeCharacter{1D400}{\mathbf{A}}
\DeclareUnicodeCharacter{1D6C2}{\pmb{\alpha}}
\DeclareUnicodeCharacter{1D6C3}{\pmb{\beta}}
\DeclareUnicodeCharacter{1D6C4}{\pmb{\gamma}}
\DeclareUnicodeCharacter{1D6C5}{\pmb{\delta}}
\DeclareUnicodeCharacter{1D6C6}{\pmb{\epsilon}}
\DeclareUnicodeCharacter{1D6C8}{\pmb{\eta}}
\DeclareUnicodeCharacter{1D6DA}{\pmb{\omega}}

\newcommand\mydots{\makebox[0.6em][c]{.\hfil.\hfil.}}
\DeclareUnicodeCharacter{2026}{\mydots}
\DeclareUnicodeCharacter{226B}{\gg}

\DeclareUnicodeCharacter{22C9}{\ltimes}
\DeclareUnicodeCharacter{22CA}{\rtimes}

\usepackage[only,fatsemi]{stmaryrd}
\newcommand{\unicodeRelationalComposition}{\fatsemi}
\DeclareUnicodeCharacter{2A3E}{\unicodeRelationalComposition}
\DeclareUnicodeCharacter{2983}{\llparenthesis}
\DeclareUnicodeCharacter{2984}{\rrparenthesis}

\DeclareUnicodeCharacter{2014}{\,---\,} %

\IfFileExists{noappendix.token}%
{\usepackage[appendix=strip,bibliography=common]{apxproof}}%
{\usepackage[appendix=append,bibliography=common]{apxproof}}

\usepackage{xcolor}

\definecolor{nordred}{HTML}{bf616a}
\definecolor{bordeaux}{HTML}{821529}
\definecolor{bluelink}{HTML}{003399}
\definecolor{nordred}{HTML}{bf616a}
\definecolor{nordblue}{HTML}{81a1c1}
\definecolor{norddarkblue}{HTML}{5e81ac}
\definecolor{nordgreen}{HTML}{a3be8c}
\definecolor{nordnight}{HTML}{4c566a}
\newcommand{\cRed}[1]{{\color{nordred}#1}}
\newcommand{\cBlu}[1]{{\color{norddarkblue}#1}}

\AtEndPreamble{
  \RequirePackage{hyperref}
  \RequirePackage{cleveref}
  \crefname{theorem}{Theorem}{Theorems}
  \crefname{lemma}{Lemma}{Lemmas}
  \crefname{proposition}{Proposition}{Propositions}
  \crefname{corollary}{Corollary}{Corollaries}
  \crefname{definition}{Definition}{Definitions}
  \crefname{example}{Example}{Examples}

  \hypersetup{
    breaklinks = true,
    linktocpage,
    colorlinks = true,
    linkcolor = nordnight,
    citecolor = nordgreen,
    urlcolor = nordblue
  }
} %
\makeatletter
\newcommand{\nicelinktarget}[1]{\Hy@raisedlink{\hypertarget{#1}{}}}
\makeatother

\usepackage{hyphenat}
\renewcommand\Set{\hyperlink{linkSet}{\mathbf{Set}}}

\newcommand\id{\mathrm{id}}

\usepackage[xcolor,no patch,hyperref,quotation=false,electronic]{knowledge}
\knowledgeconfigure{notion}
\knowledgestyle{notion}{color=nordnight}
\knowledgestyle{intro notion}{emphasize,color=nordnight}

\knowledge{notion}
| normal magmad
| normal magmads
| Normal magmad
| Normal magmads
| normal monad
| normal monads
| Normal monad
| Normal monads

\knowledge{notion}
| string diagrams for commutative magmoids
| String diagrams for commutative magmoids

\knowledge{notion}
| right do-notation
| Right do-notation

\knowledge{notion}
| left do-notation
| Left do-notation

\knowledge{notion}
| disintegration axiom
| disintegration 
| Disintegration axiom
| Disintegration 

\knowledge{notion}
| uniformity
| Uniformity
| uniformity axiom
| Uniformity axiom
| uniformly traced distributive multicategory
| uniformly traced distributive multicategories
| Uniformly traced distributive multicategory
| Uniformly traced distributive multicategories

\knowledge{notion}
| associating morphism
| Associating morphism
| associating morphisms
| Associating morphisms
| associating
| Associating

\knowledge{notion}
| distributive sesquilaw
| distributive sesquilaws
| Distributive sesquilaw
| Distributive sesquilaws
| sesquilaw
| sesquilaws
| Sesquilaw
| Sesquilaws

\knowledge{notion}
| associative
| Associative

\knowledge{notion}
| sesquilaw homomorphism
| sesquilaw homomorphisms
| Sesquilaw homomorphism
| Sesquilaw homomorphisms

\knowledge{notion}
| monoidal sesquilaw
| monoidal sesquilaws
| Monoidal sesquilaw
| Monoidal sesquilaws

\knowledge{notion}
| copy-discard magmad
| Copy-discard magmad
| copy-discard magmads
| Copy-discard magmads
| affine-relevant magmad
| Affine-relevant magmad
| affine-relevant magmads
| Affine-relevant magmads

\knowledge{notion}
| normalized distribution
| normalized distributions
| Normalized distribution
| Normalized distributions

\knowledge{notion}
| quasitotal magmad
| Quasitotal magmad
| quasitotal magmads
| Quasitotal magmads

\knowledge{notion}
| post-selection
| Post-selection
\newcommand{\psel}{\kl[post-selection]{\ensuremath{\mathsf{p}}}}
\newcommand{\isel}{\kl[post-selection]{\ensuremath{\mathsf{i}}}}

\knowledge{notion}
| normalized probabilistic channels
| normalized probabilistic channel
| Normalized probabilistic channels
| Normalized probabilistic channel
| Normalized kernels
| Normalized kernel
| normalized kernel
| normalized kernels

\knowledge{notion}
| Continuous normalized kernels
| Continuous normalized kernel
| continuous normalized kernel
| continuous normalized kernels
\newcommand{\BorelNorm}{\kl[continuous normalized kernels]{\mathbf{BorelNorm}}}

\knowledge{notion}
| bracketed division
| Bracketed division
\newcommand{\bdiv}[2]{\frac{#1}{#2}\!{\raisebox{0.0pt}{\hspace{0.0pt}$\oslash$}}}

\knowledge{notion}
| monoidal natural transformation
| monoidal natural transformations
| Monoidal natural transformation
| Monoidal natural transformations
| monoidal transformation
| monoidal transformations
| Monoidal transformation
| Monoidal transformations

\knowledge{notion}
| monad
| monads
| Monad
| Monads

\knowledge{notion}
| distribution
| Distribution
| distributions
| Distributions

\knowledge{notion}
| non-empty powerset monad
| Non-empty powerset monad
| non-empty powerset
| Non-empty powerset
| Affine powerset monad
| affine powerset monad
| Affine powerset
| affine powerset
\DeclareUnicodeCharacter{1D411}{\kl[non-empty powerset monad]{\mathbf{R}}} %

\knowledge{notion}
| non-empty relation
| non-empty relations
| Non-empty relation
| Non-empty relations
| affine relation
| affine relations
| Affine relation
| Affine relations

\knowledge{notion}
| subaffine relation
| subaffine relations
| Subaffine relation
| Subaffine relations

\knowledge{notion}
| partial affine relation
| partial affine relations
| Partial affine relation
| Partial affine relations

\knowledge{notion}
| monoidal monad
| monoidal monads
| Monoidal monad
| Monoidal monads

\knowledge{notion}
| monoidal non-associative monad
| monoidal non-associative monads
| Monoidal non-associative monad
| Monoidal non-associative monads

\knowledge{notion}
| natural transformation
| natural transformations
| Natural transformation
| Natural transformations

\knowledge{notion}
| symmetric monoidal categories
| symmetric monoidal category
| Symmetric monoidal categories
| Symmetric monoidal category

\knowledge{notion}
| distributive law
| distributive laws
| Distributive law
| Distributive laws

\knowledge{notion}
| monoidal distributive law
| monoidal distributive laws
| Monoidal distributive law
| Monoidal distributive laws

\knowledge{notion}
| almost distributive law
| almost distributive laws
| Almost distributive law
| Almost distributive laws
| almost-distributive law
| almost-distributive laws
| Almost-distributive law
| Almost-distributive laws

\knowledge{notion}
| term
| terms
| Term
| Terms

\knowledge{notion}
| interchange
| Interchange
| interchange axiom
| Interchange axiom

\knowledge{notion}
| magmad
| Magmad
| magmads
| Magmads

\knowledge{notion}
| monoidal magmad
| Monoidal magmad
| monoidal magmads
| Monoidal magmads

\knowledge{notion}
| commutative magmad
| Commutative magmad
| commutative magmads
| Commutative magmads
| commutative monad
| Commutative monad
| commutative monads
| Commutative monads

\knowledge{notion}
| left-exchanging magmad
| Left-exchanging magmad
| left-exchanging magmads
| Left-exchanging magmads
| left-exchanging monad
| Left-exchanging monad
| left-exchanging monads
| Left-exchanging monads
| left-exchange
| Left-exchange

\knowledge{notion}
| affine magmoid
| Affine magmoid
| affine magmoids
| Affine magmoids

\knowledge{notion}
| affine
| Affine
| affine monad
| Affine monad

\knowledge{notion}
| relevant
| Relevant
| relevant monad
| Relevant monad

\knowledge{notion}
| copy-discard sesquilaw
| Copy-discard sesquilaw
| copy-discard sesquilaws
| Copy-discard sesquilaws

\knowledge{notion}
| sesquilaw magmoid
| Sesquilaw magmoid
| sesquilaw magmoids
| Sesquilaw magmoids

\knowledge{notion}
| correctness triple
| Correctness triple
| correctness triples
| Correctness triples

\knowledge{notion}
| posetal imperative category
| Posetal imperative category
| posetal imperative categories
| Posetal imperative categories

\knowledge{notion}
| posetal imperative multicategory
| Posetal imperative multicategory
| posetal imperative multicategories
| Posetal imperative multicategories

\knowledge{notion}
| full support
| Full support

\knowledge{notion}
| hypergraph
| Hypergraph
| hypergraph magmoid
| hypergraph magmoids
| Hypergraph magmoid
| Hypergraph magmoids

\knowledge{notion}
| guard combinators
| Guard combinators
| guard combinator
| Guard combinator

\knowledge{notion}
| command combinators
| Command combinators
| command combinator
| Command combinator

\knowledge{notion}
| monoidal category
| monoidal categories
| Monoidal category
| Monoidal categories

\knowledge{notion}
| sesquifunctor
| Sesquifunctor
| sesquifunctors
| Sesquifunctors

\knowledge{notion}
| premonoidal category
| Premonoidal category
| premonoidal categories
| Premonoidal categories
| premonoidality
| Premonoidality
| premonoidal
| Premonoidal

\knowledge{notion}
| central
| Central
| central morphism
| Central morphism
| central morphisms
| Central morphisms

\knowledge{notion}
| Symmetric premonoidal category
| Symmetric premonoidal categories
| symmetric premonoidal category
| symmetric premonoidal categories

\knowledge{notion}
| deterministic
| Deterministic
| deterministic morphism
| Deterministic morphism
| deterministic morphisms
| Deterministic morphisms

\knowledge{notion}
| total
| total morphism
| total morphisms
| Total
| Total morphism
| Total morphisms

\knowledge{notion}
| cocartesian multicategory
| Cocartesian multicategory
| cocartesian multicategories
| Cocartesian multicategories

\knowledge{notion}
| predistributive multicategory
| Predistributive multicategory
| predistributive multicategories
| Predistributive multicategories

\knowledge{notion}
| posetal distributive copy-discard multicategory
| Posetal distributive copy-discard multicategory
| posetal distributive copy-discard multicategories
| Posetal distributive copy-discard multicategories

\knowledge{notion}
| index
| indices
| Index
| Indices

\knowledge{notion}
| state
| states
| State
| States

\knowledge{notion}
| state combinator
| state combinators
| State combinator
| State combinators

\knowledge{notion}
| multicategory
| Multicategory
| Multicategories
| multicategories

\knowledge{notion}
| distributive multicategory
| distributive multicategories
| Distributive multicategory
| Distributive multicategories

\knowledge{notion}
| substitution
| Substitution
| substitutions
| Substitutions
| substituting
| Substituting

\knowledge{notion}
| guard
| guards
| Guard
| Guards

\knowledge{notion}
| predicate
| predicates
| Predicates
| Predicate

\knowledge{notion}
| predicate combinators
| Predicate combinators

\knowledge{notion}
| command
| commands
| Commands
| Command

\knowledge{notion}
| fresh
| Fresh

\knowledge{notion}
| alpha-equivalent
| alpha-equivalence
| Alpha-equivalent
| Alpha-equivalence

\knowledge{notion}
| variable
| variables
| Variable
| Variables

\knowledge{notion}
| traced distributive copy-discard multicategory
| Traced distributive copy-discard multicategory
| traced distributive copy-discard multicategories
| Traced distributive copy-discard multicategories

\knowledge{notion}
| traced distributive multicategory
| Traced distributive multicategory
| traced distributive multicategories
| Traced distributive multicategories

\knowledge{notion}
| traced predistributive copy-discard multicategory
| Traced predistributive copy-discard multicategory
| traced predistributive copy-discard multicategories
| Traced predistributive copy-discard multicategories

\knowledge{notion}
| predistributive copy-discard multicategory
| Predistributive copy-discard multicategory
| predistributive copy-discard multicategories
| Predistributive copy-discard multicategories

\knowledge{notion}
| distributive copy-discard multicategory
| Distributive copy-discard multicategory
| distributive copy-discard multicategories
| Distributive copy-discard multicategories

\knowledge{notion}
| variable substitution
| variable substitutions
| Variable substitution
| Variable substitutions

\knowledge{notion}
| posetal distributive signature
| posetal distributive signatures
| Posetal distributive signature
| Posetal distributive signatures

\knowledge{notion}
| posetal uniformity
| Posetal uniformity

\knowledge{notion}
| posetal distributive copy-discard category
| posetal distributive copy-discard categories
| Posetal distributive copy-discard category
| Posetal distributive copy-discard categories

\knowledge{notion}
| posetal uniform trace
| posetal uniform traces
| Posetal uniform traces
| Posetal uniform trace
| posetal uniform traced monoidal category
| Posetal uniform traced monoidal category
| posetal uniform traced monoidal categories
| Posetal uniform traced monoidal categories

\knowledge{notion}
| label substitution
| label substitutions
| Label substitution
| Label substitutions

\knowledge{notion}
| commutative Lawvere theory
| Commutative Lawvere theory
| commutative Lawvere theories
| Commutative Lawvere theories

\knowledge{notion}
| Lawvere theory
| Lawvere theories
| clone
| Clone
| clones
| Clones

\knowledge{notion}
| anchor
| anchors
| Anchor
| Anchors
| label
| Label
| Labels
| labels

\knowledge{notion}
| conditional
| conditionals
| Conditional
| Conditionals

\knowledge{notion}
| conditional composition
| conditional compositions
| Conditional composition
| Conditional compositions

\knowledge{notion}
| context
| Context
| contexts
| Contexts

\knowledge{notion}
| copy-discard category
| Copy-discard category
| copy-discard categories
| Copy-discard categories
| copy-discard
| Copy-discard

\knowledge{notion}
| copy-discard monoidal category
| Copy-discard monoidal category
| copy-discard monoidal categories
| Copy-discard monoidal categories

\knowledge{notion}
| copy-discard premonoidal category
| Copy-discard premonoidal category
| copy-discard premonoidal categories
| Copy-discard premonoidal categories
| premonoidal copy-discard category
| Premonoidal copy-discard category
| premonoidal copy-discard categories
| Premonoidal copy-discard categories

\knowledge{notion}
| commutative imperative category
| Commutative imperative category
| commutative imperative categories
| Commutative imperative categories

\knowledge{notion}
| imperative category
| Imperative category
| imperative categories
| Imperative categories

\knowledge{notion}
| imperative multicategory
| Imperative multicategory
| imperative multicategories
| Imperative multicategories

\knowledge{notion}
| multiquiver
| multiquivers
| Multiquiver
| Multiquivers

\knowledge{notion}
| tricocycloid
| tricocycloids
| Tricocycloid
| Tricocycloids
| symmetric tricocycloid
| symmetric tricocycloids
| Symmetric tricocycloid
| Symmetric tricocycloids

\knowledge{notion}
| right do-notation
| Right do-notation
| right do notation
| Right do notation

\knowledge{notion}
| normalized distribution magmad
| Normalized distribution magmad
| magmad of normalized distributions
| Magmad of normalized distributions

\knowledge{notion}
| normalized distributions magmoid
| normalized distribution magmoid
| normalization magmoid
| Normalized distributions magmoid
| Normalized distribution magmoid
| Normalization magmoid
| monoidal magmoid of normalized finitary distributions
| normalized finitary distribution magmoid
\newcommand{\Norm}{\kl[Normalization magmoid]{\mathbf{Norm}}}

\knowledge{notion}
| finite normalized kernels
| finite normalized kernel
| Finite normalized kernels
| Finite normalized kernel

\knowledge{notion}
| left-relevant magmoid
| left-relevant magmoids
| Left-relevant magmoid
| Left-relevant magmoids

\knowledge{notion}
| quasitotal magmoid
| quasitotal magmoids
| Quasitotal magmoid
|  magmoids

\knowledge{notion}
| monoidal magmoid
| monoidal magmoids
| Monoidal magmoid
| Monoidal magmoids

\knowledge{notion}
| symmetric monoidal magmoid
| symmetric monoidal magmoids
| Symmetric monoidal magmoid
| Symmetric monoidal magmoids

\knowledge{notion}
| symmetric monoidal non-associative category
| Symmetric monoidal non-associative category
| symmetric monoidal non-associative categories
| Symmetric monoidal non-associative categories

\knowledge{notion}
| tensor schema
| tensor schemas
| Tensor schema
| Tensor schemas

\knowledge{notion}
| string diagram
| string diagrams
| String diagram
| String diagrams

\knowledge{notion}
| commuting magmoid
| commuting magmoids
| Commuting magmoid
| Commuting magmoids
| commutative magmoid
| commutative magmoids
| Commutative magmoid
| Commutative magmoids
| left commutative magmoid
| left commutative magmoids
| Left commutative magmoid
| Left commutative magmoids
\knowledge{notion}
| right commutative magmoid
| right commutative magmoids
| Right commutative magmoid
| Right commutative magmoids

\knowledge{notion}
| commutative non-associative monad
| commutative non-associative monads
| Commutative non-associative monad
| Commutative non-associative monads
| commuting non-associative monad
| commuting non-associative monads
| Commuting non-associative monad
| Commuting non-associative monads

\knowledge{notion}
| strict monoidal magmoid
| strict monoidal magmoids
| Strict monoidal magmoid
| Strict monoidal magmoids

\knowledge{notion}
| bimonoidally strict
| Bimonoidally strict
| strict
| Strict

\knowledge{notion}
| distributive signature
| distributive signatures
| Distributive signature
| Distributive signatures

\knowledge{notion}
| distributive category
| distributive categories
| Distributive category
| Distributive categories
| distributive copy-discard category
| distributive copy-discard categories
| Distributive copy-discard category
| Distributive copy-discard categories

\knowledge{notion}
| basic type
| basic types
| Basic type
| Basic types

\knowledge{notion}
| generator
| generators
| Generator
| Generators

\knowledge{notion}
| effectful triple
| effectful triples
| Effectful triple
| Effectful triples

\knowledge{notion}
| double signature
| double signatures
| Double signature
| Double signatures

\knowledge{notion}
| forgetful double signature
| forgetful double signatures
| Forgetful double signature
| Forgetful double signatures

\knowledge{notion}
| double signature morphism
| double signature morphisms
| Double signature morphism
| Double signature morphisms

\knowledge{notion}
| pinwheel double category
| pinwheel double categories
| Pinwheel double category
| Pinwheel double categories
| double category with pinwheels
| double categories with pinwheels
| Double category with pinwheels
| Double categories with pinwheels

\knowledge{notion}
| bigraph
| Bigraph
| bigraphs
| Bigraphs
| 2-graph
| 2-graphs
| 2-Graph
| 2-Graphs

\knowledge{notion}
| 2-graph morphism
| 2-graph morphisms

\knowledge{notion}
| 2-graph category
| 2-graph categories

\knowledge{notion}
| double category
| double categories
| Double category
| Double categories

\knowledge{notion}
| weighted polygraph
| weighted polygraphs
| Weighted polygraph
| Weighted polygraphs
| timed polygraph
| timed polygraphs
| Timed polygraph
| Timed polygraphs

\knowledge{notion}
| morphism of timed polygraphs
| Morphism of timed polygraphs
| timed polygraph morphism
| Timed polygraph morphism

\knowledge{notion}
| tilted bicategory signature
| tilted bicategory signatures
| Tilted bicategory signature
| Tilted bicategory signatures
| tilted bigraph
| Tilted bigraph
| tilted bigraphs
| Tilted bigraphs
| tilted 2-graph
| Tilted 2-graph
| tilted 2-graphs
| Tilted 2-graphs

\knowledge{notion}
| path
| paths
| Path
| Paths
| composable path
| composable paths

\knowledge{notion}
| signature
| signatures
| Signature
| Signatures

\knowledge{notion}
| congruence
| congruences
| Congruence
| Congruences

\knowledge{notion}
| symmetry
| symmetries
| Symmetry
| Symmetries

\knowledge{notion}
| profunctor
| profunctors
| Profunctor
| Profunctors

\knowledge{notion}
| left inclusion
| Left inclusion

\knowledge{notion}
| right inclusion
| Right inclusion

\knowledge{notion}
| left whiskering
| Left whiskering

\knowledge{notion}
| right whiskering
| Right whiskering

\knowledge{notion}
| inclusion and whiskering
| Inclusion and whiskering
| whiskering interacts with inclusion
| Whiskering interacts with inclusion

\knowledge{notion}
| whiskering
| Whiskering

\knowledge{notion}
| rewiring preserves whiskering
| Rewiring preserves whiskering
| whiskering preserves rewiring
| Whiskering preserves rewiring

\knowledge{notion}
| rewiring preserves interchange
| Rewiring preserves interchange

\knowledge{notion}
| rewiring preserves composition
| Rewiring preserves composition

\knowledge{notion}
| rewiring
| Rewiring 
| definition of rewiring

\knowledge{rewiring is an action}[]{notion}

\knowledge{notion}
| composition
| Composition
| term composition
| Term composition

\knowledge{notion}
| term composition is associative
| Term composition is associative

\knowledge{notion}
| term composition is unital
| Term composition is unital

\knowledge{notion}
| tensoring
| Tensoring

\knowledge{notion}
| arrow notation preterm
| Arrow notation preterm
| arrow notation preterms
| Arrow notation preterms

\knowledge{notion}
| arrow notation term
| Arrow notation term
| arrow notation terms
| Arrow notation terms

\knowledge{notion}
| almost monad
| almost monads
| Almost monad
| Almost monads

\knowledge{notion}
| category
| categories
| Category
| Categories

\knowledge{notion}
| polyquiver
| polyquivers
| Polyquiver
| Polyquivers

\knowledge{notion}
| monoidal almost distributive law
| Monoidal almost distributive law
| monoidal almost distributive laws
| Monoidal almost distributive laws

\knowledge{notion}
| symmetric monoidal magmoid
| Symmetric monoidal magmoid
| symmetric monoidal magmoids
| Symmetric monoidal magmoids

\knowledge{notion}
| observe-arrow notation
| Observe-arrow notation
| observe-arrow notation term
| Observe-arrow notation term
| observe-arrow notation terms
| Observe-arrow notation terms

\knowledge{notion}
| subdistribution
| Subdistribution
| subdistributions
| Subdistributions
| subdistributional
| substochastic kernel
| substochastic kernels
| Substochastic kernel
| Substochastic kernels

\knowledge{notion}
| rescaling
| Rescaling

\knowledge{notion}
| normalisation
| normalisations
| Normalisation
| Normalisations
| normalization
| normalizations
| Normalization
| Normalizations
| normalised
| Normalised
\newcommand{\norm}[1]{{\kl[normalization]{\mathsf{n}}}{#1}}

\knowledge{notion}
| continuous normalization
| Continuous normalization
\newcommand{\cnorm}[1]{{\kl[continuous normalization]{\mathsf{N}}}{#1}}

\knowledge{notion}
| observe-arrow notation copy-discard-compare category of terms

\knowledge{notion}
| category of arrow notation terms

\knowledge{notion}
| copy-discard-compare categories
| copy-discard-compare category
| Copy-discard-compare categories
| Copy-discard-compare category

\knowledge{notion}
| Kleisli extension
| Kleisli extensions

\knowledge{notion}
| Kleisli category
| Kleisli categories

\knowledge{notion}
| partially additive
| partially additive monad
| Partially additive
| Partially additive monad

\knowledge{notion}
| sets
| category of sets
\renewcommand{\Set}{\kl[sets]{\mathsf{Set}}}

\knowledge{notion}
| powerset
| powerset monad

\knowledge{notion}
| category of relations
| Category of relations

\knowledge{notion}
| maybe
| maybe monad
| Maybe
| Maybe monad
\DeclareUnicodeCharacter{1D40C}{\kl[Maybe monad]{\mathbf{M}}}

\knowledge{notion}
| continuous maybe
| continuous maybe monad
| Continuous maybe
| Continuous maybe monad
\DeclareUnicodeCharacter{1D4DC}{\kl[continuous maybe monad]{\mathcal{M}}}

\knowledge{notion}
| partial stochastic kernels
| Partial stochastic kernels
| partial stochastic kernel
| Partial stochastic kernel

\knowledge{notion}
| category of partial functions
| Category of partial functions

\knowledge{notion}
| Markov category
| Markov categories

\knowledge{notion}
| Markov magmoid
| Markov magmoids

\knowledge{notion}
| discrete Markov magmoid
| discrete Markov magmoids
| Discrete Markov magmoid
| Discrete Markov magmoids

\knowledge{notion}
| exact observation
| Exact observation
| exact observations
| Exact observations

\knowledge{notion}
| discrete copy-discard magmoid
| discrete copy-discard magmoids
| Discrete copy-discard magmoid
| Discrete copy-discard magmoids
| Frobenius equation

\knowledge{notion}
| quasitotal
| Quasitotal
| quasitotality
| Quasitotality
| quasitotal morphism
| Quasitotal morphism
| quasitotal morphisms
| Quasitotal morphisms

\knowledge{notion}
| partial Markov category
| partial Markov categories
| Partial Markov category
| Partial Markov categories

\knowledge{notion}
| quasi-Markov category
| quasi-Markov categories
| Quasi-Markov category
| Quasi-Markov categories

\knowledge{notion}
| monoidal almost-distributive law
| monoidal almost-distributive laws
| Monoidal almost-distributive law
| Monoidal almost-distributive laws

\knowledge{notion}
| black-hole semantics
| Black-hole semantics
| black-hole
| Black-hole

\knowledge{notion}
| almost-distributive law
| almost-distributive laws
| Almost-distributive law
| Almost-distributive laws

\knowledge{notion}
| magmoid
| Magmoid
| magmoids
| Magmoids
| Unital magmoid
| unital magmoid
| Unital magmoids
| unital magmoids
\knowledge{notion}
| non-associative category
| Non-associative category
| non-associative categories
\knowledge{notion}
| monoidal non-associative category
| Monoidal non-associative category
| monoidal non-associative categories
| Monoidal non-associative categories

\knowledge{notion}
| Giry monad
\newcommand{\unicodecalD}{\ensuremath{\kl[Giry monad]{\mathcal{D}}}}
\DeclareUnicodeCharacter{1D4D3}{\unicodecalD}

\knowledge{notion}
| finitary distribution monad
| distribution monad
| Finitary distribution monad
| Distribution monad
\DeclareUnicodeCharacter{1D403}{\kl[distribution monad]{\mathbf{D}}} %

\knowledge{notion}
| discrete subdistribution
| discrete subdistributions
| discrete subdistribution monad

\knowledge{notion}
| category of discrete stochastic channels
| Category of discrete stochastic channels

\knowledge{notion}
| support
| Support
\newcommand{\supp}{\kl[support]{\mathrm{supp}}}

\knowledge{notion}
| Dirac distribution

\knowledge{notion}
| measurable subdistribution
| measurable subdistributions
| measurable subdistribution monad

\knowledge{notion}
| standard Borel spaces
| Standard Borel spaces

\knowledge{notion}
| quantale-valued sets

\knowledge{notion}
| assertion-correctness triple
| assertion-correctness triples
| Assertion-correctness triple
| Assertion-correctness triples

\knowledge{notion}
| state-correctness triple
| state-correctness triples
| State-correctness triple
| State-correctness triples

\knowledge{notion}
| predicate-correctness triple
| predicate-correctness triples
| Predicate-correctness triple
| Predicate-correctness triples

\knowledge{notion}
| assertion-incorrectness triple
| assertion-incorrectness triples
| Assertion-incorrectness triple
| Assertion-incorrectness triples

\knowledge{notion}
| state-incorrectness triple
| state-incorrectness triples
| State-incorrectness triple
| State-incorrectness triples

\knowledge{notion}
| predicate-incorrectness triple
| predicate-incorrectness triples
| Predicate-incorrectness triple
| Predicate-incorrectness triples

\knowledge{notion}
| relational assertion-correctness triple
| relational assertion-correctness triples
| Relational assertion-correctness triple
| Relational assertion-correctness triples

\knowledge{notion}
| relational state-correctness triple
| relational state-correctness triples
| Relational state-correctness triple
| Relational state-correctness triples

\knowledge{notion}
| relational predicate-correctness triple
| relational predicate-correctness triples
| Relational predicate-correctness triple
| Relational predicate-correctness triples

\knowledge{notion}
| relational assertion-incorrectness triple
| relational assertion-incorrectness triples
| Relational assertion-incorrectness triple
| Relational assertion-incorrectness triples

\knowledge{notion}
| relational state-incorrectness triple
| relational state-incorrectness triples
| Relational state-incorrectness triple
| Relational state-incorrectness triples

\knowledge{notion}
| relational predicate-incorrectness triple
| relational predicate-incorrectness triples
| Relational predicate-incorrectness triple
| Relational predicate-incorrectness triples

\knowledge{notion}
| extension
| extensions
| Extension
| Extensions

\knowledge{notion}
| non-associative monad
| non-associative monads
| Non-associative monad
| Non-associative monads

\knowledge{notion}
| copy-discard magmoid
| copy-discard magmoids
| copy discard magmoid
| copy discard magmoids
| Copy-discard magmoid
| Copy-discard magmoids
| Copy discard magmoid
| Copy discard magmoids

\theoremstyle{plain}
\newtheoremrep{theorem}{Theorem}[section]
\newtheorem*{theorem*}{Theorem}
\newtheoremrep{proposition}[theorem]{Proposition}
\newtheoremrep{lemma}[theorem]{Lemma}
\newtheoremrep{corollary}[theorem]{Corollary}

\theoremstyle{definition}
\newtheorem{definition}[theorem]{Definition}

\theoremstyle{remark}
\newtheorem{remark}[theorem]{Remark}
\newtheorem{example}[theorem]{Example}
\usepackage{xparse}

\NewDocumentCommand{\multicite}{ogogogog}{%
    \citetext{%
        \IfValueT{#2}{%
            \IfValueT{#1}{\citealp[#1]{#2}}%
            \IfNoValueT{#1}{\citealp{#2}}%
        }%
        \IfValueT{#4}{%
            ;
            \IfValueT{#3}{\citealp[#3]{#4}}%
            \IfNoValueT{#3}{\citealp{#4}}%
        }%
        \IfValueT{#6}{%
            ;
            \IfValueT{#5}{\citealp[#5]{#6}}%
            \IfNoValueT{#5}{\citealp{#6}}%
        }%
        \IfValueT{#8}{%
            ;
            \IfValueT{#7}{\citealp[#7]{#8}}%
            \IfNoValueT{#7}{\citealp{#8}}%
        }%
    }%
} %
\definecolor{nord0}{RGB}{46, 52, 64}
\definecolor{nord1}{RGB}{59, 66, 82}
\definecolor{nord2}{RGB}{67, 76, 94}
\definecolor{nord3}{RGB}{76, 86, 106}
\definecolor{nord4}{RGB}{216, 222, 233}
\definecolor{nord5}{RGB}{229, 233, 240}
\definecolor{nord6}{RGB}{236, 239, 244}
\definecolor{nord7}{RGB}{143, 188, 187}
\definecolor{nord8}{RGB}{136, 192, 208}
\definecolor{nord9}{RGB}{129, 161, 193}
\definecolor{nord10}{RGB}{94, 129, 172}
\definecolor{nord11}{RGB}{191, 97, 106}
\definecolor{nord12}{RGB}{208, 135, 112}
\definecolor{nord13}{RGB}{235, 203, 139}
\definecolor{nord14}{RGB}{163, 190, 140}
\definecolor{nord15}{RGB}{180, 142, 173}

\usepackage{listings}

\lstdefinelanguage{Racket}{
  sensitive = true,
  alsoletter = {<,>,-,'},
  keywords={if, equal?, return, rDo, do, lDo, accDo, leftDo, rightDo, Norm, frontDoor, backDoor, <-},
  otherkeywords={>, ==},
  keywords = [2]{match, @get, @post, @put, @delete, require, define, ensure, otherwise, call, lang, racket, define-syntax, syntax-rules},
  keywords = [3]{'left, 'middle, 'right, '(), uniform, uniformDoor, host, observeLeft, prevalence, test, uncertainty, observe, 'gene, 'nogene, 'smoker, 'nonsmoker, 'tar, 'notar, list, <-},
  keywordstyle={\bfseries\color{nord3}},%
  keywordstyle=[2]\color{nord12},%
  keywordstyle=[3]{\color{nord3}},%
  keywordstyle=[4]\color{nord15},
  numbers=left,
  basicstyle={\small\ttfamily\color{nord0}},
  columns=flexible,
  keepspaces=true,
  numberstyle={\tiny\ttfamily},
  stepnumber=1,
  numbersep=6pt,
  showstringspaces=false,
  breaklines=true,
  frame=single,
  frameround=tttt,
  comment=[l]{;},
  morecomment=[s]{;;*}{*/},
  commentstyle={\color{nord14}\ttfamily},
  stringstyle={\color{nord10}\ttfamily},
  morestring=[b]",
  literate={†}{$\smash{{}^{\dagger}}$}1 {λ}{$\lambda$}1 {η}{$\eta$}1
}
\newcommand*\dif{\mathop{}\!\mathrm{d}} 
\AddToHook{env/proposition/begin}{\crefalias{theorem}{proposition}}
\AddToHook{env/lemma/begin}{\crefalias{theorem}{lemma}}
\AddToHook{env/corollary/begin}{\crefalias{theorem}{corollary}}
\AddToHook{env/definition/begin}{\crefalias{theorem}{definition}}
\AddToHook{env/remark/begin}{\crefalias{theorem}{remark}}

\newcommand{\ketb}[1]{\ket{\cBlu{#1}}}

\newcommand{\er}{\mathscr{r}}

\AtBeginDocument{}

\DeclareRobustCommand{\NOOP}[1]{}

\newcommand{\blackSpider}{\smash{\raisebox{-0.25ex}{\includegraphics[height=0.8em]{fig-mini-black-spider.pdf}}}}
\newcommand{\whiteSpider}{\smash{\raisebox{-0.25ex}{\includegraphics[height=0.8em]{fig-mini-white-spider.pdf}}}}
\newcommand{\miniFrob}{\smash{\raisebox{-0.25ex}{\includegraphics[height=0.9em]{fig-mini-frobenius.pdf}}}}
\newcommand{\miniComm}{\smash{\raisebox{-0.25ex}{\includegraphics[height=0.9em]{fig-mini-commuting.pdf}}}}
\newcommand{\miniStrings}[1]{{{\hyperref[#1]{{{\miniFrob}}}}}}
\newcommand{\miniCommute}[1]{{{\hyperref[#1]{{{\miniComm}}}}}}

\let\section\savedsection
\makeatletter
\let\ACM@origbaselinestretch\baselinestretch
\makeatother
  
\allowdisplaybreaks

\title[Normalized Probabilistic Semantics is Not Associative]%
{Normalized Probabilistic Semantics is Not Associative}

\author{Elena Di Lavore}
\affiliation{
  \institution{Tallinn University of Technology}
  \country{Estonia}
}
 
\author{Mario Rom\'an}
\affiliation{
  \institution{Tallinn University of Technology}
  \country{Estonia}}

\author{M\'ark Sz\'eles}
\affiliation{
  \institution{Radboud University Nijmegen}
  \country{The Netherlands}}
\date{\today}

\usepackage{etoolbox}
\makeatletter
\patchcmd{\ACM@mk@linecount}{\color{red}}{\color{gray}}{}{}
\patchcmd{\ACM@mk@linecount}{\color{red}}{\color{gray}}{}{}
\makeatother

\begin{document}

\begin{abstract}
  Normalization, $\mathbf{D}(X + 1) → \mathbf{D}(X) + 1$, fails to form a distributive law but still yields a non-associative monad—i.e., a \emph{magmad}. We introduce a first normalized-by-construction probabilistic semantics: the magmad of normalized distributions. Non-associativity allows this semantics to express both observations and interventions, both Pearl's and Jeffrey's updates, and both evidential and causal interpretations. Within this framework, we derive results from causal inference, including Pearl's front-door and back-door criteria.
\end{abstract} %
\maketitle 
\medskip

\section{Introduction}
\label{sec:introduction}
\newcommand{\cst}[1]{\mbox{\sffamily\color{nord3}{#1}}}

\kl{Normalization} is essential in probabilistic inference: Bayes' rule for updating a prior with a likelihood requires rescaling to compute a valid posterior,
$$\mathit{posterior} \propto \mathit{likelihood} · \mathit{prior}.$$ 
However, probabilistic programming semantics is not necessarily \emph{normalized-by-construction}: \kl{normalization} typically appears as a primitive and, when updating, denotational models rely on subdistributions~\cite{panangaden1999category,borgstrom16} or unnormalized distributions~\cite{kozen81,staton2017commutative,culpepper18,ehrhard2017measurable}.

\medskip

We introduce the first denotational normalized-by-construction probabilistic semantics: the non-associative \kl{monad}—or, \emph{\kl{magmad}}—of \kl{normalized distributions}. 
Perhaps surprisingly, non-associativity is forced and, far from a limitation, it grants expressive power.  Let us illustrate this semantics and its non-associativity with two examples of progressive complexity.

\subsection{Example—The Monty Hall problem}
Imagine solving the famous \emph{Monty Hall problem} \cite{selvin75:montyhall}: 
in a game show, you choose from three closed doors for a chance of winning the prize behind one of them; however, after choosing, the host opens one of the empty doors and—with two closed doors standing—invites you to switch your guess. \emph{Should you switch?}

\newcommand{\scancel}[1]{\smash{\cancel{\texttransparent{0.5}{#1}}}}
\newcommand{\dprize}[1]{\ensuremath{\cBlu{#1}}}
\newcommand{\dhost}[1]{\ensuremath{\cRed{#1}}}
\medskip

{
\newcommand{\pL}{\ket{\dprize{L}}}
\newcommand{\pM}{\ket{\dprize{M}}}
\newcommand{\pR}{\ket{\dprize{R}}}
\newcommand{\dL}{\ket{\dhost{L}}}
\newcommand{\dM}{\ket{\dhost{M}}}
\newcommand{\dR}{\ket{\dhost{R}}}
\newcommand{\tdL}{\texttransparent{0.5}{\ket{\dhost{L}}}}
Let us compute.\footnote{We use \emph{ket-notation} to denote formal sums. The generic distribution over a set $X$ is written as $λ_1\ketb{x_1} + \dots + λ_n\ketb{x_n}$ for some elements $\cBlu{x_1}, \dots, \cBlu{x_n} ∈ X$ and some scalars $λ_1, \dots,λ_n ∈ ℝ^{+}$ adding up to the unit, $λ_1 + \dots + λ_n = 1$.}
We \emph{(i)} consider a prior uniform probability that a prize is behind any of three doors (\dprize{L}, \dprize{M}, \dprize{R}); then, \emph{(ii)} after choosing, say, the middle door, the host will randomly and uniformly open a door (\dhost{L}, \dhost{M}, \dhost{R}) that cannot be neither the chosen door nor the one with the prize; \emph{(iii)} we then observe that the host opens, e.g., the left door (\dhost{L}) precluding any other possibility; and, upon this, \emph{(n)} we renormalize the remaining probabilities and \emph{(r)} conclude we should switch to the right door (\dprize{R}): it doubles our chances of getting the prize—the posterior distribution is $\tfrac{1}{3} \pM + \tfrac{2}{3} \pR$.
\begin{alignat*}{4}
 & \mbox{\emph{(i)}}\qquad && 
   \tfrac{1}{3} \pL + 
   \tfrac{1}{3} \pM +
   \tfrac{1}{3} \pR 
 \\[-0.2em]
 & \mbox{\emph{(ii)}} && 
     \tfrac{1}{3} \pL\!\dR + 
     \tfrac{1}{6} \pM\!\dL + 
     \tfrac{1}{6} \pM\!\dR + 
     \tfrac{1}{3} \pR\!\dL 
   \\[-0.2em]
   &\mbox{\emph{(iii)}} && 
     \scancel{\tfrac{1}{3} \pL\!\dR} + 
     \tfrac{1}{6} \pM\!\dL +
     \scancel{\tfrac{1}{6} \pM\!\dR} + 
     \tfrac{1}{3} \pR\!\dL
   \\[-0.2em]
   &\mbox{\emph{(n)}}  && 
     \tfrac{1}{3} \pM\!\dL + 
     \tfrac{2}{3} \pR\!\dL
  \\[-0.2em]
   &\mbox{\emph{(r)}}  && 
     \tfrac{1}{3} \pM + 
     \tfrac{2}{3} \pR.
\end{alignat*}

\noindent
However, the \emph{Monty Hall problem} admits another parenthesization. Let us consider everything since the action of the host a parenthesized subproblem: we now \emph{(n$_2$)} normalize internally before \emph{(n$_1$)} normalizing globally, and the posterior is $\tfrac{1}{2}\pM + \tfrac{1}{2}\pR$—switching door stops mattering.
\begin{alignat*}{4}
   & \emph{(i)}\quad 
   && \tfrac{1}{3} \pL + \tfrac{1}{3} \pM + \tfrac{1}{3} \pR 
   \\[-0.2em]
   &\quad \emph{(ii)} 
   && \quad 
    \tfrac{1}{3} \pL\![\tfrac{1}{1} \ket{\dhost{R}}] + 
    \tfrac{1}{3} \pM\![\tfrac{1}{2} \ket{\dhost{L}} + 
                      \tfrac{1}{2} \ket{\dhost{R}}] + 
    \tfrac{1}{3} \pR\![\tfrac{1}{1} \ket{\dhost{L}}]
   \\[-0.2em]
   &\quad \emph{(iii)}
   && \quad 
   \tfrac{1}{3} \pL\![\scancel{\tfrac{1}{1} \ket{\dhost{R}}}] + 
   \tfrac{1}{3} \pM\![\tfrac{1}{2}\ket{\dhost{L}} + 
                     \scancel{\tfrac{1}{2}\ket{\dhost{R}}}] + 
   \tfrac{1}{3} \pR\![\tfrac{1}{1} \dL]
   \\[-0.2em]
   &  \quad \emph{(n$_2$)}
   && \quad 
   \scancel{\tfrac{1}{3} \pL\![0]} + 
   \tfrac{1}{3} \pM\![\tfrac{1}{1} \dL] + 
   \tfrac{1}{3} \pR\![\tfrac{1}{1} \dL]
   \\[-0.2em]
   &  \emph{(n$_1$)} 
   && \tfrac{1}{2}\pM\!\dL + \tfrac{1}{2}\pR\!\dL
   \\[-0.2em]
   &  \emph{(r)} 
   && \tfrac{1}{2}\pM + \tfrac{1}{2}\pR.
\end{alignat*}
An interpretation for this case is that \emph{(iii)} we \emph{intervene} to force the host to open the left door $(\dhost{L})$, and the game show halts otherwise. Because it is forced, the host's decision stops carrying any inferential information: we are equally likely to see it no matter where the prize is.

\newsavebox{\codeMontyHallA}
\begin{lrbox}{\codeMontyHallA}
\begin{minipage}{0.4\textwidth}
\begin{lstlisting}[language=Racket,numbers=none,mathescape]
(do (prize) <- uniformDoor
    (door) <- (host prize)
    () <- (observeLeft door)
    return (prize))
\end{lstlisting}
\end{minipage}
\end{lrbox}
\newsavebox{\codeMontyHallB}
\begin{lrbox}{\codeMontyHallB}
\begin{minipage}{0.4\textwidth}
\begin{lstlisting}[language=Racket,numbers=none,mathescape]
(do (prize) <- uniformDoor
    (door) <- (do
       (door) <- (host prize)
       () <- (observeLeft door)
       return (door))
    return (door))
\end{lstlisting}
\end{minipage}
\end{lrbox}
Monty Hall is a famously controversial problem \cite{selvin75:montyhall,vosSavant}. We will boil down controversy to a difference of parenthesization: \emph{when to normalize} does matter and—when parentheses normalize—parentheses matter. Monty Hall is a %
counterexample to associativity. %
Indeed, our two computations mirror two programs that should be equal for any associative semantics (\Cref{fig:montyhall}).
\begin{figure}[!h]
  \centering
\[\usebox{\codeMontyHallA}\quad ≠\quad \usebox{\codeMontyHallB}
\]\vspace{-1em}
\caption{Two interpretations of the Monty Hall problem.\label{fig:montyhall}}
\end{figure}
\smallskip

\newcommand{\uniformDoor}{\cst{uniformDoor}}
\newcommand{\host}{\cst{host}}
\newcommand{\observeLeft}{\cst{observeLeft}}
These programs match our computations when \emph{(i)} we define  $\uniformDoor = \tfrac{1}{3}\ketb{L} + \tfrac{1}{3}\ketb{M} + \tfrac{1}{3}\ketb{R}$;  \emph{(ii)} we define $\host(\cBlu{L}) = \dR$, with $\host(\cBlu{M}) = \tfrac{1}{2}\dR + \tfrac{1}{2}\dL$ and $\host(\cBlu{R}) = \dL$; and  \emph{(iii)} we define $\observeLeft(\cBlu{L}) = 1\ket{}$ and $\observeLeft(\cBlu{M}) = \observeLeft(\cBlu{R}) = 0$.
}

\newcommand{\dprev}[1]{\ensuremath{\cBlu{#1}}}
\newcommand{\dtest}[1]{\ensuremath{\cRed{#1}}}

\subsection{Example—The Unclear Test problem\label{sec:unclear-test}}
{
\newcommand{\sfrac}[2]{\tfrac{#1}{#2}}
\newcommand{\I}{\ket{\dprev{I}}}
\renewcommand{\H}{\ket{\dprev{H}}}
\renewcommand{\P}{\ket{\dtest{P}}}
\newcommand{\N}{\ket{\dtest{N}}}
\emph{Does this transfer to numerical disagreements?} Let us explicitly solve a standard \emph{testing problem}. Consider an illness with a $33\%$ prevalence ($\cst{prevalence} = \sfrac{1}{3}\I + \sfrac{2}{3} \H$) for which  we have a test with $75\%$ sensitivity ($\cst{test}(\cBlu{I}) = \sfrac{3}{4} \P + \sfrac{1}{4} \N$) and $50\%$ specificity ($\cst{test}(\cBlu{H}) = \sfrac{1}{2} \P + \sfrac{1}{2} \N$). What is the posterior probability of having this illness after a positive test?
\begin{alignat*}{4}
   &\emph{(i)}\quad 
   && \tfrac{1}{3} \I + \tfrac{1}{3} \H
   \\[-0.2em]
   &\emph{(ii)}\quad  
   && \tfrac{1}{4} \I\!\P + \tfrac{1}{12} \I\!\N 
    + \tfrac{1}{3} \H\!\P + \tfrac{1}{3} \H\!\N
   \\[-0.2em]
   &\emph{(iii)}\quad
   && \tfrac{1}{4} \I\!\P + \scancel{\tfrac{1}{12} \I\!\N}
    + \tfrac{1}{3} \H\!\P + \scancel{\tfrac{1}{3} \H\!\N}
   \\[-0.2em]
   &\emph{(n)}\quad  
   && \tfrac{3}{7} \I\!\P + \tfrac{4}{7} \H\!\P.
   \\[-0.2em]
   &\emph{(r)}\quad  
   && \tfrac{3}{7} \I + \tfrac{4}{7} \H.
\end{alignat*}
}%
\newcommand{\I}{\ket{\dprev{I}}}%
\renewcommand{\H}{\ket{\dprev{H}}}%
\renewcommand{\P}{\ket{\dtest{P}}}%
\newcommand{\N}{\ket{\dtest{N}}}%
\newcommand{\tP}{\texttransparent{0.5}{\ket{\dtest{P}}}}%
\newcommand{\tN}{\texttransparent{0.5}{\ket{\dtest{N}}}}%
Answering $\tfrac{3}{7}\ketb{\dprev{I}}$ is hopefully uncontroversial, but a slight modification of the problem lends itself to two interpretations.
Let us now consider an \emph{unclear testing problem}. The statement is the same as before, but the test result is itself uncertain—e.g., we got an intermediate result that means positive $75\%$ of the time—and thus, observing a positive test yields a distribution over the actual result ($\cst{uncertain} = \tfrac{3}{4}\ket{\dtest{P}} + \tfrac{1}{4}\ket{\dtest{N}}$).
A first solution states all of the information \emph{(i,ii,iii,iv)} and normalizes at the end \emph{(n)}. 
\begin{alignat*}{4}
   &\emph{(i)}\quad 
   && \tfrac{1}{3} \I + \tfrac{1}{3} \H
   \\[-0.2em]
   &\emph{(ii)}\quad  
   && \tfrac{1}{4} \I\!\P + \tfrac{1}{12} \I\!\N 
    + \tfrac{1}{3} \H\!\P + \tfrac{1}{3} \H\!\N
   \\[-0.2em]
   &\emph{(iii)}\quad  
   && \tfrac{3}{16} \I\!\P\!\P 
    + \tfrac{1}{16} \I\!\P\!\N
    + \tfrac{3}{16} \I\!\N\!\P
    + \tfrac{1}{48} \I\!\N\!\N + \\[-0.2em]
   &&& \quad \tfrac{1}{4} \H\!\P\!\P 
    + \tfrac{1}{12} \H\!\P\!\N
    + \tfrac{1}{4} \H\!\N\!\P
    + \tfrac{1}{12} \H\!\N\!\N
  \\[-0.2em]
   &\emph{(iv)}\quad  
   && \tfrac{3}{16} \I\!\P\!\P 
    + \scancel{\tfrac{1}{16} \I\!\P\!\N}
    + \scancel{\tfrac{3}{16} \I\!\N\!\P}
    + \tfrac{1}{48} \I\!\N\!\N + \\[-0.2em]
   &&& \quad \tfrac{1}{4} \H\!\P\!\P 
    + \scancel{\tfrac{1}{12} \H\!\P\!\N}
    + \scancel{\tfrac{1}{4} \H\!\N\!\P}
    + \tfrac{1}{12} \H\!\N\!\N
  \\[-0.2em]
   &\emph{(n)}\quad  
   && \tfrac{9}{26} \I\!\P\!\P + \tfrac{1}{26} \I\!\N\!\N + \tfrac{6}{13} \H\!\P\!\P + \tfrac{2}{13}\H\!\N\!\N
   \\[-0.2em]
   &\emph{(r)}\quad  
  && \tfrac{5}{13} \I + \tfrac{8}{13}\H.
\end{alignat*}
However, a second solution treats the original \emph{testing problem} as a subproblem \emph{(ii,iii,iv)}, normalizing internally $(n_2)$ before normalizing globally $(n_1)$. These two solutions are almost identical—a quick calculation could dismiss both as a probability of roughly $37\%$ of having the illness. Yet, they yield different numbers, arguably answering two different questions.
\begin{alignat*}{4}
   &\emph{(i)}\quad 
   && \tfrac{3}{4} \P + \tfrac{1}{4} \N
   \\[-0.2em]
   &\quad \emph{(ii)}
   && \quad 
      [\tfrac{1}{3}\I + \tfrac{2}{3}\H] \tfrac{3}{4}\!\P
    +  [\tfrac{1}{3}\I + \tfrac{2}{3}\H] \tfrac{1}{4}\!\N
   \\[-0.2em]
   &\quad \emph{(iii)}
   && \quad
      \tfrac{3}{4} \P\![\tfrac{1}{4}\I\!\P + \tfrac{1}{12}\I\!\N + \tfrac{1}{3}\H\!\P + \tfrac{1}{3}\H\!\N] +
   \\[-0.2em]
   &&& \quad 
       \quad \tfrac{1}{4} \N\![\tfrac{1}{4}\I\!\P + \tfrac{1}{12}\I\!\N + \tfrac{1}{3}\H\!\P + \tfrac{1}{3}\H\!\N]
   \\[-0.2em]
   &\quad \emph{(iv)}
   && \quad
      \tfrac{3}{4} \P\![\tfrac{1}{4}\I\!\P + \scancel{\tfrac{1}{12}\I\!\N} +   \tfrac{1}{3} \H\!\P + \scancel{\tfrac{1}{3}\H\!\N}] +
   \\[-0.2em]
   &&& \quad 
       \quad \tfrac{1}{4} \N\![\scancel{\tfrac{1}{4}\I\!\P} + \tfrac{1}{12}\I\!\N + \scancel{\tfrac{1}{3}\H\!\P} + \tfrac{1}{3}\H\!\N]
   \\[-0.2em]
   &\quad \emph{(n$_2$)}
   && \quad
      [\tfrac{3}{7}\I\!\P + \tfrac{4}{7}\H\!\P] \tfrac{3}{4}\!\P +
      [ \tfrac{1}{5}\I\!\N + \tfrac{4}{5}\H\!\P ] \tfrac{1}{4}\!\N
   \\[-0.2em]
   & \emph{(n$_1$)} \quad
   && \tfrac{9}{28} \I\!\P\!\P +
      \tfrac{3}{7} \H\!\P\!\P +
      \tfrac{1}{20} \I\!\N\!\N +
      \tfrac{1}{5} \H\!\N\!\N
   \\[-0.2em]
   & \emph{(r)} \quad
   && \tfrac{13}{35} \I +
      \tfrac{22}{35} \H.
\end{alignat*}

Again, the \emph{unclear test problem} becomes a counterexample to associativity for normalized composition: the following two programs yield different numbers (\Cref{fig:uncleartest}). %
\newsavebox{\codeUncleartestA}
\begin{lrbox}{\codeUncleartestA}
\begin{minipage}{0.45\textwidth}
\begin{lstlisting}[language=Racket,numbers=none,mathescape]
(do (patient) <- prevalence
    (outcome) <- (test patient)
    (result) <- uncertainty
    () <- (observe result outcome)
    return (patient))
\end{lstlisting}
\end{minipage}
\end{lrbox}
\newsavebox{\codeUncleartestB}
\begin{lrbox}{\codeUncleartestB}
\begin{minipage}{0.45\textwidth}
\begin{lstlisting}[language=Racket,numbers=none,mathescape]
(do (result) <- uncertainty
    (patient) <- (do
      (patient) <- prevalence
      (outcome) <- (test patient)
      () <- (observe result outcome)
      return (patient))
    return (patient))
\end{lstlisting}
\end{minipage}
\end{lrbox}
\begin{figure}[h]
\vspace{-1ex}
\centering
\[\usebox{\codeUncleartestA}\quad ≠ \quad
\usebox{\codeUncleartestB}
\]
\vspace{-1em}
\caption{Two interpretations of the \emph{unclear test problem}.\label{fig:uncleartest}}
\end{figure}

It is possibly the second solution—the one that does not simply normalize at the end—the more familiar one to statisticians. Jacobs has studied these two solutions as \emph{Pearl's update} and \emph{Jeffrey's update}, respectively \cite{jacobs2019mathematics}, and we provide an interpretation for both by the end of the article.

\begin{remark}
  In unnormalized probabilistic semantics, encoding these updates requires an additional operator—usually called $\mathsf{normalize}(•)$—used to break associativity. We avoid growing our semantics to unnormalized distributions, defining an associative composition, and then breaking it with a primitive operator. Instead, we accept that composition was not associative to begin with: this allows our semantics to be normalized-by-construction.
\end{remark}

\subsection{Accepting non-associativity}

\emph{Normalized composition is not associative: is there an associative and normalized probabilistic semantics?}
Normalized probabilistic semantics mixes two kinds of effect—stochasticity and partiality—that we encode with the \kl{distribution monad} ($𝐃$) and the \kl{maybe monad} $(𝐌)$, respectively. 

\begin{toappendix}
\begin{definition}[Monad]
  \AP A \intro{monad}, $(𝐓,μ^{𝐓},η^{𝐓})$, consists of an endofunctor, $𝐓 ፡ ℂ → ℂ$, together with \kl{natural transformations}, $η^{𝐓} ፡ X → 𝐓X$ and $μ^{𝐓} ፡ 𝐓 𝐓 X → 𝐓 X$, making the following three diagrams commute.
  \begin{equation*}
\begin{tikzcd}[ampersand replacement=\&]
	{𝐓𝐓𝐓X} \& {𝐓𝐓X} \\
	{𝐓𝐓X} \& {𝐓X}
	\arrow["{𝐓\mu^{𝐓}}", from=1-1, to=1-2]
	\arrow[""{name=0, anchor=center, inner sep=0}, "{\mu^{𝐓}𝐓}"', from=1-1, to=2-1]
	\arrow[""{name=1, anchor=center, inner sep=0}, "{\mu^{𝐓}}", from=1-2, to=2-2]
	\arrow["{\mu^{𝐓}}"', from=2-1, to=2-2]
	\arrow["{①}"{description}, draw=none, from=0, to=1]
\end{tikzcd}
\qquad
\begin{tikzcd}
	{𝐓X} & {𝐓X} \\
	{𝐓𝐓X}
	\arrow[from=1-1, to=1-2]
	\arrow["{\eta^{𝐓}𝐓}"', from=1-1, to=2-1]
	\arrow[""{name=0, anchor=center, inner sep=0}, "{\mu^{𝐓}}"', from=2-1, to=1-2]
	\arrow["{②}"{description}, draw=none, from=1-1, to=0]
\end{tikzcd}
\qquad
\begin{tikzcd}[ampersand replacement=\&]
	{𝐓X} \& {𝐓X} \\
	{𝐓𝐓X}
	\arrow[""{name=0, anchor=center, inner sep=0}, from=1-1, to=1-2]
	\arrow["{𝐓\eta^{𝐓}}"', from=1-1, to=2-1]
	\arrow["{\mu^{𝐓}}"', from=2-1, to=1-2]
	\arrow["{③}"{description}, draw=none, from=0, to=2-1]
\end{tikzcd}
\end{equation*}
\end{definition}
\end{toappendix}

\begin{toappendix}
\begin{definition}[Monoidal monad]
  \AP A \intro{monoidal monad}, $(𝐓, μ^{𝐓}, η^{𝐓}, u, v)$, on a monoidal category $(ℂ,⊗,I)$ consists of a \kl{monad} $(𝐓,μ^{𝐓},η^{𝐓})$, whose underlying functor is lax monoidal: i.e., there exist structural natural transformations,
  \[
  u_{X,Y} ፡ 𝐓 X ⊗ 𝐓 Y → 𝐓(X ⊗ Y)\ \mbox{ and }\ 
  v ፡ I → 𝐓 I,
  \]
  satisfying associativity, $(u_{X,Y} ⊗ \id_{Z}) ⨾ u_{X ⊗ Y, Z} = (\id_{X} ⊗ u_{Y,Z}) ⨾ u_{X, Y ⊗ Z}$, and unitality, $(\id_{𝐓 X} ⊗ v) ⨾ u = \id_{𝐓 X}$ and $(v ⊗ \id_{𝐓 X}) ⨾ u = \id_{𝐓 X}$. In other words, the following diagrams must commute.
  \[
  \begin{tikzcd}[column sep=large, ampersand replacement=\&]
    𝐓 X ⊗ 𝐓 Y ⊗ 𝐓 Z 
    \dar[swap]{u ⊗ 𝐓 Z}
    \rar{𝐓 X ⊗ u} \& 
    𝐓X ⊗ 𝐓(Y ⊗ Z) \dar{u} \\
    𝐓(X ⊗ Y) ⊗ 𝐓Z \rar{u} \& 𝐓(X ⊗ Y ⊗ Z)
  \end{tikzcd}
  \qquad
  \begin{tikzcd}[ampersand replacement=\&]
    𝐓X \rar{\id_{𝐓X} ⊗ v} \dar{v ⊗ \id} \& 𝐓X ⊗ 𝐓I \dar{u_{X,I}} \\
    𝐓I ⊗ 𝐓X \rar{u_{I,X}} \& 𝐓X
  \end{tikzcd}
  \]
  Moreover, the unit and multiplication of the monad must be \kl{monoidal
  natural transformations}, meaning that the following diagrams must commute.
  \[\begin{tikzcd}[ampersand replacement=\&, column sep=small]
    𝐓𝐓 X ⊗ 𝐓𝐓 Y \dar{u} \rar{μ ⊗ μ} \& 𝐓 X ⊗ 𝐓 Y \rar{u} \& 𝐓(X ⊗ Y) \\
    𝐓(𝐓 X ⊗ 𝐓 Y) \rar{𝐓 u} \& 𝐓𝐓(X ⊗ Y) \urar[swap]{μ} \&
  \end{tikzcd}
  \ 
  \begin{tikzcd}[ampersand replacement=\&, column sep=small]
    X ⊗ Y \rar{η} \dar[swap]{η ⊗ η}\& 𝐓(X ⊗ Y) \\
    𝐓 X ⊗ 𝐓 Y \urar[swap]{u} \&
  \end{tikzcd}
  \ 
  \begin{tikzcd}[ampersand replacement=\&, column sep=small]
    I \dar[swap]{\id} \rar{η} \& 𝐓 I \\
    I \urar[swap]{v} \&
  \end{tikzcd}
  \ 
  \begin{tikzcd}[ampersand replacement=\&, column sep=small]
    I \rar{v} \dar[swap]{v} \& 𝐓 I  \\
    𝐓 I \rar{𝐓 v} \& 𝐓𝐓 I \uar[swap]{μ_I}
  \end{tikzcd}\]
\end{definition}
\end{toappendix}

\begin{toappendix}
\begin{definition}[Monoidal distributive law]
  \label{ax:def:monoidal-distributive-law}
  \AP A \kl{monoidal distributive law} between \kl{monoidal monads}, $(𝐒,μ^{𝐒},η^{𝐒},u^{𝐒},v^{𝐒})$ and $(𝐓,μ^{𝐓},η^{𝐓},u^{𝐓},v^{𝐓})$, is a \kl{distributive law}, $ψ_X ፡ 𝐓𝐒X → 𝐒𝐓X$, whose transformation is monoidal, meaning that the two following diagrams must commute.
  \[\begin{tikzcd}[ampersand replacement=\&]
    𝐓𝐒X ⊗ 𝐓𝐒Y \dar[swap]{ψ_X ⊗ ψ_Y} \rar{u^{𝐓}} \& 𝐓(𝐒X ⊗ 𝐒Y) \rar{𝐓u^{𝐒}} \& 𝐓𝐒(X ⊗ Y) \dar{ψ_{X ⊗ Y}} \\
    𝐒𝐓X ⊗ 𝐒𝐓Y \rar{u^{𝐒}} \& 𝐒(𝐓X ⊗ 𝐓Y) \rar{𝐒u^{𝐓}} \& 𝐒𝐓(X ⊗ Y)
  \end{tikzcd}
  \quad
  \begin{tikzcd}[ampersand replacement=\&]
    I \dar[swap]{v} \rar{v} \& 𝐒I \rar{𝐒v} \& 𝐒𝐓I  \\
    𝐓I \rar{𝐓v} \& 𝐒𝐓I \urar[swap]{ψ_{I}} \& 
  \end{tikzcd}\]
\end{definition}  
\end{toappendix}

\begin{definition}[Distribution]
  \AP A \intro{distribution} over a set is a finite formal sum\footnote{A finite formal sum can be understood as an element of the free $ℝ^{+}$-module over the set.} over its elements whose coefficients are positive and add up to exactly $1$.
  The (finitary) \intro[finitary distribution monad]{distribution monad}, $𝐃 ፡ \Set → \Set$, assigns, to each set, the set of distributions over it.
  \[
  𝐃X = 
  \left\{ 
    \sum_{i=0}^n λ_i\!\ket{\cBlu{x_i}} \ \middle|\ 
    \cBlu{x_i} ∈ X, λ_i ∈ \mathbb{R}^{+}, \sum_{i=0}^n λ_i = 1 
  \right\}.\]
  Its \emph{multiplication}, $\smash{μ^{𝐃}} ፡ 𝐃𝐃X → 𝐃X$, is defined by $μ^{𝐃}(\sum\nolimits_{i} λ_{i} \ket{\sum\nolimits_{j} μ_{ij} \ket{x_{ij}} }) = \sum\nolimits_{i} \sum\nolimits_{j} λ_{i}μ_{ij} \ket{  x_{ij} }$; its \emph{unit}, $\smash{η^{𝐃}} ፡ X → 𝐃X$, is defined by $\smash{η^{𝐃}}(x) = \ket{x}$.
\end{definition}

\begin{definition}[Maybe monad]
  \AP The \intro{maybe monad} (sometimes called the \emph{option monad}), $𝐌 ፡
  \Set → \Set$,  assigns to each set $X$, the same set with an extra disjoint
  element usually denoted by $⊥ ∈ 𝐌X$. That is, $$𝐌X = X + \{⊥\}.$$ 
  Its \emph{multiplication}, $μ^{𝐌} ፡ 𝐌𝐌X → 𝐌X$, using $𝐌𝐌X ≅ X + \{⊥₁\} + \{⊥₂\}$, is defined by $μ(x) = x$ for $x ∈ X$, by $μ(⊥₁) = ⊥$ and by $μ(⊥₂) = ⊥$; its \emph{unit}, $η^{𝐌} ፡ X → 𝐌X$, is defined by $η(x) = x$.
\end{definition}

Monads compose when there is a \kl{distributive law} between them; and, under mild compatibility assumptions, a \kl{distributive law} is moreover needed to compose monads \cite{beck1969distributive}.  A \kl{distributive law} of stochasticity over partiality is any consistent way of converting \kl{subdistributions} into \kl{normalized distributions}. For instance, there exists a degenerate \kl{distributive law} of stochasticity over partiality providing the so-called \emph{black-hole semantics}. Black-hole semantics equates any non-full distribution to the empty distribution; explicitly, it is the transformation $(-)^{⊥} ፡ 𝐃𝐌X → 𝐌𝐃X$ defined by
\[
\left( λ_1 \ket{\cBlu{x_1}} + \dots + λ_n \ket{\cBlu{x_n}} \right)^{⊥} =
\begin{cases}
  λ_1 \ket{\cBlu{x_1}} + \dots + λ_n \ket{\cBlu{x_n}},  & 
  \mbox{ if } \sum_{i=0}^{n} λ_i = 1, \\
  0,  & \mbox{ otherwise. } \\
\end{cases}
\]

Alas, \emph{black-hole semantics} is not helpful for our purposes: by equating any probability of failure to failure, it misses the solution to any problem involving \kl{subdistributions}. \emph{Could we find another \kl{distributive law} with normalized semantics?}

We prove that there is no \kl{distributive law} for normalized semantics.\footnote{An equivalent formulation of this fact—in terms of single-point extensions of convex algebras—appeared first in the work of Sokolova and Woracek \cite{sokolova18:termination}. See the Appendix for both a direct proof and an equivalence to this result.} There is no monad structure distributing stochasticity over partiality. Any such semantics must drop some axiom.

\begin{theorem*}[No-go result, \Cref{thm:no-go}]%
  \label{intro:cor:black-hole-unique}
  \AP \kl{Black-hole semantics} constitutes the only \kl{distributive law} between the \kl{distribution monad} and the \kl{maybe monad}, $𝐃𝐌 → 𝐌𝐃$.
\end{theorem*}

Still, \kl{normalized distributions} have a rich structure: both in the discrete and the continuous case, distributions form a \kl{magmad}—a non-associative monad—that arises from a weakening of a \kl{distributive law} interacting with an actual \kl{distributive law}. Associativity is unnecessary: the usual semantics re-emerges from systematically associating to the left (\Cref{sec:magmads}); and associativity may be a limitation, after all: non-associativity enables a \kl{magmadic metalanguage} with multiple modes of update where we can reason about results in causal inference (\Cref{sec:reasoning-metalanguage}).

\begin{recycle}
\end{recycle}

\smallskip

\subsection{Related work}

\emph{Probabilistic programming semantics} employs \kl{subdistributions}, starting with Kozen's \cite{kozen81} and Panangaden's substochastic variant of the Giry monad \cite{giry1982categorical,panangaden1999category}. Since these, both operational \cite{park08,P03:acpl,dLZ12:poslc,culpepper18} and denotational semantics \cite{introProbWood18,zamdzhiev21,vakar2019domain} predominantly account for rejection with \kl[subdistribution]{subdistributional} \cite{borgstrom16,fr19:lambdaprob} or unnormalized semantics \cite{ehrhard2017measurable,ept11:coherencespaces,staton2017commutative,dash23}, with normalization sporadically justifying program equations \cite{staton_et_al_2016}.
\emph{Probabilistic programming languages} either normalize as a program transformation \cite{narayanan16:hakaru}, or let their inference algorithms handle it \cite{GMRBT12:church,toplinWood16:anglican}.

Categorical probability theory, in a line of work starting from Golubtsov, Cho and Jacobs, and Fritz, has abstracted stochastic kernels into \emph{Markov categories} \cite{golubtsov1999axiomatic,cho2019disintegration,fritz2020synthetic}. Further work has abstracted substochastic kernels and partial stochastic kernels into \emph{partial Markov categories} \cite{dilavore23:evidential} and \emph{quasi-Markov categories}  \cite{quasimarkov25}, respectively.  
In particular, string diagrammatic methods can model Pearl's causal interventions \cite{pearl2009causality} by syntactic substitution \cite{fong2013causal,jacobs21:causal,dseparation23}. At the same time, multiple string diagrammatic axiomatizations of \kl{normalization} have been proposed: we highlight those in terms of \emph{normalization boxes}~\cite{lorenz2023causal,jacobs2025compositional} and \emph{partial Markov categories} \cite{dilavore23:evidential}.
Simpson's \emph{probability sheaves} constitute another approach to synthetic probability theory \cite{simpson17,simpson24}; for which Stein recently proposed a comparison \cite{stein25}. Jacobs' \emph{hypernormalization}~\cite{jacobs17} and Garner's \emph{tricocycloids}~\cite{garner2018abstract} form yet another abstract approach to \kl{normalization} that, however, does not address its non-associativity.

\emph{Magmads} and failing distributive laws are relatively infrequent: Munch-Maccagnoni \cite{MunchMaccagnoni13} proposed a non\hyp{}natural monad\hyp{}comonad distributive law \cite{mangel2025classical} and its Kleisli category to unify call\hyp{}by\hyp{}name and call\hyp{}by\hyp{}value. In probability, mass and chance interpretations determine a non-natural distributive law, $𝐃𝐃 → 𝐃𝐃$ \cite{tsai25chancemass}. Weak distributive laws, instead, \cite{street2009weak,bohm2010weak,garner2020vietoris} appear in imprecise probability \cite{goy-petrisan:weak-distributive,sarkis21,liell24:imprecise}, as distributive laws are insufficient for systems with non-determinism and probability \cite{varacca06}. In this non-deterministic case, with applications to probabilistic trace semantics \cite{stark2000complete,silvaS11,cirstea2025}, Bonchi, Sokolova, and Vignudelli distinguished possibilistic monads for \emph{may}, \emph{must}, and \emph{may-must} semantics \cite{dengGHM09,bonchiSV22}; via \kl{support} (c.f.~\cite{fritz2021probability}), we relate these to \kl{normalized semantics}, \kl{partial stochastic semantics}, and \kl{substochastic semantics}, respectively.

\subsection{Contributions}

\kl{Normalization} is not a \kl{distributive law} (\Cref{prop:normalization-is-not-a-distributive-law}) and, moreover, there is no \kl{distributive law} of the \kl{distribution monad} over the \kl{maybe monad} (\Cref{thm:no-go}); still, \kl{normalization} satisfies three of the axioms of \kl{distributive laws} (\Cref{prop:normalization-almost-distributive-law}).

We introduce ``\kl{sesquilaws}'' (\Cref{def:sesquilaw})—\emph{one and a half \kl{distributive laws}} that interact appropriately—and we prove that \kl{normalization} is a \kl{sesquilaw} (\Cref{thm:normalization-sesquilaw}). In this abstract setting, we derive basic program simplifications (e.g.,~\Cref{thm:renormalization}). Examples of \kl{sesquilaw} include possibilistic post-selection (\Cref{prop:postselection}) and \kl{normalization} in standard Borel spaces (\Cref{thm:continuous-normalization-is-a-sesquilaw}); we relate these via \kl{sesquilaw morphisms}, proving that \kl{support}—known to be a monad morphism (c.f.~\cite{fritz2021probability})—is also a \kl{sesquilaw morphism} (\Cref{prop:support-is-a-sesquilaw-homomorphism,prop:finitary-continuous}).

We then study ``\kl{magmads}'' (\Cref{def:magmad}) and we prove that a \kl{sesquilaw} induces a \kl{magmad}, a \kl{monad}, and an action of the latter on the former (\Cref{lemma:right-monad-action}); moreover, in a \kl{sesquilaw}, some effects do always associate (\Cref{lemma:relevance}). We introduce two \kl{magmadic metalanguages}, depending on associativity (\Cref{def:donotation,def:left-do-notation}); and we use them to characterize Bayes' update, Pearl's update, and Jeffrey's update (\Cref{prop:bayes-removal,prop:pearl-jeffrey}).

Finally, we introduce reasoning principles for the metalanguage: representing equality checks (\Cref{sec:partialfrobenius}), conditionals (\Cref{sec:disintegration}), copyable and discardable effects (\Cref{sec:affine-relevant}), and a variant of commutativity that arises only in the non-associative case (\Cref{prop:monad-left-exchange}). We use these to prove the front-door and door-door criteria from causal inference (\Cref{thm:frontdoor,prop:backdoorAdjustment}).

\subsection{Synopsis}

\Cref{sec:distributivelaw} recalls \kl{distributive laws} and proves that there is a single \kl{distributive law} of the \kl{distribution monad} over the \kl{maybe monad}. \Cref{sec:sesquilaw} introduces \kl{sesquilaws}, with examples in continuous probabilistic semantics and possibilistic semantics.
\Cref{sec:magmads-magmoids} constructs \kl{magmads} from \kl{sesquilaws}, introduces the \kl{magmadic metalanguage}, and distinguishes \kl{commutative magmads} and \kl{left-exchanging magmads}. \Cref{sec:update} studies Bayes' update, Pearl's update and Jeffrey's update. \Cref{sec:reasoning-metalanguage} picks some rules of reasoning and derives the front-door and back-door criteria. We provide detailed proofs in the Appendix.

\section{Distributive Laws}
\label{sec:distributivelaw}

\kl{Distributive laws} (\Cref{def:distributiveLaw}) have well-known uses and limitations \cite{beck1969distributive,zwart:nogo}. Briefly, the composition of two monads is not a monad again, but \kl{distributive laws} are sufficient—and necessary under mild conditions—to endow this composition with a monad structure (\Cref{prop:monad-from-distributive-law}). 

\begin{definition}[Distributive law {{\cite{beck1969distributive}}}]
  \label[definition]{def:distributiveLaw}
  \AP A \intro{distributive law} of a \kl{monad}, $(𝐓,\smash{μ^{𝐓}},\smash{η^{𝐓}})$, over a monad on the same category,  $(𝐒,\smash{μ^{𝐒}},\smash{η^{𝐒}})$, is a natural transformation, $ψ ፡ 𝐓𝐒X → 𝐒𝐓X$, that makes the following four diagrams commute.
  \label[definition]{def:distributivelaw}
  \begin{equation*}
    \begin{tikzcd}[ampersand replacement=\&]
      {𝐓𝐓𝐒X} \& {𝐓𝐒X} \& {𝐒𝐓X} \\
      {𝐓𝐒𝐓X} \& {𝐒𝐓𝐓X}
      \arrow["{\mu^{𝐓}𝐒}", from=1-1, to=1-2]
      \arrow[""{name=0, anchor=center, inner sep=0}, "{𝐓\psi}"', from=1-1, to=2-1]
      \arrow["\psi", from=1-2, to=1-3]
      \arrow["{\psi𝐓}"', from=2-1, to=2-2]
      \arrow[""{name=1, anchor=center, inner sep=0}, "{𝐒\mu^{𝐓}}"', from=2-2, to=1-3]
      \arrow["{①}"{description}, draw=none, from=0, to=1]
    \end{tikzcd}
    \qquad
\begin{tikzcd}[ampersand replacement=\&]
	{𝐒X} \& {𝐒𝐓X} \\
	{𝐓𝐒X}
	\arrow[""{name=0, anchor=center, inner sep=0}, "{𝐒\eta^{𝐓}}", from=1-1, to=1-2]
	\arrow["{\eta^{𝐓}𝐒}"', from=1-1, to=2-1]
	\arrow["\psi"', from=2-1, to=1-2]
	\arrow["{②}"{description}, draw=none, from=0, to=2-1]
\end{tikzcd}
  \end{equation*}
  \begin{equation*}
\begin{tikzcd}[ampersand replacement=\&]
	{𝐓𝐒𝐒X} \& {𝐓𝐒X} \& {𝐒𝐓X} \\
	{𝐒𝐓𝐒X} \& {𝐒𝐒𝐓X}
	\arrow["{𝐓\mu^{𝐒}}", from=1-1, to=1-2]
	\arrow[""{name=0, anchor=center, inner sep=0}, "{\psi𝐒}"', from=1-1, to=2-1]
	\arrow["\psi", from=1-2, to=1-3]
	\arrow["{\psi𝐓}"', from=2-1, to=2-2]
	\arrow[""{name=1, anchor=center, inner sep=0}, "{\mu^{𝐒}𝐓}"', from=2-2, to=1-3]
	\arrow["{③}"{description}, draw=none, from=0, to=1]
\end{tikzcd}
    \qquad
\begin{tikzcd}[ampersand replacement=\&]
	{𝐓X} \& {𝐒𝐓X} \\
	{𝐓𝐒X}
	\arrow[""{name=0, anchor=center, inner sep=0}, "{\eta^{𝐒}𝐓}", from=1-1, to=1-2]
	\arrow["{𝐓\eta^{𝐒}}"', from=1-1, to=2-1]
	\arrow["\psi"', from=2-1, to=1-2]
	\arrow["{④}"{description}, draw=none, from=0, to=2-1]
\end{tikzcd}
  \end{equation*}
  These four axioms of a \kl{distributive law} are called ① $𝐓$-multiplicativity, ② $𝐓$-unitality, ③ $𝐒$-multiplicativity, and ④ $𝐒$-unitality, respectively.
\end{definition}

\begin{proposition}[{{\cite[\S9.2. Proposition 2.6]{barr2000toposes}}}]
  \label[proposition]{prop:monad-from-distributive-law}
  A \kl{distributive law}, $ψ_X ፡ 𝐓𝐒X → 𝐒𝐓X$, between two monads, $(𝐒,\smash{μ^{𝐒}},\smash{η^{𝐒}})$ and $(𝐓,\smash{μ^{𝐓}},\smash{η^{𝐓}})$, is necessary and sufficient to induce a monad structure on the composite functor, $(𝐒𝐓,\smash{μ^{𝐒𝐓}_ψ}, \smash{η^{𝐒𝐓}_ψ})$, such that \emph{(i)} the inclusions of the two monads, $𝐒η^{𝐓} ፡ 𝐒X → 𝐒𝐓X$ and $η^{𝐒}𝐓 ፡ 𝐓X → 𝐒𝐓X$ are monad morphisms, \emph{(ii)} the unit arises from tensoring, $\smash{η^{𝐒𝐓}_ψ} = \smash{η^{𝐒}η^{𝐓}}$, \emph{(iii)} and the \emph{middle unitary equation} holds, $\smash{𝐒η^{𝐓}η^{𝐒}𝐓} ⨾ 𝐒ψ𝐓 ⨾ \smash{μ^{𝐒}μ^{𝐓}} = 𝐒𝐓$.
\end{proposition} 

This section first recalls  the standard \kl{distributive law} of subdistributional semantics: a \kl{distributive law} of partiality over stochasticity (\Cref{sec:subdistributions}). The problem becomes to find a \kl{distributive law} in the other direction: a \kl{distributive law} of stochasticity over partiality. We recall the formalization of \emph{black-hole semantics} as a \kl{distributive law} (\Cref{sec:black-hoke-semantics}). Finally, we prove first that \kl{normalization} is not a \kl{distributive law} (\Cref{prop:normalization-is-not-a-distributive-law}) and then the stronger result that there is no \kl{distributive law} but \emph{black-hole semantics} (\Cref{thm:no-go}).

\begin{remark}[Monoidal monads]
  All of the monads we employ are \kl{monoidal monads} (also known as \emph{commutative monads}), and we will consider \kl{monoidal distributive laws} between them. Detailed definitions can be found in the Appendix, but are not necessary for this section.
\end{remark}

\subsection{Example—Subdistributional semantics}
\label{sec:subdistributions}

Let us recall the \kl{distributive law} of partiality over stochasticity (\Cref{prop:partiality-over-stochasticity})—a particular case of Beck's \kl{distributive law} of partiality over any effect \cite[\S4.2]{beck1969distributive}. The monad arising from this \kl{distributive law} can be characterized as the \kl{monad} of \kl{subdistributions} (\Cref{def:subdistribution}).

\begin{proposition}[Partiality over stochasticity]
  \label[proposition]{prop:partiality-over-stochasticity}
  \AP There exists a \kl{distributive law} of the \kl{maybe monad} over the \kl{distribution monad}, $(-)^{•} ፡ 𝐌 𝐃 X → 𝐃 𝐌 X$, defined by $(⊥)^{•} = 0$ and the inclusion $(\sum\nolimits_i λ_i \ketb{x_i})^{•} = \sum\nolimits_i λ_i \ketb{x_i}$.
\end{proposition}

\begin{definition}[Subdistribution]
  \label[definition]{def:subdistribution}
  \AP A \intro{subdistribution} over a set is a finite formal sum over its elements whose coefficients are positive and add up to less or equal than $1$. The (finitary) \intro[subdistribution monad]{subdistribution monad}, $(𝐃𝐌, \smash{μ^{𝐃𝐌}_{•}}, \smash{η^{𝐃𝐌}_{•}})$, assigns, to each set, the set of \kl{subdistributions} over it.
  \[
  𝐃𝐌 X ≅ 
  \left\{ 
    \sum_{i=0}^n λ_i\!\ket{\cBlu{x_i}} \ \middle|\ 
    \cBlu{x_i} ∈ X, λ_i ∈ \mathbb{R}^{+}, \sum_{i=0}^n λ_i ≤ 1 
  \right\}.\]
  Its \emph{multiplication}, $\smash{μ^{𝐃𝐌}_{•}} ፡ 𝐃𝐌𝐃𝐌X → 𝐃𝐌X$, is given by $μ^{𝐃𝐌}_{•}(\sum\nolimits_{i} λ_{i} \ket{\sum\nolimits_{j} μ_{ij} 
  \ket{\cBlu{x_{ij}}} }) = \sum\nolimits_{i} \sum\nolimits_{j} λ_{i}μ_{ij} \ket{  \cBlu{x_{ij}} }$; its \emph{unit}, $\smash{η^{𝐃𝐌}_{•}} ፡ X → 𝐃𝐌X$, is given by $\smash{η^{𝐃𝐌}_{•}}(\cBlu{x}) = \ket{\cBlu{x}}$.
\end{definition}

\subsection{Example—Black-hole semantics}
\label{sec:black-hoke-semantics}

Perhaps surprisingly, there does exist a \kl{distributive law} of stochasticity over partiality \cite{sokolova18:termination}. Let us recall it here (\Cref{prop:partiality-over-stochasticity}). The monad arising from this \kl{distributive law} assigns, to each set, the set of \kl{normalized distributions} over it (\Cref{def:normalized-distribution}). Alas, its multiplication does not normalize (e.g.,~\cite{mohammed2025partializations}), and this monad is not suitable for normalized stochastic semantics.

\begin{proposition}[Stochasticity over partiality]
  \label[proposition]{prop:partiality-over-stochasticity}
  \AP There exists a \kl{distributive law} of the \kl{distribution monad} over the \kl{maybe monad}, $(-)^{⊥} ፡ 𝐃 𝐌 X → 𝐌 𝐃 X$, defined by $(\sum\nolimits_{i=0}^n λ_i\ketb{x_i})^{⊥} = \sum\nolimits_{i=0}^n λ_i\ketb{x_i}$ when $\sum\nolimits_{i=0}^n λ_i = 1$; and $(\sum\nolimits_{i=0}^n λ_i\ketb{x_i})^{⊥} = 0$ when $\sum\nolimits_{i=0}^n λ_i ≠ 1$.
\end{proposition}

\begin{definition}%
  \label[definition]{def:normalized-distribution}%
  \AP A \intro{normalized distribution} over a set is a finite formal sum over its elements whose coefficients are positive and add up to exactly $1$ or $0$. The \intro{partial distribution monad}, $(𝐌𝐃, \smash{μ_{⊥}^{𝐌𝐃}}, \smash{η_{⊥}^{𝐌𝐃}})$, assigns, to each set, the set of \kl{normalized distributions} over it.
  \[
  𝐌𝐃 X ≅ 
  \left\{ 
    \sum_{i=0}^n λ_i\!\ket{\cBlu{x_i}} \ \middle|\ 
    \cBlu{x_i} ∈ X, λ_i ∈ \mathbb{R}^{+}, \sum_{i=0}^n λ_i = 1 \mbox{ or } 0
  \right\}.\]
  Its \emph{multiplication}, $\smash{μ^{𝐌𝐃}_{⊥}} ፡ 𝐌𝐃𝐌𝐃X → 𝐌𝐃X$, is given by $\smash{μ^{𝐌𝐃}_{⊥}}(λ_i \ket{\sum\nolimits_{j} μ_{ij} \ket{\cBlu{x_{ij}}} }) = \sum\nolimits_{i} \sum\nolimits_{j} λ_{i}μ_{ij} \ket{  \cBlu{x_{ij}} }$ whenever $\sum\nolimits_{j} \mu_{ij} = 1$ for each $i = 0,...,n$, and $\smash{μ^{𝐌𝐃}_{⊥}}(λ_i \ket{\sum\nolimits_{j} μ_{ij} \ket{\cBlu{x_{ij}}} }) = 0$ otherwise; its \emph{unit}, $\smash{η^{𝐌𝐃}_{⊥}} ፡ X → 𝐌𝐃X$, is given by $\smash{η^{𝐌𝐃}_{⊥}}(\cBlu{x}) = \ket{\cBlu{x}}$.
\end{definition}

\begin{remark}
  \label[remark]{rem:partialstochastickernels-literature}
  \kl{Partial stochastic kernels} are the paradigmatic example of \emph{quasi-Markov category} \cite{fritz2025empirical,mohammed2025partializations}. While this \kl{quasi-Markov category} plays a role in black-hole semantics, it does not provide updating semantics nor it addresses the problem of normalization: indeed, it is useful precisely because it marks with failure whenever a normalization problem is encountered.
\end{remark}

\subsection{Counterexample—Normalized semantics}

\kl{Normalization} (\Cref{def:normalization}) satisfies all of the axioms of a \kl{distributive law}, except for one: the $𝐃$-multiplicativity axiom (\Cref{prop:normalization-is-not-a-distributive-law}).  We prove a stronger result: there is a single \kl{distributive law} of the \kl{distribution monad} over the \kl{maybe monad} (\Cref{thm:no-go}, c.f.~\cite[Theorem 5.3]{sokolova18:termination}). Next section develops \kl{normalization} not as a \kl{distributive law}, but as a \kl{sesquilaw}.

\begin{definition}[Normalization]
  \label[definition]{def:normalization}
  \AP \intro{Normalization}, $(-)^{∘} ፡ 𝐃𝐌 X → 𝐌𝐃X$, is the natural transformation defined by \emph{(i)} division with the mass of a \kl{subdistribution} whenever not zero, and \emph{(ii)} zero otherwise.
  \[
  \mbox{\emph{(i)}}\quad
  (λ_1\!\ketb{x_1} + \dots + λ_n\!\ketb{x_n})^{∘} = 
  \frac{λ_1}{\sum_i λ_i}\!\ketb{x_1} + \dots + \frac{λ_n}{\sum_i λ_i}\!\ketb{x_n}
  \quad\mbox{ and }\quad\mbox{\emph{(ii)}}\quad
  (0)^{∘} = 0.
  \]
\end{definition}

\begin{theorem}
  \label[theorem]{prop:normalization-is-not-a-distributive-law}
  \kl{Normalization}, $(-)^{∘} ፡ 𝐃𝐌 X → 𝐌𝐃 X$, is not a \kl{distributive law} of the \kl{distribution monad} over the \kl{maybe monad}.
\end{theorem}
\begin{proof}
  {
  \newcommand\x{\cBlu{x}}
  \newcommand\y{\cBlu{y}}
  We prove that the axiom of $𝐃$-multiplicativity fails. Consider the set $X = \{\cBlu{x}, \cBlu{y}\}$ and pick the following distribution of subdistributions:
  \[
  \tfrac{1}{2} \ket{\tfrac{1}{3}\ket{\x} + \tfrac{2}{3} \ket{\cBlu{⊥}}} + 
  \tfrac{1}{2} \ket{\tfrac{2}{3}\ket{\y} + \tfrac{1}{3} \ket{\cBlu{⊥}}} ∈ 𝐃𝐃𝐌X.
  \]

  Let us compute both sides of the $𝐃$-multiplicativity axiom (①, Definition \ref{def:distributivelaw}). On the first side, \emph{(i)} given the \kl{distribution}, we \emph{(ii)} \kl{normalize} internally, \emph{(iii)} \kl{normalize} externally—which happens to be trivial in this case—and \emph{(iv)} multiply.
\begin{align*}
    &\emph{(i)} && 
    \tfrac{1}{2} \ket{\tfrac{1}{3} \ket{\x} + \tfrac{2}{3} \ket{\cBlu{⊥}}} 
    + 
    \tfrac{1}{2} \ket{\tfrac{2}{3} \ket{\cBlu{y}} + \tfrac{1}{3} \ket{\cBlu{⊥}}}
    \\[-0.4em]
    &\emph{(ii)} && 
    \tfrac{1}{2}\ket{\ket{\x}} + \tfrac{1}{2}\ket{\ket{\y}}
    \\[-0.4em]
    &\emph{(iii)} && 
    \tfrac{1}{2}\ket{\ket{\x}} + \tfrac{1}{2}\ket{\ket{\y}} 
    \\[-0.4em]
    &\emph{(iv)} && 
    \tfrac{1}{2}\ket{\x} + \tfrac{1}{2}\ket{\y}.
  \end{align*}
  On the other side, \emph{(i)} given the \kl{distribution}, \emph{(ii)}  we multiply and then \emph{(iii)} \kl{normalize}.
  \begin{align*}
    &\emph{(i)} && 
    \tfrac{1}{2}\ket{ \tfrac{1}{3} \ket{\x} + \tfrac{2}{3} \ket{⊥}} + \tfrac{1}{2}\ket{\tfrac{2}{3} \ket{\y} + \tfrac{1}{3} \ket{⊥}} 
    \\[-0.4em]
    &\emph{(ii)} && 
    \tfrac{1}{6} \ket{\x} + \tfrac{1}{3} \ket{⊥} +\tfrac{1}{3} \ket{\y} + \tfrac{1}{6} \ket{⊥}
    \\[-0.4em]
    &\emph{(iii)} && 
    \tfrac{1}{3}\ket{\x} + \tfrac{2}{3}\ket{\y}.
\end{align*}
  However, $\tfrac{1}{3}\ket{\x} + \tfrac{2}{3}\ket{\y} ≠ \tfrac{1}{2}\ket{\x} + \tfrac{1}{2}\ket{\y}$, contradicting $𝐃$-multiplicativity.
  }
\end{proof}

\begin{theorem}[No-go result]%
  \label{cor:black-hole-unique}
  \label[theorem]{thm:no-go}
  Black-hole semantics constitutes the only \kl{distributive law} between the \kl{distribution monad} and the \kl{maybe monad}, $𝐃𝐌 → 𝐌𝐃$.
\end{theorem}
\begin{proof}
  We will prove that any \kl{distributive law}, $ψ ፡ 𝐃𝐌X → 𝐌𝐃X$, must coincide with black-hole semantics. Given any non-unit scalar $ρ ∈ (0,1)$, by naturality, it must be the case that $ψ(ρ\ketb{x}) = 0$ or $ψ(ρ\ketb{x}) = \ketb{x}$, for all $\cBlu{x} ∈ X$. 

  Let us prove that it must be zero by contradiction. The following two distributions of subdistributions multiply to the same subdistribution: by the $𝐃$-multiplicativity axiom, applying the distributive law internally and externally should yield the same result on
  \[
  \frac{1}{1+ρ}\ket{ρ \ketb{x}} + \frac{ρ}{1+ρ}\ket{\ketb{y}}
  \quad\mbox{and}\quad
  \frac{1}{1+ρ}\ket{\ketb{x}} + \frac{ρ}{1+ρ}\ket{ρ \ketb{y}}.
  \]
  However, we can \emph{(i)} take both distributions of subdistributions; \emph{(i)} normalize internally, using the assumption that $ψ(ρ\ketb{x}) = \ketb{x}$ and that $ψ(\ketb{x}) = 1$ by the $𝐃$-unitality axiom; \emph{(ii)} normalize externally—which leaves the result unchanged—using that 
  $ψ(\sum\nolimits_{i=0}^n λ_i\ketb{x_i}) = \sum\nolimits_{i=0}^n \ketb{x_i}$ whenever $\sum\nolimits_{i=0}^n λ_i = 1$ by the $𝐌$-unitality axiom; and \emph{(iii)} multiply, obtaining two different results.
  \begin{align*}
    \emph{(i)}\quad 
    && \tfrac{1}{1+ρ}\ket{ρ \ketb{x}} + \tfrac{ρ}{1+ρ}\ket{\ketb{y}} 
    & \quad\mbox{ and }\quad 
      \tfrac{ρ}{1+ρ}\ket{\ketb{x}} + \tfrac{1}{1+ρ}\ket{ρ \ketb{y}} \\
    \emph{(ii)}\quad 
    && \tfrac{1}{1+ρ}\ket{\ketb{x}} + \tfrac{ρ}{1+ρ}\ket{\ketb{y}} 
    & \quad\mbox{ and }\quad 
      \tfrac{ρ}{1+ρ}\ket{\ketb{x}} + \tfrac{1}{1+ρ}\ket{ρ \ketb{y}} \\
    \emph{(iii)}\quad 
    && \tfrac{1}{1+ρ}\ket{\ketb{x}} + \tfrac{ρ}{1+ρ}\ket{\ketb{y}} 
    & \quad\mbox{ and }\quad 
      \tfrac{ρ}{1+ρ}\ket{\ketb{x}} + \tfrac{1}{1+ρ}\ket{ρ \ketb{y}} \\
    \emph{(iv)}\quad 
    && \tfrac{1}{1+ρ}\ketb{x} + \tfrac{ρ}{1+ρ}\ketb{y}
    & \quad\mbox{ and }\quad 
      \tfrac{ρ}{1+ρ}\ketb{x} + \tfrac{1}{1+ρ}\ketb{y}.
  \end{align*}
  Because $\smash{\tfrac{ρ}{1+ρ} ≠ \tfrac{1}{1+ρ}}$ unless $ρ = 1$, these two normalized distributions are different, and we must conclude $ψ(ρ\ketb{x}) = 0$. 
  
  Let us prove that the image of any proper \kl{subdistribution} must be zero. Given an arbitrary non-zero and non-full \kl{subdistribution}, $λ_1\ketb{x_1} + \dots + λ_n\ketb{x_n}$, let us define $λ = λ_1 + \dots + λ_n$ and let us consider the following \kl{distribution} of \kl{subdistributions} that multiplies to it.
  \[
  λ \Ket{ \tfrac{λ_1}{λ} \ketb{x_1} + \dots + \tfrac{λ_n}{λ} \ketb{x_n} } + (1 - λ)\Ket{0}
  \]
  By the $𝐃$-multiplicativity axiom, applying the \kl{distributive law} must yield the same result as applying it internally, applying it externally, and multiplying. Applying the \kl{distributive law} internally changes nothing because it is applied to a full and an empty \kl{subdistributions}, respectively. However, applying it externally, by naturality and the previous result, must yield zero. 
  
  We conclude that the distributive law, when applied to any non-full and non-zero \kl{subdistribution} must return zero; and that, by unitality, when applied to a full or zero \kl{subdistribution}, it must return these unchanged. We conclude it coincides with \emph{black-hole semantics}.
\end{proof}

\begin{remark}
   It can be shown that this no-go theorem is equivalent to Sokolova and Woracek's result on the impossibility of one-point extensions of the distributions with black-hole semantics \cite[Theorem 5.3]{sokolova18:termination}. We detail this alternative proof in the Appendix.
\end{remark}

\begin{remark}
  Thus, \kl{normalization} fails one of the axioms of \kl{distributive lax}, but still satisfies most of the axioms and some further properties. Next section abstracts all of the structure of \kl{normalization} into an alternative notion: that of \kl{sesquilaw}.
\end{remark} %
\section{Sesquilaws}
\label{sec:sesquilaw}

\kl{Sesquilaws} are a particular form of \kl{almost-distributive law}. This section defines \kl{sesquilaws} to be sections to an actual \kl{distributive law} that satisfy all \kl{distributive law} axioms except for the first one, which is only satisfied up to idempotent (\Cref{subsec:sesquilaw}). In particular, because they fail one of the axioms, they are \kl{almost-distributive laws}; let us describe this more general concept first (\Cref{sec:almost-distributive-law}).

\subsection{Almost-distributive laws}
\label{sec:almost-distributive-law}

An \intro{almost-distributive law} could be any candidate \kl{distributive law} failing one of the axioms. Specifically, we could define non-$𝐒$-multiplicative, non-$𝐒$-unital, non\hyp{}$𝐓$\hyp{}multiplicative, and non\hyp{}$𝐓$\hyp{}unital \kl{almost-distributive laws}, respectively. In this terminology, a \emph{weak distributive law} \cite{street2009weak,bohm2010weak,garner2020vietoris,goy-petrisan:weak-distributive} would be a non-$𝐓$-unital \kl{almost-distributive law} or, sometimes, a non-$𝐒$-unital \kl{almost-distributive law}. For the rest of the text, however, we focus on non-$𝐓$-multiplicative \kl{almost-distributive laws}, and simply call them \kl[almost-distributive law]{almost\hyp{}distributive laws}.

\begin{definition}[Almost-distributive law]
  \label[definition]{def:almost-distributive-law}
  \AP An \intro{almost-distributive law} of a \kl{monad}, $(𝐓,\smash{μ^{𝐓}},\smash{η^{𝐓}})$, over a \kl{monad} on the same category,  $(𝐒,\smash{μ^{𝐒}},\smash{η^{𝐒}})$, is a natural transformation, $ψ ፡ 𝐓𝐒X → 𝐒𝐓X$, that makes the following three diagrams commute.
  \begin{equation*}
  \begin{tikzcd}[ampersand replacement=\&]
    {𝐓𝐒𝐒X} \& {𝐓𝐒X} \& {𝐒𝐓X} \\
    {𝐒𝐓𝐒X} \& {𝐒𝐒𝐓X}
    \arrow["{𝐓\mu^{𝐒}}", from=1-1, to=1-2]
    \arrow[""{name=0, anchor=center, inner sep=0}, "{\psi𝐒}"', from=1-1, to=2-1]
    \arrow["\psi", from=1-2, to=1-3]
    \arrow["{\psi𝐓}"', from=2-1, to=2-2]
    \arrow[""{name=1, anchor=center, inner sep=0}, "{\mu^{𝐒}𝐓}"', from=2-2, to=1-3]
    \arrow["{①}"{description}, draw=none, from=0, to=1]
  \end{tikzcd}
    \quad
    \begin{tikzcd}[ampersand replacement=\&]
      {𝐒X} \& {𝐒𝐓X} \\
      {𝐓𝐒X}
      \arrow[""{name=0, anchor=center, inner sep=0}, "{𝐒\eta^{𝐓}}", from=1-1, to=1-2]
      \arrow["{\eta^{𝐓}𝐒}"', from=1-1, to=2-1]
      \arrow["\psi"', from=2-1, to=1-2]
      \arrow["{②}"{description}, draw=none, from=0, to=2-1]
    \end{tikzcd}
    \quad
  \begin{tikzcd}[ampersand replacement=\&]
    {𝐓X} \& {𝐒𝐓X} \\
    {𝐓𝐒X}
    \arrow[""{name=0, anchor=center, inner sep=0}, "{\eta^{𝐒}𝐓}", from=1-1, to=1-2]
    \arrow["{𝐓\eta^{𝐒}}"', from=1-1, to=2-1]
    \arrow["\psi"', from=2-1, to=1-2]
    \arrow["{③}"{description}, draw=none, from=0, to=2-1]
  \end{tikzcd}
  \end{equation*}
\end{definition}

\begin{proposition}
  \label{prop:normalization-almost-distributive-law}
  \kl{Normalization}, $(-)^{∘} ፡ 𝐃 𝐌 X → 𝐌 𝐃 X$, forms an \kl{almost-distributive law}.
\end{proposition}
\begin{proof}
  We must prove that the following three diagrams commute. 
  \begin{equation*}
    \begin{tikzcd}[ampersand replacement=\&,column sep=scriptsize]
    𝐃𝐌𝐌 X \dar[swap]{ψ 𝐌} \rar{𝐃 μ^{𝐌}} \& 𝐃𝐌X \rar{ψ} \& 𝐌𝐃X \\
    𝐌𝐃𝐌 X \rar{𝐌ψ} \& 𝐌𝐌𝐃X \urar[swap]{μ^{𝐌} 𝐃} \&
    \end{tikzcd}
    \quad
  \begin{tikzcd}[ampersand replacement=\&,column sep=scriptsize]
    𝐃 X \dar[swap]{𝐃 η^{𝐌}} \rar{η^{𝐌} 𝐃} \& 𝐌𝐃 X \\
    𝐃𝐌 X \urar[swap]{ψ} \& 
  \end{tikzcd}
\quad
\begin{tikzcd}[ampersand replacement=\&,column sep=scriptsize]
  𝐌 X \dar[swap]{η^{𝐃} 𝐌} \rar{𝐌 η^{𝐃}} \& 𝐌 𝐃 X \\
  𝐃 𝐌 X \urar[swap]{ψ} \& 
  \end{tikzcd}
\end{equation*}

  For $𝐌$\hyp{}multiplicativity, let us note that an element of $𝐃𝐌𝐌X$ is equivalently a \kl{distribution} over $X + \{⊥_1\} + \{⊥_2\}$. Let us pick a generic element, $λ₁\ketb{x₁} + \dots + λ_n\ketb{x_n} + μ₁\ketb{⊥₁} + μ₂\ketb{⊥₂}$, and let us reason by cases. If $μ₁ + μ₂ = 1$, then normalizing in any order yields the zero distribution, in the first case because $μ₁ + μ₂ = 1$ and in the second case because, whenever $μ₁ ≠ 1$, it must be the case that $μ₂/μ₂ = 1$. Let us now assume that $μ₁ + μ₂ ≠ 1$ and let us define $λ = λ_1 + \dots + λ_n$. On the one side, we \emph{(i)} consider the generic element, \emph{(ii)} multiply using the \kl{maybe monad}, and \emph{(iii)} \kl{normalize},
  \begin{align*}
    \emph{(i)}\quad &
    λ₁\ketb{x₁} + \dots + λ_n\ketb{x_n} + μ_1\ketb{⊥_1} + μ_2\ketb{⊥_2} \\[-0.4em]
    \emph{(ii)}\quad &
    λ₁\ketb{x₁} + \dots + λ_n\ketb{x_n} + (μ_1 + μ_2)\ketb{⊥} \\[-0.4em]
    \emph{(iii)}\quad &
    \tfrac{λ₁}{λ}\ketb{x₁} + \dots + \tfrac{λ_n}{λ}\ketb{x_n}.
  \end{align*}
  On the other side, we \emph{(i)} consider the generic element, \emph{(ii)} normalize at $⊥₂$, \emph{(iii, iv)} normalize at $⊥₁$ and multiply using the \kl{maybe monad}. We use the fact that $(λ_i / λ + μ_1)(λ + μ_1 / λ) = λ_i / λ$.
  \begin{align*}
    \emph{(i)}\quad &
      λ_1\ketb{x_1} + \dots + λ_n\ketb{x_n}
      + μ_1\ketb{⊥_1} + μ_2\ketb{⊥_2} \\[-0.4em]
    \emph{(ii)}\quad &
      \tfrac{λ₁}{λ + μ_1}\ketb{x₁} + \dots + \tfrac{λ_n}{λ + μ_1}\ketb{x_n} + \tfrac{μ_1}{λ + μ_1}\ketb{⊥_1} \\[-0.4em]
    \emph{(iii)}\quad &
    \tfrac{λ₁}{λ}\ketb{x₁} + \dots + \tfrac{λ_n}{λ}\ketb{x_n}.
  \end{align*}

  For $𝐌$\hyp{}unitality, we must check that the normalization of a normalized distribution is itself: we note that, in that case, $λ_1 + \dots + λ_n = 1$.
  Finally, for $𝐃$\hyp{}unitality, we reason by cases on $𝐌X$. We directly check that $(η^{𝐃}{(⊥)})^{∘} = (\ket{⊥})^{∘} = 0$ and
  $(η^{𝐃}{(x)})^{∘} = (\ketb{x})^{∘} = \ketb{x}$.
\end{proof}

\kl{Normalization}, the \kl{almost-distributive law}, induces the
Kleisli \kl{magmoid} of \kl{normalized kernels}, $\Norm$. 
Together with the \kl{distributive law} describing \kl{subdistributions},
we have \emph{``one and a half''} \kl{distributive laws} that interact with
each other: a \kl{sesquilaw}.

\subsection{Sesquilaws}
\label{subsec:sesquilaw}

To recap, \kl{normalization}, $\norm{(-)} ፡ 𝐃𝐌 → 𝐌𝐃$, satisfies all
\kl{distributive law} axioms except for the $𝐃$-mul\-ti\-pli\-ca\-tivity axiom.
Still, \kl{normalization} satisfies an equation resembling this missing
multiplicativity: $((f)^{∘•} 𑊩 g)^{∘} = (f 𑊩 g)^{∘}$, for any two
\kl{substochastic kernels}. Careful inspection reveals that the
$𝐃$-multiplicativity axiom holds \emph{up-to-an-idempotent}: the
\kl{distributive law} of \kl{subdistributions}, $(-)^{•} ፡ 𝐌𝐃 → 𝐃𝐌$, is the
partial inverse inducing this idempotent. 

This section introduces \kl{sesquilaws}, which abstract this situation with a single extra equation (\Cref{def:distributive-swap}). This equation consists of multiplicativity up to the idempotent determined by a section-retraction pair of two \kl{distributive law} candidates. \kl{Normalization} is a \kl{sesquilaw} (\Cref{thm:normalization-sesquilaw}), and we use the \kl{sesquilaw} axioms to derive a result on \kl{normalization} (\Cref{thm:renormalization}).

\begin{definition}[Sesquilaw]%
  \label[definition]{def:distributive-swap}%
  \label[definition]{def:sesquilaw}%
  \AP A \intro{sesquilaw}, $(𝐒,𝐓,m,n)$, of a \kl{monad} $(𝐓,\smash{μ^𝐓},\smash{η^𝐓})$ over a \kl{monad} on the same category $(𝐒,\smash{μ^𝐒},\smash{η^𝐒})$ consists of a \kl{distributive law}, $m ፡ 𝐒 𝐓 X → 𝐓 𝐒 X$, and an \kl{almost distributive law}, $n ፡ 𝐓 𝐒 X → 𝐒 𝐓 X$, additionally making the following two diagrams commute.
  \begin{equation*}
  \begin{tikzcd}[ampersand replacement=\&]
    {𝐓𝐒𝐓X} \& {𝐒𝐓𝐓X} \& {𝐒𝐓X} \\
    {𝐓𝐓𝐒X} \& {𝐓𝐒X}
    \arrow["{n𝐓}", from=1-1, to=1-2]
    \arrow["{𝐓m}"', from=1-1, to=2-1]
    \arrow["{𝐒\mu^{𝐓}}", from=1-2, to=1-3]
    \arrow[""{name=0, anchor=center, inner sep=0}, "{\mu^{𝐓}𝐒}"', from=2-1, to=2-2]
    \arrow["n"', from=2-2, to=1-3]
    \arrow["{①}"{description}, draw=none, from=1-2, to=0]
  \end{tikzcd} 
  \qquad
\begin{tikzcd}[ampersand replacement=\&]
	{𝐒𝐓X} \& {𝐒𝐓X} \\
	{𝐓𝐒X}
	\arrow[no head, from=1-1, to=1-2]
	\arrow["m"', from=1-1, to=2-1]
	\arrow[""{name=0, anchor=center, inner sep=0}, "n"', from=2-1, to=1-2]
	\arrow["{②}"{description}, draw=none, from=1-1, to=0]
\end{tikzcd}
\end{equation*}
\end{definition}

\begin{theorem}
  \label[theorem]{thm:normalization-sesquilaw}
  Normalized semantics, $(-)^{∘} ፡ 𝐃𝐌X → 𝐌𝐃X$, and subdistributional semantics, $(-)^{•} ፡ 𝐌𝐃X → 𝐃𝐌X$, form a \kl{sesquilaw} of the \kl{distribution monad} over the \kl{maybe monad}.
\end{theorem}
\begin{proof}
  We already know that \kl{normalization} and the inclusion into
  \kl{subdistributions} are mutual inverses, and that \kl{normalization} forms an \kl{almost distributive law} (Proposition~\ref {prop:normalization-almost-distributive-law}).
  Let us prove the axiom of \kl{sesquilaws}.
  \[\begin{tikzcd}[ampersand replacement=\&]
	𝐃𝐌𝐃X \dar[swap]{𝐃 m} \rar{n 𝐃} \& 𝐌𝐃𝐃X \rar{𝐌 μ} \& 𝐌𝐃X \\
  𝐃𝐃𝐌X \rar{μ 𝐌} \& 𝐃𝐌X \urar[swap]{n} 
  \end{tikzcd}\]
  
  Let us note that the generic element of $𝐃𝐌𝐃X$ is, equivalently, a \kl{subdistribution} of \kl{distributions}, $λ_1\ket{μ_{11}\ketb{x_{11}} + \dots + μ_{1{m_1}}\ketb{x_{1m_1}}} + \dots + λ_n\ket{μ_{n1}\ketb{x_{n1}} + \dots + μ_{n{m_n}}\ketb{x_{nm_n}}}$, where $μ_{i1} + \dots + μ_{in_i} = 1$ for each $i = 0, ..., n$.
  Let us first \emph{(i)} consider the generic element, \emph{(ii)} normalize and flatten the \kl{distributions}.
  \begin{align*}
    \emph{(i)}\quad &
    λ_1\ket{μ_{11}\ketb{x_{11}} + \dots + μ_{1{m_1}}\ketb{x_{1m_1}}} + \dots + λ_n\ket{μ_{n1}\ketb{x_{n1}} + \dots + μ_{n{m_n}}\ketb{x_{nm_n}}} \\[-0.4em]
    \emph{(ii)}\quad &
    \tfrac{λ_1}{λ}\ket{μ_{11}\ketb{x_{11}} + \dots + μ_{1{m_1}}\ketb{x_{1m_1}}} + \dots + \tfrac{λ_n}{λ}\ket{μ_{n1}\ketb{x_{n1}} + \dots + μ_{n{m_n}}\ketb{x_{nm_n}}} \\[-0.4em]
    \emph{(iii)}\quad &
    \tfrac{λ_1μ_{11}}{λ}\ketb{x_{11}} + \dots + \tfrac{λ_1μ_{1m_1}}{λ}\ketb{x_{1m_1}} + \dots + \tfrac{λ_nμ_{n1}}{λ} \ketb{x_{n1}} + \dots + \tfrac{λ_nμ_{nm_n}}{λ}\ketb{x_{nm_n}}.
  \end{align*}
  Let us now \emph{(i)} consider the generic element, \emph{(ii)} flattening while regarding the \kl{distributions} as \kl{subdistributions}, and \emph{(iii)} normalizing.
  \begin{align*}
    \emph{(i)}\quad &
    λ_1\ket{μ_{11}\ketb{x_{11}} + \dots + μ_{1{m_1}}\ketb{x_{1m_1}}} + \dots + λ_n\ket{μ_{n1}\ketb{x_{n1}} + \dots + μ_{n{m_n}}\ketb{x_{nm_n}}} \\
    \emph{(ii)}\quad &
    λ_1μ_{11}\ketb{x_{11}} + \dots + λ_1μ_{1m_1}\ketb{x_{1m_1}} + \dots + λ_nμ_{n1} \ketb{x_{n1}} + \dots + λ_nμ_{nm_n} \ketb{x_{nm_n}} \\[-0.4em]
    \emph{(iii)}\quad &
    \tfrac{λ_1μ_{11}}{\smash{\sum\nolimits_{ij} λ_i μ_{ij}}}\ketb{x_{11}} + \dots + \tfrac{λ_1μ_{1m_1}}{\smash{\sum\nolimits_{ij} λ_i μ_{ij}}}\ketb{x_{1m_1}} + \dots + \tfrac{λ_nμ_{n1}}{\smash{\sum\nolimits_{ij} λ_i μ_{ij}}} \ketb{x_{n1}} + \dots + \tfrac{λ_nμ_{nm_n}}{\smash{\sum\nolimits_{ij} λ_i μ_{ij}}}\ketb{x_{nm_n}}.
  \end{align*}
  The last step is to see that $\sum\nolimits_{ij} λ_i μ_{ij} = \sum\nolimits_{ij} λ_i (\sum\nolimits_{j} μ_{ij}) = \sum\nolimits_{ij} λ_i$, which follows from the fact that $μ_{i1} + \dots + μ_{in_i} = 1$ for each $i = 0, ..., n$.
\end{proof}

\kl{Sesquilaw} contain enough structure to abstractly prove useful facts about \kl{normalization}. For instance, following rule—which we call \emph{renormalization} (\Cref{thm:renormalization})—states that normalizing after each multiplication is equivalent to normalizing at the end. It will serve us later to justify the soundness of normalizing and resampling after each transformation (\Cref{lemma:right-monad-action}).

\begin{lemmarep}[Renormalization]
  \label{thm:renormalization}
  \AP Any \kl{sesquilaw}, $(𝐒,𝐓,m,n)$, induces an idempotent, $k = (n ⨾ m) ፡ 𝐓𝐒 → 𝐓𝐒$.  This idempotent is left-absorptive, meaning that the following diagram commutes.
\[\begin{tikzcd}[ampersand replacement=\&]
	{𝐓𝐒𝐓𝐒 X} \& {𝐓𝐒X} \& {𝐒𝐓X} \\
	{𝐒𝐓𝐓𝐒 X} \& {𝐓𝐒𝐓𝐒 X} \& {𝐓𝐒X}
	\arrow["{\mu^{𝐓𝐒}_m}", from=1-1, to=1-2]
	\arrow[""{name=0, anchor=center, inner sep=0}, "{n𝐓𝐒}"', from=1-1, to=2-1]
	\arrow["n", from=1-2, to=1-3]
	\arrow["{m𝐓𝐒}", from=2-1, to=2-2]
	\arrow["{\mu^{𝐓𝐒}_m}", from=2-2, to=2-3]
	\arrow[""{name=1, anchor=center, inner sep=0}, "n"', from=2-3, to=1-3]
	\arrow[draw=none, from=0, to=1]
\end{tikzcd}\]
\end{lemmarep}
\begin{proofsketch}
  The proof uses naturality, the definition of the monad arising from a \kl{distributive law}, the $𝐒$-multiplicativity of the \kl{almost-distributive law} $n$, the first \kl{sesquilaw} axiom, and the second \kl{sesquilaw} axiom. See the Appendix for details.
\end{proofsketch}
\begin{proof}
  We employ commutative diagrams. We use ⓪ naturality, ① the definition of the monad arising from a \kl{distributive law}, ② the $𝐒$-multiplicativity of the \kl{almost-distributive law} $n$, ③ the first \kl{sesquilaw} axiom, and ④ the second \kl{sesquilaw} axiom.
  \[\begin{tikzcd}[ampersand replacement=\&]
    \&\&\& {𝐓𝐒𝐓𝐒 X} \& {𝐓𝐒X} \&\& \\
    \&\& {𝐓𝐒𝐓𝐒 X} \& {𝐓𝐓𝐒𝐒 X} \& {𝐓𝐒X} \\
    \& {𝐓𝐒𝐓𝐒 X} \& {𝐓𝐓𝐒𝐒 X} \& {𝐓𝐒𝐒 X} \& {𝐓𝐒X} \& {𝐒𝐓X} \\
    \& {𝐓𝐒𝐓𝐒 X} \& {𝐓𝐓𝐒𝐒 X} \& {𝐓𝐒𝐒 X} \& {𝐒𝐓𝐒X} \& {𝐒𝐒𝐓X} \& {𝐒𝐓X} \\
    {𝐓𝐒𝐓𝐒 X} \& {𝐒𝐓𝐓𝐒 X} \& {𝐓𝐒𝐓𝐒 X} \& {𝐒𝐓𝐓𝐒X} \& {𝐒𝐓𝐒X} \& {𝐒𝐒𝐓X} \& {𝐒𝐓X} \\
    \& {𝐓𝐒𝐓𝐒 X} \& {𝐒𝐓𝐓𝐒 X} \&\& {𝐒𝐓𝐒X} \& {𝐒𝐒𝐓X} \& {𝐒𝐓X} \\
    \& {𝐓𝐒𝐓𝐒 X} \& {𝐓𝐓𝐒𝐒 X} \& {𝐓𝐒𝐒 X} \& {𝐒𝐓𝐒X} \& {𝐒𝐒𝐓X} \& {𝐒𝐓X} \\
    \&\&\&\& {𝐓𝐒𝐒 X} \& {𝐓𝐒 X} \& {𝐒𝐓X}
    \arrow["{\mu^{𝐓𝐒}}", from=1-4, to=1-5]
    \arrow[""{name=0, anchor=center, inner sep=0}, no head, from=1-4, to=2-3]
    \arrow[""{name=1, anchor=center, inner sep=0}, no head, from=1-5, to=2-5]
    \arrow["{𝐓m𝐒}", from=2-3, to=2-4]
    \arrow[no head, from=2-3, to=3-2]
    \arrow["{\mu^{𝐓}\mu^{𝐒}}", from=2-4, to=2-5]
    \arrow[""{name=2, anchor=center, inner sep=0}, no head, from=2-4, to=3-3]
    \arrow[""{name=3, anchor=center, inner sep=0}, no head, from=2-5, to=3-5]
    \arrow["{𝐓m𝐒}", from=3-2, to=3-3]
    \arrow[no head, from=3-2, to=4-2]
    \arrow["{\mu^{𝐓}𝐒𝐒}", from=3-3, to=3-4]
    \arrow[no head, from=3-3, to=4-3]
    \arrow["{𝐓\mu^{𝐒}}", from=3-4, to=3-5]
    \arrow[""{name=4, anchor=center, inner sep=0}, no head, from=3-4, to=4-4]
    \arrow["n", from=3-5, to=3-6]
    \arrow[""{name=5, anchor=center, inner sep=0}, no head, from=3-6, to=4-7]
    \arrow["{𝐓m𝐒}", from=4-2, to=4-3]
    \arrow[""{name=6, anchor=center, inner sep=0}, no head, from=4-2, to=5-3]
    \arrow["{\mu^{𝐓}𝐒𝐒}", from=4-3, to=4-4]
    \arrow["{n𝐒}", from=4-4, to=4-5]
    \arrow["{𝐒n}", from=4-5, to=4-6]
    \arrow[""{name=7, anchor=center, inner sep=0}, no head, from=4-5, to=5-5]
    \arrow["{\mu^{𝐒}𝐓}", from=4-6, to=4-7]
    \arrow[no head, from=4-6, to=5-6]
    \arrow[no head, from=4-7, to=5-7]
    \arrow["{n𝐓𝐒}", from=5-1, to=5-2]
    \arrow[no head, from=5-1, to=6-2]
    \arrow["{m𝐓𝐒}", from=5-2, to=5-3]
    \arrow[""{name=8, anchor=center, inner sep=0}, no head, from=5-2, to=6-3]
    \arrow["{n𝐓𝐒}", from=5-3, to=5-4]
    \arrow["{𝐒\mu^{𝐓}𝐒}", from=5-4, to=5-5]
    \arrow[""{name=9, anchor=center, inner sep=0}, no head, from=5-4, to=6-3]
    \arrow["{𝐒n}", from=5-5, to=5-6]
    \arrow[no head, from=5-5, to=6-5]
    \arrow["{\mu^{𝐒}𝐓}", from=5-6, to=5-7]
    \arrow[no head, from=5-6, to=6-6]
    \arrow[no head, from=5-7, to=6-7]
    \arrow["{n𝐓𝐒}", from=6-2, to=6-3]
    \arrow[""{name=10, anchor=center, inner sep=0}, no head, from=6-2, to=7-2]
    \arrow["{𝐒\mu^{𝐓}𝐒}", from=6-3, to=6-5]
    \arrow["{𝐒n}", from=6-5, to=6-6]
    \arrow[""{name=11, anchor=center, inner sep=0}, no head, from=6-5, to=7-5]
    \arrow["{\mu^{𝐒}𝐓}", from=6-6, to=6-7]
    \arrow[no head, from=6-6, to=7-6]
    \arrow[no head, from=6-7, to=7-7]
    \arrow["{𝐓m𝐒}", from=7-2, to=7-3]
    \arrow["{\mu^{𝐓}𝐒𝐒}", from=7-3, to=7-4]
    \arrow["{n𝐒}", from=7-4, to=7-5]
    \arrow[""{name=12, anchor=center, inner sep=0}, no head, from=7-4, to=8-5]
    \arrow["{𝐒n}", from=7-5, to=7-6]
    \arrow["{\mu^{𝐒}𝐓}", from=7-6, to=7-7]
    \arrow[""{name=13, anchor=center, inner sep=0}, no head, from=7-7, to=8-7]
    \arrow["{𝐓\mu^{𝐒}}", from=8-5, to=8-6]
    \arrow["n", from=8-6, to=8-7]
    \arrow["{①}"{description}, draw=none, from=0, to=1]
    \arrow["{⓪}"{description}, draw=none, from=2, to=3]
    \arrow["{②}"{description}, draw=none, from=4, to=5]
    \arrow["{③}"{description}, draw=none, from=6, to=7]
    \arrow["{④}"{description}, draw=none, from=8, to=9]
    \arrow["{③}"{description}, draw=none, from=10, to=11]
    \arrow["{②}"{description}, draw=none, from=12, to=13]
  \end{tikzcd}\]
  This concludes the proof.
\end{proof}

\subsection{Example — Normalized boolean semantics is associative}

Possibilistic semantics is encoded by the \kl{affine powerset monad} (\Cref{def:affine-powerset}), which distributes over the \kl{maybe monad} with a \kl{sesquilaw} (\Cref{prop:postselection}). This \kl{sesquilaw} is a simple form of normalization—removing the possibility of failure—that we call \emph{postselection} (\Cref{def:postselection}).

\begin{definition}[Affine powerset monad]
  \label[definition]{def:affine-powerset}
  \AP The \intro{affine powerset monad} (or, \emph{non-empty finitary powerset
  monad}), $𝐑 ፡ \Set → \Set$, assigns, to each set, its set of non-empty finite subsets,
  \[ 𝐑 X = \{ U ⊆ X \mid U \mbox{ finite, and } U ≠ \varnothing \}. \]
  Its \emph{multiplication}, $\smash{μ^{𝐑}} ፡ 𝐑𝐑 X → 𝐑 X$, is defined by $μ^{𝐑}(A) = \bigcup\nolimits_{X ∈ A} X$; its \emph{unit}, $\smash{η^{𝐑}} ፡ X → 𝐑 X$, is defined by $\smash{η^{𝐑}}(x) = \{x\}$.
\end{definition}

\begin{definition}[Post-selection]
  \label[definition]{def:postselection}
  \AP \kl{Post-selection}, $\psel(-) ፡ 𝐑𝐌X → 𝐌𝐑X$, is the natural transformation that removes failure if it can and fails only if it must. It is defined by 
  \[\psel(A \cup \{⊥\}) = \psel(A) = A, \mbox{ for each $A ∈ 𝐑X$, and } \psel(\{⊥\}) = ∅.\]
  In other words, it is the join-preserving function satisfying $\psel(\{x\}) = \{x\}$ and $\psel(\{⊥\}) = ∅$. \kl{Post-selection} is a retraction of the inclusion $\isel ፡ 𝐌𝐑X → 𝐑𝐌X$ defined by $\isel(⊥) = \{∅\}$ and $\isel(A) = A$ for each $A ∈ 𝐑X$.
\end{definition}

\begin{propositionrep}[Normalized boolean semantics]
  \label{prop:postselection}
  \kl{Post-selection} and inclusion form a \kl{sesquilaw}, $(𝐌,𝐑,\isel,\psel)$, of the \kl{affine powerset monad} over the \kl{maybe monad}.
\end{propositionrep}
\begin{proof}
  A better description of \kl{post-selection} start by regarding non-empty subsets, the elements of $𝐑X$, as \emph{non-null predicates} on $X$: predicates that are true on at least one element. Moreover, we define a zero-handling division $(\oslash)$ that answers with the failure element when the denominator is zero. We take this point of view for this proof, as we shall see it simplifies it considerably. %

  \kl{Post-selection} of a \emph{predicate}, $α ፡ X → \{0,1\}$ or $α ∈ 𝐑𝐌X$, is defined by $α(x)$ except when it is null: in that case, we return the extra element of the \kl{maybe monad},
  \[\psel(α) = \bdiv{α(x)}{\bigvee_{x ∈ X} α(x)}.\]

  Let us check that \kl{post-selection} is a \kl{monoidal natural transformation}. It suffices to note that the following formula holds.
  \[
  \bdiv{α(x)}{\bigvee_{x ∈ X} α(x)} \land \bdiv{β(x)}{\bigvee_{x ∈ X} β(x)} = 
  \bdiv{α(x) \land β(x)}{\bigvee_{x ∈ X} α(x) \land \bigvee_{x ∈ X} β(x)}.
  \] 

  Let us check the \kl{sesquilaw} axiom. It is well-known that there exists a distributive law of the maybe monad over any other $\Set$-monad.
  Let $ψ ∈ 𝐑 𝐌𝐑 X$ be, equivalently, a predicate on
  non-null predicates. On the lower path, we use the monad multiplication and \kl{normalization}.
  \begin{align*}
    𝐌 μ(\psel{(ψ)})(x) = 
    \bigvee_{α ∈ 𝐑 X} \psel{(ψ)}(α) \land α(x) = 
    \bigvee_{α ∈ 𝐑 X} \bdiv{ψ(α)}{\bigvee_{β ∈ 𝐑 X} ψ(β)} \land α(x).
  \end{align*}
  On the right hand side, we use the \kl{normalization} and monad multiplication.
  \begin{align*}
    \psel{(μ(ψ(-^{•})))}(x) 
    &= \bdiv{μ(ψ(-^{•}))(x)}{\bigvee_{y ∈ X} μ(ψ(-^{•}))(y)}
    = \bdiv{\bigvee_{α ∈ 𝐃X} ψ(α^{•}) \land α(x)}{\bigvee_{y ∈ X} \bigvee_{α ∈ 𝐑 X} ψ(α^{•}) \land α(y)} \\
    &= \bdiv{\bigvee_{α ∈ 𝐃X} ψ(α^{•}) \land α(x)}{\bigvee_{α ∈ 𝐃X} ψ(α^{•}) \land \bigvee_{y ∈ X} ψ(y)} 
    = \bdiv{\bigvee_{α ∈ 𝐃X} ψ(α^{•}) \land α(x)}{\bigvee_{α ∈ 𝐃X} ψ(α^{•})}.
  \end{align*}
  The last step uses that all $α ∈ 𝐑 X$ are non-null predicates: $\bigvee_{y ∈ X} α(y)$ must be exactly $1$. The rest of the proof is analogous.
\end{proof}

\begin{remark}
  The \kl{magmad} associated to the \kl{sesquilaw}, $(𝐌𝐑,\smash{μ^{𝐌𝐑}_{\psel}}, \smash{η^{𝐌𝐑}_{\psel}})$ happens to be a \kl{monad}—it is equivalent to the powerset monad. The \kl{monad} associated to the \kl{sesquilaw}, $(𝐑𝐌, \smash{μ^{𝐑𝐌}_{\isel}}, \smash{η^{𝐑𝐌}_{\isel}})$, is the \kl{monad} of subaffine relations used for \emph{may-must semantics} (e.g.,~\cite{bonchiSV22}). Moreover, there exists another sesquilaw determining possibilistic \emph{black-hole semantics} \cite{sarkis21,bonchiSV22}.
\end{remark}

\subsection{Example — Normalized continuous probabilistic semantics is not associative}

\kl{Sesquilaws} are abstract enough to unify the discrete and the continuous case; this section details \kl{normalization} in standard Borel spaces. We introduce the  Giry monad (\Cref{def:giry}) and Panangaden's substochastic variant (\Cref{def:panangaden}).

\begin{definition}[Giry monad, {{\cite{lawvere1962category,giry1982categorical}}}]
  \label[definition]{def:giry}
  \AP\phantomintro{Giry monad}The functor taking a Borel space to the set of probability measures over it, $𝓓 ፡ \mathbf{BorelMes} → \mathbf{BorelMes}$, forms a monad, $(𝓓,\smash{η^{𝓓}},\smash{μ^{𝓓}})$, whose unit, $\smash{η^{𝓓}} ፡ X → 𝓓X$, is the Dirac delta, $\smash{η^{𝓓}}(x) = δ(x;-)$, and whose multiplication, $\smash{μ^{𝓓}} ፡ 𝓓𝓓X → 𝓓X$, is given by Lebesgue integration,
  \[
    μ^{𝓓}(Ω)(U) = \int_{ν ∈ 𝓓X} ν(U) · Ω(\dif ν).
  \]
\end{definition}

\begin{definition}[Panangaden monad, {{\cite[\S 3]{panangaden1999category}}}]
  \AP\label{def:panangaden}\phantomintro{Continuous maybe monad} Standard Borel spaces admit coproducts and thus a Maybe monad, $𝓜 ፡ \mathbf{BorelMes} → \mathbf{BorelMes}$. There exists a \kl{distributive law}, $(-)^{•}_X ፡ 𝓜𝓓X → 𝓓𝓜X$, yielding a monad structure on $(𝓓𝓜, \smash{μ_{•}^{𝓓𝓜}}, \smash{η^{𝓓𝓜}_{•}})$.
\end{definition}

\begin{definition}[Continuous normalization]
  \AP \intro[continuous normalization]{Normalization}, $\cnorm{(-)} ፡ 𝓓𝓜X → 𝓜𝓓X$, is the \kl[natural transformation]{natural transformation} defined by \emph{(i)} division with the mass of a \kl{subdistribution} whenever not zero, and \emph{(ii)} zero otherwise.
  \[\mbox{\emph{(i)}}
  \quad
  \cnorm{ν}(U) = \frac{ν(U)}{ν(X)}, 
  \mbox{ when } ν(X) ≠ 0, \mbox{ and }
  \quad
  \mbox{\emph{(ii)}}
  \quad \cnorm{ν}(U) = 0
  \mbox{ when } ν(X) = 0.\]
\end{definition}

\begin{theoremrep}
  \label{thm:continuous-normalization-is-a-sesquilaw} \kl[Continuous normalization]{Normalization} in Borel spaces, $\cnorm{(-)} ፡ 𝓓𝓜X → 𝓜𝓓X$,  induces a \kl{sesquilaw}.
\end{theoremrep}
\begin{proofsketch}
  Monoidality follows from linearity of the Lebesgue integral. The \kl{sesquilaw} axiom follows ultimately from the monotone convergence theorem, employed as for the associativity of the Giry monad \cite{panangaden1999category}. See the Appendix for details.
\end{proofsketch}
\begin{proof}
  During this proof, we use zero-handling division ($\oslash$), defined by $x \oslash 0 = 0$ and $x \oslash y = x / y$ whenever $y ≠ 0$.
  Let us first prove that \kl[continuous normalization]{normalization} is
  \kl[monoidal natural transformation]{monoidal}. Consider two measures, $ν_1 ∈
  𝓓X$ and $ν_2 ∈ 𝓓Y$. We use \emph{(i)} the definition of \kl{normalization}
  and the tensor of measures, \emph{(ii)} the definition of the indicator
  function, \emph{(iii)} linearity of the Lebesgue integral, \emph{(iv)}
  linearity of the Lebesgue integral, and \emph{(v)} the definition of the
  tensor of measures.
  \begin{align*}
    \cnorm{(ν_1 ⊗ ν_2)}(U)
    \overset{\emph{(i)}}{=} &\quad \bdiv{∫_{y ∈ Y} \left( ∫_{x ∈ X} \xi_U(x,y) · ν_1(\dif x) \right) · ν_2(\dif y) }{∫_{y ∈ Y} \left( ∫_{x ∈ X} \xi_{X × Y}(x,y) · ν_1(\dif x) \right) · ν_2(\dif y)} \\
    \overset{\emph{(ii)}}{=} &\quad  \bdiv{∫_{y ∈ Y} \left( ∫_{x ∈ X} \xi_U(x,y) · ν_1(\dif x) \right) · ν_2(\dif y) }{∫_{y ∈ Y} \left( ∫_{x ∈ X} 1 · ν_1(\dif x) \right) · ν_2(\dif y)} \\
    \overset{\emph{(iii)}}{=} & \quad  \bdiv{∫_{y ∈ Y} \left( ∫_{x ∈ X} \xi_U(x,y) · ν_1(\dif x) \right) · ν_2(\dif y) }{ν_1(X) · ν_2(Y)} \\
    \overset{\emph{(iv)}}{=} & \quad  ∫_{y ∈ Y} \left( ∫_{x ∈ X} \xi_U(x,y) · \bdiv{ν_1(\dif x)}{ν_1(X)} \right) · \bdiv{ν_2(\dif y)}{ν_2(Y)}
    \overset{\emph{(v)}}{=}  (\cnorm{(ν_1)} ⊗ \cnorm{(ν_2)})(U).
  \end{align*}
  
  Let us now prove that \kl[continuous normalization]{normalization}
  determines an \kl{almost distributive law}. We must prove that the following three diagrams commute.
  \begin{equation*}
  \begin{tikzcd}[ampersand replacement=\&,column sep=scriptsize]
  𝓓𝓜𝓜X \dar[swap]{\cnorm{} 𝓜} \rar{𝓓μ^{𝓜}} \& 𝓓𝓜X \rar{\cnorm{}} \& 𝓜𝓓 X \\
  𝓜𝓓𝓜 X \rar{𝓜 \cnorm{}} \& 𝓜𝓜𝓓 X \urar[swap]{μ^{𝓜} 𝓓} 
  \end{tikzcd}\quad
  \begin{tikzcd}[ampersand replacement=\&,column sep=scriptsize]
    𝓓 X \dar[swap]{𝓓 η^{𝓜}} \rar{η^{𝓜} 𝓓} \& 𝓜 𝓓 X \\
    𝓓𝓜 X \urar[swap]{\cnorm{}}
  \end{tikzcd}
  \quad
  \begin{tikzcd}[ampersand replacement=\&,column sep=scriptsize]
    𝓜 X \dar[swap]{η^{𝓓} 𝓜} \rar{𝓜 η^{𝓓}} \& 𝓜𝓓 X \\
    𝓓𝓜 X \urar[swap]{\cnorm{}}
  \end{tikzcd}
\end{equation*}

For $𝓜$\hyp{}multiplicativity, note that an element $d ∈ 𝓓𝓜𝓜 X$ is a measure
over $X + \{⊥\} + \{⊥'\}$. Let us reason by cases. If $d(\{⊥,⊥'\}) = 1$, then
$\cnorm{(μ^{𝓜}(d))} = ⊥ = μ^{𝓜}(\cnorm{(\cnorm{(d)})})$, either because $d(⊥') =
1$, or because $d(⊥') ≠ 1$ but then $d(⊥) > 0$ and $\cnorm{(d)}(⊥) = 1$. Assume,
thus, that $d(\{⊥,⊥'\}) ≠ 1$. We then need to prove that normalizing both
separately is the same as normalizing after identifying both.
  \begin{align*}
  μ^{𝓜}(\cnorm{(\cnorm{(d)})})(U) =
  \bdiv{
    \cnorm{(d)}(U)}
    {\cnorm{(d)}(X)}  =
  \bdiv{
    \bdiv{d(U)}
         {d(X + \{⊥\} )}}
    { \bdiv{d(X)}{d(X + \{⊥\} )} }
   = \bdiv{d(U)}{d(X)} 
   = \norm{(μ^{𝐌}(d))}.  
  \end{align*}

  For $𝓜$\hyp{}unitality, we must check that the normalization of a normalized
  distribution is itself: by $η^{𝓜}{(d)}(U) = d(U)$, we conclude that
  $η^{𝐌}{(d)}(X) = d(X) = 1$, and thus $\cnorm{(η^{𝓜}{(d)})} = η^{𝓜}{(d)}$.
  For $𝓓$\hyp{}unitality, we reason by cases on $𝓜X$. We check that
  $\cnorm{(η^{𝐃}{(⊥)})} = \cnorm{(\ket{⊥})} = ⊥$ and that
  $\cnorm{(η^{𝐃}{(x)})} = \cnorm{(\ket{x})} = x$.

  Let us now prove that \kl{normalization} and \kl{subdistributions} form a
  \kl{sesquilaw}. 
  
  We first need to prove that \kl{normalization} and the inclusion into
  subdistributions are a section-retraction pair. Given any normalized
  distribution, $ν ∈ 𝓜𝓓X$, the measure of the whole set must either be
  zero or one, $ν(X) = 1$ or $ν(X) = 0$. If it is zero, it follows that the
  measure of any subset, $U ⊆ X$, must also be zero, $ν(U) = 0$; as a consequence,
  in any of the two cases,
  \[
  \cnorm{(ν)}(U) = \bdiv{ν(U)}{ν(X)} = ν(U).
  \]

  Let us prove that the second formulation of the axiom of \kl{distributive
  sesquilaws} holds.
  \[\begin{tikzcd}[ampersand replacement=\&]
	𝓓𝓜𝓓 X \dar[swap]{𝓓 m} \rar{\cnorm{} 𝓓} \& 𝓜𝓓𝓓 X \rar{𝓜 μ} \& 𝓜𝓓 X \\
  𝓓𝓓𝓜 X \rar{μ 𝓜} \& 𝓓𝓜 X \urar[swap]{\cnorm{}} 
  \end{tikzcd}\]
  
  Let $ν ∈ 𝓓𝓜𝓓 X$ be a \kl{subdistribution} of \kl{distributions}. We must
  prove that normalizing and flattening the distributions is the same as, while
  regarding the distributions as subdistributions, flattening and then
  normalizing. In other words, we seek to prove
  \[
  𝓜 μ(\cnorm{(ν)}) = \cnorm{(μ(ν(-^{•})))}.
  \]
  On the left hand side, we use \emph{(i)} the definition of monad multiplication for the Giry monad,
  and \emph{(ii)} the definition of \kl{normalization}.
  \begin{align*}
    𝓜 μ(\cnorm{(ν(-^{•}))})(U) 
    \overset{\emph{(i)}}{=} \int_{α ∈ 𝓓 X} α(U) · \cnorm{(ν(-^{•}))}(\dif α)
    \overset{\emph{(ii)}}{=} \int_{α ∈ 𝓓 X} α(U) · \bdiv{ν(\dif α)}{ν(𝓓 X)}.
  \end{align*}
  On the right hand side, we use \emph{(i)} the definition of \kl{normalization}, \emph{(ii)} the
  definition of monad multiplication for the Giry monad and, finally, \emph{(iii)} that the $α ∈ 𝓓 X$
  elements are full distributions, meaning that $α(X) = 1$.
  \begin{align*}
    \cnorm{(μ(ν))}(U) 
    \overset{\emph{(i)}}{=} 
    \bdiv{μ(ν)(U)}{μ(ν)(X)} 
    \overset{\emph{(ii)}}{=} \bdiv{\int_{α ∈ 𝓓 X} α(U) · ν(\dif α)}{\int_{α ∈ 𝓓 X} α(X) · ν(\dif α)} 
    \overset{\emph{(iii)}}{=} \bdiv{\int_{α ∈ 𝓓 X} α(U) · ν(\dif α)}{ν(𝓓 X)}
  \end{align*}
  Lastly, we may divide by cases: whether $ν(𝓓X) = 0$ or $ν(𝓓X) ≠ 0$, we may use linearity of the Lebesgue integral to equate the two sides.
\end{proof}

\subsection{Sesquilaw homomorphisms}

All these notions of \kl{normalization}—discrete, possibilistic, and continuous—form \kl{sesquilaws}; and we can translate between these \kl{sesquilaws}. This section introduces \kl{sesquilaw morphisms} (\Cref{def:sesquilaw-morphism}), a formal notion for these translations. We start by recalling \kl{monad morphisms} (\Cref{def:monad-morphism}) and we present two examples of these translations at the end of the section.

It may more natural for this section to use double-categorical pasting diagrams, but we present it using standard commutative diagrams for consistency. %

\begin{definition}[Monad morphism]
  \label{def:monad-morphism}
  \AP A \intro{monad morphism} from a \kl{monad} $(𝐒,\smash{μ^{𝐒}},\smash{η^{𝐒}})$ on a category $𝔸₁$ to another \kl{monad} $(𝐓,\smash{μ^{𝐓}},\smash{η^{𝐓}})$ on a category $𝔸₂$ consists of a functor $𝐅 ፡ 𝔸₁ → 𝔸₂$ and a natural transformation $α ፡ 𝐅𝐒 → 𝐓𝐅$ such that the following two diagrams commute.
  \[
\begin{tikzcd}[ampersand replacement=\&]
	{𝐅𝐒𝐒} \& {𝐅𝐒} \& {𝐓𝐅} \\
	{𝐓𝐅𝐒} \& {𝐓𝐓𝐅}
	\arrow["{𝐅\mu^{𝐒}}", from=1-1, to=1-2]
	\arrow[""{name=0, anchor=center, inner sep=0}, "{\alpha 𝐒}"', from=1-1, to=2-1]
	\arrow["\alpha", from=1-2, to=1-3]
	\arrow["{𝐓\alpha}"', from=2-1, to=2-2]
	\arrow[""{name=1, anchor=center, inner sep=0}, "\alpha"', from=2-2, to=1-3]
	\arrow["{①}"{description}, draw=none, from=1, to=0]
\end{tikzcd}
\quad 
\begin{tikzcd}[ampersand replacement=\&]
	{𝐅} \& {𝐓𝐅} \\
	{𝐅𝐒}
	\arrow[""{name=0, anchor=center, inner sep=0}, "{\eta^{𝐓}𝐅}", from=1-1, to=1-2]
	\arrow["{𝐅\eta^{𝐒}}"', from=1-1, to=2-1]
	\arrow["\alpha"', from=2-1, to=1-2]
	\arrow["{②}"{description}, draw=none, from=2-1, to=0]
\end{tikzcd}
\]
\end{definition}

\begin{definition}[Sesquilaw morphism]
  \label{def:sesquilaw-morphism}
  \AP A \intro{sesquilaw morphism} from a \kl{sesquilaw} $(𝐒₁,𝐓₁,m₁,n₁)$ on a category $𝔸₁$ to another \kl{sesquilaw} $(𝐒₂,𝐓₂,m₂,n₂)$ on a category $𝔸₂$ is a pair of \kl{monad morphisms} $α ፡ 𝐅 𝐒₁ → 𝐒₂ 𝐅$ and $β ፡ 𝐅 𝐓₁ → 𝐓₂ 𝐅$ such that the following two diagrams commute. 
\[
\begin{tikzcd}[cramped]
	{𝐅𝐒_1𝐓_1} & {𝐅𝐓_1𝐒_1} & {𝐓_2𝐅𝐒_1} \\
	{𝐒_2𝐅𝐓_1} & {𝐒_2𝐓_2𝐅} & {𝐓_2𝐒_2𝐅}
	\arrow["{𝐅m_1}", from=1-1, to=1-2]
	\arrow[""{name=0, anchor=center, inner sep=0}, "{\alpha𝐓_1}"', from=1-1, to=2-1]
	\arrow["{\beta𝐒_1}", from=1-2, to=1-3]
	\arrow[""{name=1, anchor=center, inner sep=0}, "{𝐓_2\alpha}", from=1-3, to=2-3]
	\arrow["{𝐒_2\beta}"', from=2-1, to=2-2]
	\arrow["{m_2𝐅}"', from=2-2, to=2-3]
	\arrow["{①}"{description}, draw=none, from=1, to=0]
\end{tikzcd}
\quad
\begin{tikzcd}[cramped]
	{𝐅𝐓_1𝐒_1} & {𝐅𝐒_1𝐓_1} & {𝐒_2𝐅𝐓_1} \\
	{𝐓_2𝐅𝐒_1} & {𝐓_2𝐒_2𝐅} & {𝐒_2𝐓_2𝐅}
	\arrow["{𝐅n_1}", from=1-1, to=1-2]
	\arrow[""{name=0, anchor=center, inner sep=0}, "{\beta𝐒_1}"', from=1-1, to=2-1]
	\arrow["{\alpha𝐓_1}", from=1-2, to=1-3]
	\arrow[""{name=1, anchor=center, inner sep=0}, "{𝐒_2\beta}", from=1-3, to=2-3]
	\arrow["{𝐓_2\alpha}"', from=2-1, to=2-2]
	\arrow["{n_2𝐅}"', from=2-2, to=2-3]
	\arrow["{②}"{description}, draw=none, from=1, to=0]
\end{tikzcd}
\]
\end{definition}

\begin{example}[Support]
  \AP The \intro{support} of a \kl{distribution}, is its subset of elements. 
  \kl{Support} determines a \kl{natural transformation}, $\supp_X ፡ 𝐃 X → 𝐑 X$, defined by
  \[
  \supp(λ₁\!\ketb{x₁} + \dots + λ_n\!\ketb{x_n}) = \{ \cBlu{x_1}, ..., \cBlu{x_n} \}.
  \]
  It is well-known that \kl{support} is a \kl{monad morphism}. We note that it is moreover a \kl{sesquilaw morphism}: intuitively, the \kl{support} of a \kl{normalization} is the \kl{post-selection} of its \kl{support}.
\end{example}

\begin{propositionrep}[Support]
  \label{prop:support-is-a-sesquilaw-homomorphism}
  \AP \kl{Support}, $\supp_X ፡ 𝐃 X → 𝐑 X$, is a \kl{sesquilaw homomorphism}—with $𝐅 = \mathrm{Id}$ the identity functor—between \kl{normalization}  and \kl{post-selection}.
\end{propositionrep}

\begin{example}[Discrete to continuous]
  Discrete \kl{normalization} is a particular case of continuous normalization. Every set can be seen as a discrete standard Borel space, $𝐅 ፡ \mathbf{Set} → \mathbf{BorelMes}$; and there exist monad morphisms $α ፡ 𝐅𝐌 X → 𝓜𝐅X$ and $β ፡ 𝐅𝐃X → 𝓓𝐅X$ regarding each finitary distribution as a probability measure.
\end{example}

\begin{proposition}[Finitary distributions as probability measures]
  \label{prop:finitary-continuous}
  \AP Inclusion of finitary distributions into probability measures, $β ፡ 𝐅𝐃X → 𝓓𝐅X$, and the inclusion of the maybe monad, $α ፡ 𝐅𝐌X → 𝓜𝐅X$,  form a sesquilaw—with $𝐅$ the inclusion of sets into standard Borel spaces—between discrete and continuous normalization.
\end{proposition} %
\section{Magmads and magmoids}
\label{sec:magmads-magmoids}

We now introduce \kl{magmads} (\Cref{sec:magmads})—\emph{non-associative ``monads''}—and develop our main example: the \kl[normalized distribution magmad]{magmad of normalized distributions}. Each \kl{magmad} induces a Kleisli \kl{magmoid} (\Cref{sec:kleislimagmoids})—a \emph{non-associative ``category''}—and enables two \emph{magmadic metalanguages}: one associating to the right, and one associating to the left (\Cref{sec:magmoidal-metalanguage}).
We use these metalanguages to characterize Bayes' update, Pearl's update, and Jeffrey's update \Cref{prop:bayes-removal,prop:pearl-jeffrey}.

\subsection{Magmads}
\label{sec:magmads}

A \kl{magmad} is a \kl{monad} failing the associativity axiom\footnote{For this text, we convene that magmas—and any of their derived notions—are unital.} (\Cref{def:magmad}). \kl{Magmads} are infrequent in the literature. Particularly in functional programming, it is assumed that a \kl{monad} needs to obey all of its axioms to be of interest. Still, \kl{almost distributive laws} only induce \kl{magmads} (\Cref{prop:almostd-magmad})—and \kl{normalized distributions} only form a \kl{magmad}.

When the \kl{magmad} arises from a \kl{sesquilaw}—as in the case of \kl{normalized distributions}—it is moreover endowed with an action from the \kl{sesquilaw} \kl{monad} (\Cref{lemma:right-monad-action}).

\begin{definition}[Magmad]
  \label[definition]{def:magmad}
  \AP A \intro{magmad}, $(𝐓,μ^{𝐓},η^{𝐓})$, consists of an endofunctor, $𝐓 ፡ ℂ → ℂ$, together with natural transformations $η_X ፡ X → 𝐓 X$ and $μ_X ፡ 𝐓𝐓 X → 𝐓 X$ making the following two diagrams commute.
  \begin{equation*}
\begin{tikzcd}
	{𝐓X} & {𝐓X} \\
	{𝐓𝐓X}
	\arrow[from=1-1, to=1-2]
	\arrow["{\eta^{𝐓}𝐓}"', from=1-1, to=2-1]
	\arrow[""{name=0, anchor=center, inner sep=0}, "{\mu^{𝐓}}"', from=2-1, to=1-2]
	\arrow["{①}"{description}, draw=none, from=1-1, to=0]
\end{tikzcd}
\qquad
\begin{tikzcd}[ampersand replacement=\&]
	{𝐓X} \& {𝐓X} \\
	{𝐓𝐓X}
	\arrow[""{name=0, anchor=center, inner sep=0}, from=1-1, to=1-2]
	\arrow["{𝐓\eta^{𝐓}}"', from=1-1, to=2-1]
	\arrow["{\mu^{𝐓}}"', from=2-1, to=1-2]
	\arrow["{②}"{description}, draw=none, from=0, to=2-1]
\end{tikzcd}
\end{equation*}
\end{definition}

\crefalias{propositionrep}{proposition}
\begin{propositionrep}
  \label{prop:almostd-magmad}
  \AP Any \kl{almost-distributive law}, $\smash{ψ} ፡ 𝐓𝐒 X → 𝐒𝐓 X$, induces a \kl{magmad}, $(𝐒𝐓, \smash{μ^{𝐒𝐓}_ψ}, \smash{η^{𝐒𝐓}_ψ})$ defined by $\smash{μ^{𝐒𝐓}_ψ} = 𝐒 ψ 𝐓 ⨾ \smash{μ^{𝐒}μ^{𝐓}}$ and $\smash{η^{𝐒𝐓}_ψ} = \smash{η^{𝐒}η^{𝐓}}$.
\end{propositionrep}
\begin{proof}
  The unit of the composite \kl{non-associative monad} is $η = η^{𝐒}η^{𝐓}$; the multiplication is $μ = 𝐒 ψ 𝐓 𑊩 μ^{𝐒} μ^{𝐓}$. From unitality of the \kl{almost distributive law}, unitality of the \kl{non-associative monad} follows. Whenever the \kl{almost distributive law} is moreover monoidal, the unit and multiplication of the monad become monoidal by construction.
\end{proof}

\begin{theoremrep}
  \label{lemma:right-monad-action}
  Any \kl{sesquilaw}, $(𝐒,𝐓,m,n)$, induces a right action of its \kl{monad}, $𝐓𝐒$, into its \kl{magmad}, $𝐒𝐓$: i.e., a natural transformation,
  \(u_X ፡ 𝐒𝐓𝐓𝐒X → 𝐒𝐓X,\)
  making the following two diagrams commute.
  \[\begin{tikzcd}[ampersand replacement=\&]
    𝐒𝐓X \rar{η^{𝐓𝐒}}\drar[swap]{\id} \& \dar{u} 𝐒𝐓𝐓𝐒X \& 𝐒𝐓𝐓𝐒𝐓𝐒X \rar{𝐒𝐓μ^{𝐓𝐒}} \dar[swap]{u𝐓𝐒} \& 𝐒𝐓𝐓𝐒X \dar{u} \\
    \& 𝐒𝐓X \& 𝐒𝐓𝐓𝐒X \rar{u} \& 𝐒𝐓X  
  \end{tikzcd}\]
  This action is defined by either side of the following commutative diagram.
  \[\begin{tikzcd}[ampersand replacement=\&]
    𝐒𝐓𝐓𝐒X \dar[swap]{m𝐓𝐒} \rar{𝐒 μ^{𝐓} 𝐒} \&  𝐒𝐓𝐒𝐒X \rar{m𝐒} \& 𝐓𝐒𝐒X \dar{𝐓μ^{𝐒}} \\
    𝐓𝐒𝐓𝐒X \rar{𝐓m𝐒} \& 𝐓𝐓𝐒𝐒X \rar{μ^{𝐓}μ^{𝐒}} \& 𝐓𝐒X \rar{n} \& 𝐒𝐓X
  \end{tikzcd}\]
\end{theoremrep}
\begin{proofsketch}
  Multiplicativity follows from the two axioms of \kl{sesquilaws}, together  with $𝐒$-multiplicativity of the \kl{almost distributive law} and the $𝐓$-multiplicativity of the \kl{distributive law}; we also use the associativity of both \kl{monads}. Unitality follows from the unitality of both monads and the second axiom of \kl{sesquilaws}.
\end{proofsketch}
\begin{proof}
  Let us use commutative diagrams for this proof. We will show that any \kl{sesquilaw}, $(m,n,𝐒,𝐓)$, induces an action of th monad $𝐓𝐒$ into the almost monad $𝐒𝐓$. We define the action, $u ፡ 𝐒𝐓𝐓𝐒 → 𝐒𝐓$, by string diagrams, as any of the two following equivalent definitions.
  \[\begin{tikzcd}[ampersand replacement=\&]
      𝐒𝐓𝐓𝐒X \dar[swap]{m𝐓𝐒} \rar{𝐒 μ^{𝐓} 𝐒} \&  𝐒𝐓𝐒𝐒X \rar{m𝐒} \& 𝐓𝐒𝐒X \dar{𝐓μ^{𝐒}} \\
      𝐓𝐒𝐓𝐒X \rar{𝐓m𝐒} \& 𝐓𝐓𝐒𝐒X \rar{μ^{𝐓}μ^{𝐒}} \& 𝐓𝐒X \rar{n} \& 𝐒𝐓X
    \end{tikzcd}\]
Let us prove multiplicativity (in Figure \ref{fig:action-multiplicativity}). We use ⓪ naturality, ① $𝐒$-multiplicativity of the \kl{almost distributive law} $n$, ② the \kl{sesquilaw} axiom, ③ that \kl{sesquilaws} are inverses, ④ $𝐒$-multiplicativity of the \kl{distributive law} $m$, ⑤ associativity of $𝐒$, ⑥ $𝐓$-multiplicativity of the \kl{distributive law} $m$, and ⑦ associativity of the monad $𝐓$.
\begin{figure}[ht]
\[\begin{tikzcd}[ampersand replacement=\&,column sep=scriptsize]
	\&\&\&\& {𝐓𝐒𝐒X} \& {𝐓𝐒X} \& {𝐒𝐓X} \\
	{𝐓𝐒𝐓𝐒X} \& {𝐒𝐓𝐓𝐒X} \& {𝐒𝐓𝐒X} \& {𝐓𝐒𝐒X} \& {𝐒𝐓𝐒X} \& {𝐒𝐒𝐓X} \& {𝐒𝐓X} \\
	{𝐓𝐒𝐓𝐒X} \& {𝐓𝐓𝐒𝐒X} \& {𝐓𝐒𝐒X} \& {𝐒𝐓𝐒X} \& {𝐒𝐒𝐓X} \& {𝐒𝐓X} \\
	{𝐓𝐒𝐒𝐓𝐒X} \& {𝐓𝐒𝐓𝐒X} \& {𝐓𝐓𝐒𝐒X} \& {𝐓𝐒𝐒X} \& {𝐓𝐒X} \& {𝐒𝐓X} \\
	{𝐒𝐓𝐒𝐓𝐒X} \& {𝐓𝐒𝐒𝐓𝐒X} \& {𝐓𝐒𝐓𝐒𝐒X} \& {𝐓𝐓𝐒𝐒𝐒X} \& {𝐓𝐓𝐒𝐒X} \& {𝐓𝐓𝐒X} \& {𝐓𝐒X} \\
	{𝐒𝐓𝐒𝐓𝐒X} \& {𝐒𝐓𝐓𝐒𝐒X} \& {𝐓𝐒𝐓𝐒𝐒X} \& {𝐓𝐓𝐒𝐒𝐒X} \& {𝐓𝐓𝐒𝐒X} \& {𝐓𝐓𝐒X} \& {𝐓𝐒X} \\
	{𝐒𝐓𝐒𝐓𝐒X} \& {𝐒𝐓𝐓𝐒𝐒X} \& {𝐒𝐓𝐓𝐒X} \& {𝐓𝐒𝐓𝐒X} \& {𝐓𝐓𝐒𝐒X} \& {𝐓𝐒𝐒X} \& {𝐓𝐒X} \\
	{𝐒𝐓𝐓𝐒𝐓𝐒X} \& {𝐒𝐓𝐒𝐓𝐒X} \& {𝐒𝐓𝐓𝐒𝐒X} \& {𝐒𝐓𝐓𝐒X} \& {𝐒𝐓𝐒X} \& {𝐓𝐒𝐒X} \\
	{𝐒𝐓𝐓𝐒𝐓𝐒X} \& {𝐒𝐓𝐓𝐓𝐒𝐒X} \& {𝐒𝐓𝐓𝐒𝐒X} \& {𝐒𝐓𝐒𝐒X} \& {𝐒𝐓𝐒X} \\
	\& {𝐒𝐓𝐓𝐓𝐒𝐒X} \& {𝐒𝐓𝐓𝐒𝐒X} \& {𝐒𝐓𝐒𝐒X} \& {𝐒𝐓𝐒X} \\
	\& {𝐒𝐓𝐓𝐓𝐒𝐒X} \& {𝐒𝐓𝐓𝐒𝐒X} \& {𝐒𝐓𝐓𝐒X} \& {𝐒𝐓𝐒X} \\
	\& {𝐒𝐓𝐓𝐓𝐒𝐒X} \& {𝐒𝐓𝐓𝐒X}
	\arrow["{𝐓\mu^{𝐒}}", from=1-5, to=1-6]
	\arrow[""{name=0, anchor=center, inner sep=0}, no head, from=1-5, to=2-4]
	\arrow["n", from=1-6, to=1-7]
	\arrow[""{name=1, anchor=center, inner sep=0}, no head, from=1-7, to=2-7]
	\arrow["{n𝐓𝐒}", from=2-1, to=2-2]
	\arrow[""{name=2, anchor=center, inner sep=0}, no head, from=2-1, to=3-1]
	\arrow["{𝐒\mu^{𝐓}𝐒}", from=2-2, to=2-3]
	\arrow["{m𝐒}", from=2-3, to=2-4]
	\arrow[""{name=3, anchor=center, inner sep=0}, no head, from=2-3, to=3-4]
	\arrow["{n𝐒}", from=2-4, to=2-5]
	\arrow["{𝐒n}", from=2-5, to=2-6]
	\arrow[""{name=4, anchor=center, inner sep=0}, no head, from=2-5, to=3-4]
	\arrow["{\mu^{𝐒}𝐓}", from=2-6, to=2-7]
	\arrow[no head, from=2-7, to=3-6]
	\arrow["{𝐓m𝐒}", from=3-1, to=3-2]
	\arrow[no head, from=3-1, to=4-2]
	\arrow["{\mu^{𝐓}𝐒𝐒}", from=3-2, to=3-3]
	\arrow[no head, from=3-2, to=4-3]
	\arrow["{n𝐒}", from=3-3, to=3-4]
	\arrow[""{name=5, anchor=center, inner sep=0}, no head, from=3-3, to=4-4]
	\arrow["{𝐒n}", from=3-4, to=3-5]
	\arrow["{\mu^{𝐒}𝐓}", from=3-5, to=3-6]
	\arrow[""{name=6, anchor=center, inner sep=0}, no head, from=3-6, to=4-6]
	\arrow["{𝐓\mu^{𝐒}𝐓𝐒}", from=4-1, to=4-2]
	\arrow[""{name=7, anchor=center, inner sep=0}, no head, from=4-1, to=5-2]
	\arrow["{𝐓m𝐒}", from=4-2, to=4-3]
	\arrow["{\mu^{𝐓}𝐒𝐒}", from=4-3, to=4-4]
	\arrow[""{name=8, anchor=center, inner sep=0}, no head, from=4-3, to=5-5]
	\arrow["{𝐓\mu^{𝐒}}", from=4-4, to=4-5]
	\arrow["n", from=4-5, to=4-6]
	\arrow[""{name=9, anchor=center, inner sep=0}, no head, from=4-5, to=5-7]
	\arrow["{m𝐒𝐓𝐒}"', from=5-1, to=5-2]
	\arrow[""{name=10, anchor=center, inner sep=0}, no head, from=5-1, to=6-1]
	\arrow["{𝐓𝐒m𝐒}"', from=5-2, to=5-3]
	\arrow["{𝐓m𝐒𝐒}"', from=5-3, to=5-4]
	\arrow[""{name=11, anchor=center, inner sep=0}, no head, from=5-3, to=6-3]
	\arrow["{𝐓𝐓\mu^{𝐒}𝐒}"', from=5-4, to=5-5]
	\arrow[""{name=12, anchor=center, inner sep=0}, no head, from=5-4, to=6-4]
	\arrow["{𝐓𝐓\mu^{𝐒}}"', from=5-5, to=5-6]
	\arrow["{\mu^{𝐓}𝐒}"', from=5-6, to=5-7]
	\arrow[""{name=13, anchor=center, inner sep=0}, no head, from=5-6, to=6-6]
	\arrow[no head, from=5-7, to=6-7]
	\arrow["{𝐒𝐓m𝐒}", from=6-1, to=6-2]
	\arrow[no head, from=6-1, to=7-1]
	\arrow["{m𝐓𝐒𝐒}", from=6-2, to=6-3]
	\arrow[""{name=14, anchor=center, inner sep=0}, no head, from=6-2, to=7-2]
	\arrow["{𝐓m𝐒𝐒}", from=6-3, to=6-4]
	\arrow["{𝐓𝐓𝐒\mu^{𝐒}}", from=6-4, to=6-5]
	\arrow["{𝐓𝐓\mu^{𝐒}}", from=6-5, to=6-6]
	\arrow[""{name=15, anchor=center, inner sep=0}, no head, from=6-5, to=7-5]
	\arrow["{\mu^{𝐓}𝐒}", from=6-6, to=6-7]
	\arrow[""{name=16, anchor=center, inner sep=0}, no head, from=6-7, to=7-7]
	\arrow["{𝐒𝐓m𝐒}", from=7-1, to=7-2]
	\arrow[no head, from=7-1, to=8-2]
	\arrow["{𝐒𝐓𝐓\mu^{𝐒}}", from=7-2, to=7-3]
	\arrow[no head, from=7-2, to=8-3]
	\arrow["{m𝐓𝐒}", from=7-3, to=7-4]
	\arrow[""{name=17, anchor=center, inner sep=0}, no head, from=7-3, to=8-4]
	\arrow["{𝐓m𝐒}", from=7-4, to=7-5]
	\arrow["{\mu^{𝐓}𝐒𝐒}", from=7-5, to=7-6]
	\arrow["{𝐓\mu^{𝐒}}", from=7-6, to=7-7]
	\arrow[""{name=18, anchor=center, inner sep=0}, no head, from=7-6, to=8-6]
	\arrow["{𝐒\mu^{𝐓}𝐒𝐓𝐒}", from=8-1, to=8-2]
	\arrow[""{name=19, anchor=center, inner sep=0}, no head, from=8-1, to=9-1]
	\arrow["{𝐒𝐓m𝐒}", from=8-2, to=8-3]
	\arrow["{𝐒𝐓𝐓\mu^{𝐒}}", from=8-3, to=8-4]
	\arrow[""{name=20, anchor=center, inner sep=0}, no head, from=8-3, to=9-3]
	\arrow["{𝐒\mu^{𝐓}𝐒}", from=8-4, to=8-5]
	\arrow["{m𝐒}", from=8-5, to=8-6]
	\arrow[""{name=21, anchor=center, inner sep=0}, no head, from=8-5, to=9-5]
	\arrow["{𝐒𝐓𝐓m𝐒}", from=9-1, to=9-2]
	\arrow["{𝐒\mu^{𝐓}𝐓𝐒𝐒}", from=9-2, to=9-3]
	\arrow[""{name=22, anchor=center, inner sep=0}, no head, from=9-2, to=10-2]
	\arrow["{𝐒\mu^{𝐓}𝐒𝐒}", from=9-3, to=9-4]
	\arrow["{𝐒𝐓\mu^{𝐒}}", from=9-4, to=9-5]
	\arrow[""{name=23, anchor=center, inner sep=0}, no head, from=9-4, to=10-4]
	\arrow[no head, from=9-5, to=10-5]
	\arrow["{𝐒𝐓\mu^{𝐓}𝐒𝐒}", from=10-2, to=10-3]
	\arrow[no head, from=10-2, to=11-2]
	\arrow["{𝐒\mu^{𝐓}𝐒𝐒}", from=10-3, to=10-4]
	\arrow[""{name=24, anchor=center, inner sep=0}, no head, from=10-3, to=11-3]
	\arrow["{𝐒𝐓\mu^{𝐒}}", from=10-4, to=10-5]
	\arrow[""{name=25, anchor=center, inner sep=0}, no head, from=10-5, to=11-5]
	\arrow["{𝐒𝐓\mu^{𝐓}𝐒𝐒}", from=11-2, to=11-3]
	\arrow[""{name=26, anchor=center, inner sep=0}, no head, from=11-2, to=12-2]
	\arrow["{𝐒𝐓𝐓\mu^{𝐒}}", from=11-3, to=11-4]
	\arrow["{𝐒\mu^{𝐓}𝐒}", from=11-4, to=11-5]
	\arrow[""{name=27, anchor=center, inner sep=0}, no head, from=11-4, to=12-3]
	\arrow["{𝐒𝐓\mu^{𝐒}\mu^{𝐓}}", from=12-2, to=12-3]
	\arrow["{①}"{description}, draw=none, from=0, to=1]
	\arrow["{②}"{description}, draw=none, from=2, to=3]
	\arrow["{③}"{description}, draw=none, from=3, to=4]
	\arrow["{①}"{description}, draw=none, from=5, to=6]
	\arrow["{④}"{description}, draw=none, from=7, to=8]
	\arrow["{⓪}"{description}, draw=none, from=8, to=9]
	\arrow["{⓪}"{description}, draw=none, from=10, to=11]
	\arrow["{⑤}"{description}, draw=none, from=12, to=13]
	\arrow["{④}"{description}, draw=none, from=14, to=15]
	\arrow["{⓪}"{description}, draw=none, from=16, to=15]
	\arrow["{⑥}"{description}, draw=none, from=17, to=18]
	\arrow["{⓪}"{description}, draw=none, from=19, to=20]
	\arrow["{⓪}"{description}, draw=none, from=20, to=21]
	\arrow["{⑦}"{description}, draw=none, from=22, to=23]
	\arrow["{⓪}"{description}, draw=none, from=24, to=25]
	\arrow["{⓪}"{description}, draw=none, from=26, to=27]
\end{tikzcd}\]   %
  \caption{Multiplicativity for the right action of a sesquilaw.}
  \label{fig:action-multiplicativity}  
\end{figure}

The proof of unitality is more direct: it uses ① naturality, ② unitality of  the \kl{monad} $𝐓$, ③ unitality of the \kl{monad} $𝐒$, and ④ that \kl{sesquilaws} are inverses.

\[\begin{tikzcd}[ampersand replacement=\&]
	\& {𝐒𝐓X} \& {𝐒𝐓𝐓𝐒X} \& {𝐒𝐓𝐒X} \&\&\& \\
	{𝐒𝐓X} \& {𝐒𝐓𝐒X} \& {𝐒𝐓𝐓𝐒X} \& {𝐒𝐓𝐒X} \& {𝐓𝐒𝐒X} \& {𝐓𝐒X} \& {𝐒𝐓X} \\
	\& {𝐒𝐓X} \& {𝐒𝐓𝐒X} \& {𝐓𝐒𝐒X} \& {𝐓𝐒X} \& {𝐒𝐓X} \\
	\& {𝐒𝐓X} \& {𝐓𝐒X} \& {𝐓𝐒𝐒X} \& {𝐓𝐒X} \& {𝐒𝐓X} \\
	\&\& {𝐒𝐓X} \& {𝐓𝐒X} \& {𝐒𝐓X} \\
	\&\&\& {𝐒𝐓X}
	\arrow["{𝐒𝐓\eta^{𝐒}\eta^{𝐓}}", from=1-2, to=1-3]
	\arrow[""{name=0, anchor=center, inner sep=0}, no head, from=1-2, to=2-1]
	\arrow["{𝐒\mu^{𝐓}𝐒}", from=1-3, to=1-4]
	\arrow[""{name=1, anchor=center, inner sep=0}, no head, from=1-3, to=2-3]
	\arrow[no head, from=1-4, to=2-4]
	\arrow["{𝐒𝐓\eta^{𝐒}}", from=2-1, to=2-2]
	\arrow[no head, from=2-1, to=3-2]
	\arrow["{𝐒𝐓\eta^{𝐓}𝐒}", from=2-2, to=2-3]
	\arrow[""{name=2, anchor=center, inner sep=0}, no head, from=2-2, to=3-3]
	\arrow["{𝐒\mu^{𝐓}𝐒}", from=2-3, to=2-4]
	\arrow["{m𝐒}", from=2-4, to=2-5]
	\arrow[""{name=3, anchor=center, inner sep=0}, no head, from=2-4, to=3-3]
	\arrow["{𝐓\mu^{𝐒}}", from=2-5, to=2-6]
	\arrow["n", from=2-6, to=2-7]
	\arrow[no head, from=2-6, to=3-5]
	\arrow[no head, from=2-7, to=3-6]
	\arrow["{𝐒𝐓\eta^{𝐒}}", from=3-2, to=3-3]
	\arrow[""{name=4, anchor=center, inner sep=0}, no head, from=3-2, to=4-2]
	\arrow["{m𝐒}", from=3-3, to=3-4]
	\arrow[no head, from=3-4, to=2-5]
	\arrow["{𝐓\mu^{𝐒}}", from=3-4, to=3-5]
	\arrow[""{name=5, anchor=center, inner sep=0}, no head, from=3-4, to=4-4]
	\arrow["n", from=3-5, to=3-6]
	\arrow[no head, from=3-5, to=4-5]
	\arrow[no head, from=3-6, to=4-6]
	\arrow["m", from=4-2, to=4-3]
	\arrow[no head, from=4-2, to=5-3]
	\arrow["{𝐓𝐒\eta^{𝐒}}", from=4-3, to=4-4]
	\arrow[""{name=6, anchor=center, inner sep=0}, no head, from=4-3, to=5-4]
	\arrow["{𝐓\mu^{𝐒}}", from=4-4, to=4-5]
	\arrow["n", from=4-5, to=4-6]
	\arrow[""{name=7, anchor=center, inner sep=0}, no head, from=4-5, to=5-4]
	\arrow[no head, from=4-6, to=5-5]
	\arrow["m", from=5-3, to=5-4]
	\arrow[""{name=8, anchor=center, inner sep=0}, no head, from=5-3, to=6-4]
	\arrow["n", from=5-4, to=5-5]
	\arrow[""{name=9, anchor=center, inner sep=0}, no head, from=5-5, to=6-4]
	\arrow["{①}"{description}, draw=none, from=0, to=1]
	\arrow["{②}"{description}, draw=none, from=2, to=3]
	\arrow["{①}"{description}, draw=none, from=4, to=5]
	\arrow["{③}"{description}, draw=none, from=6, to=7]
	\arrow["{④}"{description}, draw=none, from=8, to=9]
\end{tikzcd}\]

With these two axioms, we have built an action. In the monoidal case, we repeat the exact same proof just considering that all natural transformations are monoidal natural transformations.
\end{proof}

\begin{corollary}
  \AP \phantomintro{normalized distribution magmad}
  \kl{Normalized distributions} form a \kl{magmad}, $(𝐌𝐃, \smash{μ^{𝐌𝐃}_{∘}}, \smash{η^{𝐌𝐃}_{∘}})$. There exists a right action of the \kl{subdistribution monad}, $(𝐃𝐌,\smash{μ^{𝐃𝐌}_{•}}, \smash{η^{𝐃𝐌}_{•}})$, on the \kl{normalized distribution} \kl{magmad}.
\end{corollary}

\subsection{Kleisli Magmoids}
\label{sec:kleislimagmoids}

\kl{Magmoids} are \emph{non-associative ``categories''}. Defining them is relatively straightforward: we repeat the definition of category without associativity (\Cref{def:magmoid}, c.f.~\cite{mangel2025classical}).  Of course, every category is a \kl{magmoid}, but we focus on the non-associative examples.

In the same sense that any \kl{monad} induces a Kleisli category, any \kl{magmad} induces a Kleisli \kl{magmoid} (\Cref{prop:kleisli-magmoids}).  Our main non-associative example is the \kl{magmoid} of \kl{normalized stochastic kernels} (\Cref{prop:magmoid-normalized-stochastic-kernels}): the Kleisli \kl{magmoid} of the \kl{normalized distribution magmad}.

\begin{definition}[Magmoid]
  \label[definition]{def:magmoid}
  \AP A \intro{magmoid}, $𝕄$, consists of morphisms, $𝕄(X;Y)$, between some collection of objects $X, Y ∈ 𝕄_{obj}$, equipped with a composition operation, $(⨾) ፡ 𝕄(X;Y) × 𝕄(Y;Z) → 𝕄(X;Z)$, and identities, $\id ∈ 𝕄(X;X)$, that are unital, meaning that $f ⨾ \id = f = \id ⨾ f$.
\end{definition}

\begin{propositionrep}[Magmoid of normalized stochastic kernels]
  \label[proposition]{prop:magmoid-normalized-stochastic-kernels}
  \AP \intro{Normalized kernels} between sets, $X → 𝐌𝐃Y$, form a \kl{magmoid}, $\Norm$. Given $f ፡ X → 𝐌𝐃Y$, we write $f(x;y)$ for the coefficient of $y ∈ Y$ in the \kl{normalized distribution} $f(x)$. Composition of two morphisms, $f ፡ X → 𝐌𝐃 Y$ and $g ፡ Y → 𝐌𝐃 Z$, is defined by
  \[
  (f ⨾ g)(x; z) =
    \frac{\sum_{v ∈ Y} f(x;v) · g(v;z)}{\sum_{v ∈ Y} \sum_{w ∈ Z} f(x;v) · g(v;w)}, \mbox{ only when } \sum_{v, w} f(x;v) · g(v;w) ≠ 0.
  \]
  
\end{propositionrep}
\begin{proof}
  We define $\id(x;x') = [x = x']$. Let us check that it is right unital. Let $f ፡ X → 𝐌𝐃Y$ be a \kl{normalized kernel}. 
  \begin{align*}
    (f ⨾ \id)(x;y) 
    = \bdiv{\sum_{y'} f(x;y') · [y' = y]}{\sum_{y',y} f(x;y') · [y' = y]}
    = \bdiv{f(x;y)}{\sum_{y} f(x;y)} = f(x;y).
  \end{align*}
  On the last step, either $\sum_{y'} f(x;y') = 0$, which implies $f(x;y) = 0$, or $\sum_{y'} f(x;y') = 1$, which simplifies the division. Left unitality is analogous.
\end{proof}

\begin{remark}
	In other words, if we consider the associated \kl{substochastic kernels}, $f^{•} ፡ X → 𝐃𝐌 Y$ and $g^{•} ፡ Y → 𝐃𝐌 Z$, it is the \kl{normalization} of their composition as subdistributions, $f ⨾ g = {(f^{•} 𑊩 g^{•})}^{∘}$.
\end{remark}

\begin{propositionrep}[Kleisli magmoids]
  \label[proposition]{prop:kleisli-magmoids}
  \AP Any \kl{magmad}, $(𝗡, \smash{μ^{𝗡}}, \smash{η^{𝗡}})$ over a category \ $ℂ$ with composition operation $(𑊩)$, induces a \kl{magmoid}, $𝕂(𝗡)$, with its same objects, with morphisms $𝕂(𝗡)(X;Y) = ℂ(X; 𝗡 Y)$, with identities $\id = η^{𝗡}$, and with compositions, $(f ⨾ g) = f 𑊩 𝗡 g 𑊩 μ^{𝗡}$.
\end{propositionrep}

\begin{propositionrep}[Sesquilaw magmoids]
  \label{lemma:relevance}
  \AP Any \kl{sesquilaw}, $(𝐒,𝐓,m,n)$, induces a \kl{Kleisli magmoid} where associativity is satisfied, $f ⨾ (g ⨾ h) = (f ⨾ g) ⨾ h$, if either
  \begin{enumerate}
    \item $f = u 𑊩 𝐒 η^{𝐓}$ for some $u ፡ X → 𝐒 Y$;
    \item $g = v 𑊩 𝐒 η^{𝐓}$ for some $v ፡ Y → 𝐒 Z$; or
    \item $h = w 𑊩 η^{𝐒} 𝐓$ for some $w ፡ Z → 𝐓 W$.
  \end{enumerate}
\end{propositionrep}
\begin{proof}
  Proving these equations for the case of normalized kernels is straightforward—we only need to check that the resulting expressions are equal—but difficult to repeat in the general case. The general case requires diagram chasing.
  These three cases reduce, respectively, to the commutativity of the following three diagrams.
\[\begin{tikzcd}[ampersand replacement=\&]
	{𝐒𝐒𝐓𝐒𝐓X} \& {𝐒𝐓𝐒𝐓𝐒𝐓X} \& {𝐒𝐓𝐒𝐓X} \\
	{𝐒𝐓𝐒𝐓𝐒𝐓X} \& {𝐒𝐓𝐒𝐓X} \& {𝐒𝐓X}
	\arrow["{𝐒\eta^{𝐓}𝐒𝐓𝐒𝐓}", from=1-1, to=1-2]
	\arrow["{𝐒\eta^{𝐓}𝐒𝐓𝐒𝐓}"', from=1-1, to=2-1]
	\arrow["{\mu^{𝐒𝐓}𝐒𝐓}", from=1-2, to=1-3]
	\arrow["{\mu^{𝐒𝐓}}", from=1-3, to=2-3]
	\arrow["{𝐒𝐓\mu^{𝐒𝐓}}"', from=2-1, to=2-2]
	\arrow["{\mu^{𝐒𝐓}}"', from=2-2, to=2-3]
\end{tikzcd}\]
\[\begin{tikzcd}
	{𝐒𝐓𝐒𝐒𝐓X} & {𝐒𝐓𝐒𝐓𝐒𝐓X} & {𝐒𝐓𝐒𝐓X} \\
	{𝐒𝐓𝐒𝐓𝐒𝐓X} & {𝐒𝐓𝐒𝐓X} & {𝐒𝐓X}
	\arrow["{𝐒𝐓𝐒\eta^{𝐓}𝐒𝐓}", from=1-1, to=1-2]
	\arrow["{𝐒𝐓𝐒\eta^{𝐓}𝐒𝐓}"', from=1-1, to=2-1]
	\arrow["{\mu^{𝐒𝐓}𝐒𝐓}", from=1-2, to=1-3]
	\arrow["{\mu^{𝐒𝐓}}", from=1-3, to=2-3]
	\arrow["{𝐒𝐓\mu^{𝐒𝐓}}"', from=2-1, to=2-2]
	\arrow["{\mu^{𝐒𝐓}}"', from=2-2, to=2-3]
\end{tikzcd}\]
\[\begin{tikzcd}
	{𝐒𝐓𝐒𝐓𝐓X} & {𝐒𝐓𝐒𝐓𝐒𝐓X} & {𝐒𝐓𝐒𝐓X} \\
	{𝐒𝐓𝐒𝐓𝐒𝐓X} & {𝐒𝐓𝐒𝐓X} & {𝐒𝐓X}
	\arrow["{𝐒𝐓𝐒𝐓\eta^{𝐒}𝐓}", from=1-1, to=1-2]
	\arrow["{𝐒𝐓𝐒𝐓\eta^{𝐒}𝐓}"', from=1-1, to=2-1]
	\arrow["{\mu^{𝐒𝐓}𝐒𝐓}", from=1-2, to=1-3]
	\arrow["{\mu^{𝐒𝐓}}", from=1-3, to=2-3]
	\arrow["{𝐒𝐓\mu^{𝐒𝐓}}"', from=2-1, to=2-2]
	\arrow["{\mu^{𝐒𝐓}}"', from=2-2, to=2-3]
\end{tikzcd}\]

  Let us consider the first case. 
  We employ ① naturality, ② unitality of the \kl{almost distributive law} $n$, ③ unitality of the \kl{monad} $𝐓$, and ④ associativity of the \kl{monad} $𝐒$.
\[\begin{tikzcd}[ampersand replacement=\&,column sep=scriptsize]
	{𝐒𝐒𝐓𝐒𝐓X} \& {𝐒𝐓𝐒𝐓𝐒𝐓X} \& {𝐒𝐒𝐓𝐓𝐒𝐓X} \& {𝐒𝐓𝐒𝐓X} \&\& \\
	{𝐒𝐒𝐓𝐒𝐓X} \& {𝐒𝐒𝐓𝐓𝐒𝐓X} \& {𝐒𝐒𝐓𝐒𝐓X} \& {𝐒𝐓𝐒𝐓X} \& {𝐒𝐒𝐓𝐓X} \& {𝐒𝐓X} \\
	\& {𝐒𝐒𝐓𝐒𝐓X} \& {𝐒𝐒𝐒𝐓𝐓X} \& {𝐒𝐒𝐓𝐓X} \& {𝐒𝐓𝐓X} \& {𝐒𝐓X} \\
	\& {𝐒𝐒𝐓𝐒𝐓X} \& {𝐒𝐒𝐒𝐓𝐓X} \& {𝐒𝐒𝐓𝐓X} \& {𝐒𝐓𝐓X} \& {𝐒𝐓X} \\
	{𝐒𝐒𝐓𝐒𝐓X} \& {𝐒𝐒𝐒𝐓𝐓X} \& {𝐒𝐒𝐓𝐓X} \& {𝐒𝐒𝐓X} \&\& {𝐒𝐓X} \\
	\& {𝐒𝐒𝐒𝐓𝐓X} \& {𝐒𝐒𝐓X} \& {𝐒𝐒𝐓𝐓X} \& {𝐒𝐒𝐓X} \& {𝐒𝐓X} \\
	{𝐒𝐒𝐓𝐒𝐓X} \& {𝐒𝐒𝐒𝐓𝐓X} \& {𝐒𝐒𝐓X} \& {𝐒𝐓𝐒𝐓X} \& {𝐒𝐒𝐓𝐓X} \& {𝐒𝐓X} \\
	{𝐒𝐒𝐓𝐒𝐓X} \& {𝐒𝐓𝐒𝐓𝐒𝐓X} \& {𝐒𝐓𝐒𝐒𝐓𝐓X} \& {𝐒𝐓𝐒𝐓X}
	\arrow["{𝐒\eta^{𝐓}𝐒𝐓𝐒𝐓}", from=1-1, to=1-2]
	\arrow[""{name=0, anchor=center, inner sep=0}, no head, from=1-1, to=2-1]
	\arrow["{𝐒n𝐓𝐒𝐓}", from=1-2, to=1-3]
	\arrow["{\mu^{𝐒}\mu^{𝐓}𝐒𝐓}", from=1-3, to=1-4]
	\arrow[""{name=1, anchor=center, inner sep=0}, no head, from=1-3, to=2-2]
	\arrow[""{name=2, anchor=center, inner sep=0}, no head, from=1-4, to=2-4]
	\arrow["{𝐒𝐒\eta^{𝐓}𝐓𝐒𝐓}", from=2-1, to=2-2]
	\arrow[""{name=3, anchor=center, inner sep=0}, no head, from=2-1, to=3-2]
	\arrow["{𝐒𝐒\mu^{𝐓}𝐒𝐓}", from=2-2, to=2-3]
	\arrow["{\mu^{𝐒}𝐓𝐒𝐓}", from=2-3, to=2-4]
	\arrow[""{name=4, anchor=center, inner sep=0}, no head, from=2-3, to=3-2]
	\arrow["{𝐒n𝐓}", from=2-4, to=2-5]
	\arrow["{\mu^{𝐒}\mu^{𝐓}}", from=2-5, to=2-6]
	\arrow["{𝐒𝐒n𝐓}"', from=3-2, to=3-3]
	\arrow[no head, from=3-2, to=4-2]
	\arrow["{𝐒\mu^{𝐒}𝐓𝐓}"', from=3-3, to=3-4]
	\arrow[""{name=5, anchor=center, inner sep=0}, no head, from=3-3, to=4-3]
	\arrow[""{name=6, anchor=center, inner sep=0}, no head, from=3-4, to=2-5]
	\arrow["{\mu^{𝐒}𝐓𝐓}"', from=3-4, to=3-5]
	\arrow["{𝐒\mu^{𝐓}}"', from=3-5, to=3-6]
	\arrow[""{name=7, anchor=center, inner sep=0}, no head, from=3-5, to=4-5]
	\arrow[""{name=8, anchor=center, inner sep=0}, no head, from=3-6, to=2-6]
	\arrow[no head, from=3-6, to=4-6]
	\arrow["{𝐒𝐒n𝐓}", from=4-2, to=4-3]
	\arrow[no head, from=4-2, to=5-1]
	\arrow["{\mu^{𝐒}𝐒𝐓𝐓}", from=4-3, to=4-4]
	\arrow[no head, from=4-3, to=5-2]
	\arrow["{\mu^{𝐒}𝐓𝐓}", from=4-4, to=4-5]
	\arrow[""{name=9, anchor=center, inner sep=0}, no head, from=4-4, to=5-3]
	\arrow["{𝐒\mu^{𝐓}}", from=4-5, to=4-6]
	\arrow[""{name=10, anchor=center, inner sep=0}, no head, from=4-6, to=5-6]
	\arrow["{𝐒𝐒n𝐓}", from=5-1, to=5-2]
	\arrow[no head, from=5-1, to=7-1]
	\arrow["{\mu^{𝐒}𝐒𝐓𝐓}", from=5-2, to=5-3]
	\arrow[""{name=11, anchor=center, inner sep=0}, no head, from=5-2, to=6-2]
	\arrow["{𝐒𝐒\mu^{𝐓}}", from=5-3, to=5-4]
	\arrow["{\mu^{𝐒}𝐓}", from=5-4, to=5-6]
	\arrow[""{name=12, anchor=center, inner sep=0}, no head, from=5-4, to=6-3]
	\arrow[""{name=13, anchor=center, inner sep=0}, no head, from=5-4, to=6-5]
	\arrow[no head, from=5-6, to=6-6]
	\arrow["{𝐒\mu^{𝐒}\mu^{𝐓}}"', from=6-2, to=6-3]
	\arrow[no head, from=6-2, to=7-2]
	\arrow["{𝐒𝐒\eta^{𝐓}𝐓}"', from=6-3, to=6-4]
	\arrow[""{name=14, anchor=center, inner sep=0}, no head, from=6-3, to=7-3]
	\arrow["{𝐒𝐒\mu^{𝐓}}"', from=6-4, to=6-5]
	\arrow[""{name=15, anchor=center, inner sep=0}, no head, from=6-4, to=7-5]
	\arrow["{\mu^{𝐒}𝐓}"', from=6-5, to=6-6]
	\arrow[""{name=16, anchor=center, inner sep=0}, no head, from=6-6, to=7-6]
	\arrow["{𝐒𝐒n𝐓}"', from=7-1, to=7-2]
	\arrow[""{name=17, anchor=center, inner sep=0}, no head, from=7-1, to=8-1]
	\arrow["{𝐒\mu^{𝐒}\mu^{𝐓}}"', from=7-2, to=7-3]
	\arrow["{𝐒\eta^{𝐓}𝐒𝐓}"', from=7-3, to=7-4]
	\arrow["{𝐒n𝐓}"', from=7-4, to=7-5]
	\arrow[""{name=18, anchor=center, inner sep=0}, no head, from=7-4, to=8-4]
	\arrow["{\mu^{𝐒}\mu^{𝐓}}"', from=7-5, to=7-6]
	\arrow["{𝐒\eta^{𝐓}𝐒𝐓𝐒𝐓}"', from=8-1, to=8-2]
	\arrow["{𝐒𝐓𝐒n𝐓}"', from=8-2, to=8-3]
	\arrow["{𝐒𝐓\mu^{𝐒}\mu^{𝐓}}"', from=8-3, to=8-4]
	\arrow["{②}"{description}, draw=none, from=0, to=1]
	\arrow["{①}"{description}, draw=none, from=1, to=2]
	\arrow["{③}"{description}, draw=none, from=3, to=4]
	\arrow["{④}"{description}, draw=none, from=5, to=7]
	\arrow["{①}"{description}, draw=none, from=6, to=4]
	\arrow["{①}"{description}, draw=none, from=8, to=6]
	\arrow["{①}"{description}, draw=none, from=9, to=10]
	\arrow["{①}"{description}, draw=none, from=11, to=12]
	\arrow["{③}"{description}, draw=none, from=12, to=13]
	\arrow["{②}"{description}, draw=none, from=14, to=15]
	\arrow["{①}"{description}, draw=none, from=15, to=16]
	\arrow["{①}"{description}, draw=none, from=17, to=18]
\end{tikzcd}\]
  The other two cases follow a similar reasoning.
\end{proof}

\begin{remark}
	The \kl{magmoid} of \kl{normalized stochastic kernels} (\Cref{prop:magmoid-normalized-stochastic-kernels}) can be equivalently defined as the \kl{Kleisli magmoid} of the \kl{normalized distribution magmad}, $𝕂(𝐌𝐃,\smash{μ_{∘}^{𝐌𝐃}},\smash{η_{∘}^{𝐌𝐃}})$.
\end{remark}

The problem with \kl{magmoids} is that their syntax is not particularly pleasant. We would like to describe problems in programming syntax, not in terms of compositions and identities. It is time to slightly adapt the monadic metalanguage \cite{moggi1991notions,wadler92}.

\subsection{The Magmadic metalanguage}
\label{sec:magmoidal-metalanguage}

Moggi's \emph{system PL} \cite{moggi1991notions} is a metalanguage for monadic computation, popularly known as \emph{do-notation}~\cite{hudak2007history}. Do-notation takes semantics over cartesian closed categories endowed with a monad. It consists of two rules, and it leaves variable management to the metatheory.
The monadic metalanguage is implicitly right associative: variables are bound to the right. For \kl{magmads}, this convention becomes a choice: we may want to associate to the left, and indeed we will see this is sometimes desirable.

We start by recalling the definition of \kl{magmad}—and \kl{monad}—in relative form (\Cref{def:magmad-relative}) and monadic do-notation (\Cref{def:donotation}). For \kl{magmads}, we also introduce left-associating do-notation (\Cref{def:left-do-notation}).

\begin{definition}[Magmad in relative form]
  \label[definition]{def:magmad-relative}
  \AP A \intro{magmad}, $(𝐓,η^{𝐓},β^{𝐓})$, consists of an assignment on sets, $𝐓X$ for each set $X$, a natural family of functions $\smash{η^{𝐓}} ፡ X → 𝐓X$ and a natural family of functions $\smash{β^{𝐓}} ፡ 𝐓 X × (X → 𝐓 Y) → 𝐓 Y$ (called \emph{binding}), satisfying the first two axioms of the following three.
  \begin{enumerate}
    \item $β^{𝐓}(η^{𝐓}(x))(λ x. f(x)) = f(x)$, left unitality;
    \item $β^{𝐓}(t)(λ x. η^{𝐓}(x)) = t$, right unitality; and
    \item $β^{𝐓}(t)(λ x. β^{𝐓}(f(x))(g)) = β^{𝐓}(β^{𝐓}(t)(λ x . f(x)))(λ y . g(y))$, associativity.
  \end{enumerate}
  A \kl{monad} is a \kl{magmad} additionally satisfying the third axiom.
\end{definition}

\newsavebox{\codeBlockDoReturn}
\begin{lrbox}{\codeBlockDoReturn}
\begin{minipage}{0.22\textwidth}
\begin{lstlisting}[language=Racket,numbers=none,mathescape]
(rightDo T
  return (x1 ...))
\end{lstlisting}
\end{minipage}
\end{lrbox}
\newsavebox{\codeBlockDoStatement}
\begin{lrbox}{\codeBlockDoStatement}
\begin{minipage}{0.21\textwidth}
\begin{lstlisting}[language=Racket,numbers=none,mathescape]
(rightDo T
  (x1 ...) <- m
  rest ...)
\end{lstlisting}
\end{minipage}
\end{lrbox}
\newsavebox{\codeBlockDoRest}
\begin{lrbox}{\codeBlockDoRest}
\begin{minipage}{0.14\textwidth}
\begin{lstlisting}[language=Racket,numbers=none,mathescape]
(rightDo T
  rest ...)
\end{lstlisting}
\end{minipage}
\end{lrbox}

\begin{definition}[Right do notation]
  \label[definition]{def:donotation}  
  \AP \phantomintro{right do-notation}\intro{Do notation} is defined inductively by the following two rules. The \emph{(i)} first rule interprets an statement as the \kl{magmad} binding. The \emph{(ii)} second rule interprets the return as the \kl{magmad} unit.
  \begin{align*}
    \usebox{\codeBlockDoStatement} \quad &\overset{\emph{(i)}}{=}\ \beta^{T}(m)\left(\lambda (x_1,\dots,x_n). \quad \usebox{\codeBlockDoRest}\quad \right);\\
	\usebox{\codeBlockDoReturn} \quad &\overset{\emph{(i)}}{=}\ \eta^{T}(x_1,\dots,x_n). 
  \end{align*}
\end{definition}

\newsavebox{\codeBlockLmo}
\begin{lrbox}{\codeBlockLmo}
\begin{minipage}{0.25\textwidth}
\begin{lstlisting}[language=Racket,numbers=none,mathescape]
(leftDo T
  (x ...) <- m1
  (y ...) <- m2
  rest ...)
\end{lstlisting}
\end{minipage}	
\end{lrbox}
\newsavebox{\codeBlockLmt}
\begin{lrbox}{\codeBlockLmt}
\begin{minipage}{0.40\textwidth}
\begin{lstlisting}[language=Racket,numbers=none,mathescape]
(leftDo T
  (x ... y ...) <- (rightDo T
          (x ...) <- m1
          (y ...) <- m2
          return (x ... y ...))
  rest ...)
\end{lstlisting}
\end{minipage}
\end{lrbox}
\newsavebox{\codeBlockLlo}
\begin{lrbox}{\codeBlockLlo}
\begin{minipage}{0.22\textwidth}
\begin{lstlisting}[language=Racket,numbers=none,mathescape]
(leftDo T
  (x ...) <- m
  return (y ...))
\end{lstlisting}
\end{minipage}
\end{lrbox}
\newsavebox{\codeBlockLlt}
\begin{lrbox}{\codeBlockLlt}
\begin{minipage}{0.22\textwidth}
\begin{lstlisting}[language=Racket,numbers=none,mathescape]
(rightDo T
  (x ...) <- m
  return (y ...))
\end{lstlisting}
\end{minipage}
\end{lrbox}
\newsavebox{\codeBlockLro}
\begin{lrbox}{\codeBlockLro}%
\begin{minipage}{0.22\textwidth}
\begin{lstlisting}[language=Racket,numbers=none,mathescape]
(leftDo T 
  return (x ...))
\end{lstlisting}
\end{minipage}
\end{lrbox}
\newsavebox{\codeBlockLrt}
\begin{lrbox}{\codeBlockLrt}%
\begin{minipage}{0.22\textwidth}
\begin{lstlisting}[language=Racket,numbers=none,mathescape]
(rightDo T 
  return (x ...))
\end{lstlisting}
\end{minipage}
\end{lrbox}

\begin{definition}[Left do notation]
  \label[definition]{def:left-do-notation}
  \AP \intro{Left do notation} is defined inducitvely by the following rules. The \emph{(i)} first rule groups the two first lines to the left; \emph{(ii)} the second rule deals with the base case of a single line; and \emph{(iii)} the third rule deals with the base case of no remaining lines.
  \begin{align*}
	\usebox{\codeBlockLmo}
	\quad &\overset{(i)}{=} \quad
	\usebox{\codeBlockLmt}\quad;\\
	\usebox{\codeBlockLlo} 
	\quad  &\overset{(ii)}{=} \quad
	\usebox{\codeBlockLlt}\quad;\\
	\usebox{\codeBlockLro} 
	\quad &\overset{(iii)}{=} \quad 
	\usebox{\codeBlockLrt}\quad.
  \end{align*}
\end{definition}

\newsavebox{\codeBlockAssocDoL}
\begin{lrbox}{\codeBlockAssocDoL}
\begin{minipage}{0.35\textwidth}
\begin{lstlisting}[language=Racket,numbers=none,mathescape]
(rightDo T
  ...
  (y ...) <- (rightDo T
	(x1 ...) <- E1
	...
	(xn ...) <- En
	return (y ...))
  ...)
\end{lstlisting}
\end{minipage}	
\end{lrbox}
\newsavebox{\codeBlockAssocDoR}
\begin{lrbox}{\codeBlockAssocDoR}
\begin{minipage}{0.35\textwidth}
\begin{lstlisting}[language=Racket,numbers=none,mathescape]
(rightDo T
  ...
  (x1 ...) <- E1
  ...
  (xn ...) <- En
  ...)
\end{lstlisting}
\end{minipage}	
\end{lrbox}

\begin{remark}[Associativity]
  Whenever the \kl{magmad} is a \kl{monad}, additionally, the following equation holds. As a consequence, left do notation and right do notation coincide for any \kl{monad}, even if they do not coincide for a \kl{magmad}.
  \begin{align*}
	\usebox{\codeBlockAssocDoL}\quad = \quad\usebox{\codeBlockAssocDoR}
  \end{align*}
\end{remark}

\newsavebox{\codeDoDefL}
\begin{lrbox}{\codeDoDefL}
\begin{minipage}{0.2\textwidth}
\begin{lstlisting}[language=Racket,numbers=none,mathescape]
(do ...)
\end{lstlisting}
\end{minipage}	
\end{lrbox}

\newsavebox{\codeDoDefR}
\begin{lrbox}{\codeDoDefR}
\begin{minipage}{0.2\textwidth}
\begin{lstlisting}[language=Racket,numbers=none,mathescape]
(leftDo MD ...)
\end{lstlisting}
\end{minipage}	
\end{lrbox}

\begin{remark}[Operational left-associativity]
  When we compute, we start from the first line and apply the next line at each step—implicitly, we associate to the left.	In functional programming, the bind operator of a monad takes the rest of the program as an argument—implicitly, it associates to the right. From now on, we use left-do notation for normalized distributions, instead of right-do notation.
  \[\usebox{\codeDoDefL}\quad = \quad\usebox{\codeDoDefR}\]
\end{remark}

\subsection{Application—Modes of update}
\label{sec:update}

Updating a prior distribution with an observation can be done in different ways, depending on the structure of our model and the observation. Following Jacobs, we call them \emph{modes of update} \cite{jacobs2019mathematics}.

This section formalizes modes of update as non-associative phenomena. The magmadic metalanguage allows us to specify different modes of update: we compare Bayes' update to the filtering of impossible cases (\Cref{prop:bayes-removal}); we compare Pearl's update to Jeffrey's update.

\paragraph{Three-statement case} Let us first clarify the anatomy of a basic Bayesian update. It contains three steps (\emph{i}) we have a prior, (\emph{ii}) we apply a stochastic channel, (\emph{iii}) and we observe some property through that channel: the association of these three steps yields two modes of update.

\begin{proposition}[Bayes update]
  \label{prop:bayes-removal}
  Let $\mathsf{prior} ፡ 1 → 𝐃X$ and $\mathsf{channel} ፡ X → 𝐃Y$ be stochastic kernels; let $\mathsf{property} ፡ Y → 𝐌1$ be a partial function. There are two modes of update: (i) Bayes' update and (ii) the removal of impossible cases.
  \begin{equation*}
  (i)\quad
  \begin{minipage}{0.35\textwidth}
\begin{lstlisting}[language=Racket,numbers=none,mathescape]
(do (x) <- prior
    (y) <- (channel x)
    () <- (property y)
    return (x))
\end{lstlisting}
  \end{minipage} 
  \qquad (ii) \quad
  \begin{minipage}{0.35\textwidth}
\begin{lstlisting}[language=Racket,numbers=none,mathescape]
(do (x) <- prior
    () <- (do
        (y) <- (channel y)
        () <- (property y)
        return ())
    return (x))    
\end{lstlisting}
\end{minipage}
\end{equation*}
\end{proposition}
\begin{proof}
  We have two ways of associating three statements—i.e., $2$ is the third Catalan number. Let us write $\mathsf{prior} = \sum\nolimits_{i=0} λ_i\ketb{x_i}$; let us write $\mathsf{channel}(\cBlu{x_i}) = \sum\nolimits_{j=0} ρ_{ij}\ketb{y_{ij}}$ and let $\mathsf{property}(\cBlu{y_{ij}}) = u_{ij}$, where $u_{ij} = 0$ or $u_{ij} = 1$. We compute explicitly the two possible updates.
  \begin{align*}
    \emph{(i)}\quad &
      {\textstyle\sum\nolimits_{i}}\ λ_{i} \ketb{x_i} &
    \emph{(i)}\quad &
      {\textstyle\sum\nolimits_{i}}\ λ_{i} \ketb{x_i}
      \\[-0.4em]
    \emph{(ii)}\quad &
      {\textstyle\sum\nolimits_{ij}}\ λ_{i}ρ_{j} \ketb{x_i} \ketb{y_{ij}} &
    \emph{(ii)}\quad &
      {\textstyle\sum\nolimits_{i}}\ λ_{i} \ketb{x_i} [{\textstyle\sum\nolimits_{j}}\ ρ_{j} \ketb{y_{ij}}]
      \\[-0.4em]
    \emph{(iii)}\quad &
      {\textstyle\sum\nolimits_{ij}}\ λ_{i}ρ_{j}u_{ij} \ketb{x_i} &
    \emph{(iii)}\quad &
      {\textstyle\sum\nolimits_{i}}\ λ_{i} \ketb{x_i} [{\textstyle\sum\nolimits_{j}}\ ρ_{j}u_{ij} \ketb{}]
      \\[-0.4em]
    \emph{(n)}\quad &
      {\textstyle\sum\nolimits_{i, u_{ij}=1}}\ λ_{i}ρ_{j}\ketb{x_i} / {{\textstyle\sum\nolimits_{i, u_{ij}=1}}\ λ_{i}ρ_{j}} &
    \emph{(n)}\quad &
      {\textstyle\sum\nolimits_{i, {\textstyle\sum\nolimits_{j}}\ ρ_{j}u_{ij} ≠ 0}}\ λ_{i} \ketb{x_i}
  \end{align*}
  Note how the first corresponds to Bayes' update, while the second filters out impossible cases. 
\end{proof}

\begin{example}[Monty Hall problem]
  The Monty Hall problem has two solutions, depending on associativity: they correspond to the usage of \kl{left do notation} ($\mathsf{leftDo}/\mathsf{do}$) or \kl{right do notation} ($\mathsf{rightDo}$). The left-associative solution applies Bayes udpate and prefers to change doors ($\tfrac{1}{3} \ketb{M} + \tfrac{2}{3} \ketb{R}$); the right-associative solution filters the impossible case of having the prize on the left door without changing the proportion of probability for the rest of the possibilities ($\tfrac{1}{2} \ketb{M} + \tfrac{1}{2} \ketb{R}$).
\begin{equation*}
\begin{minipage}{0.40\textwidth}
\begin{lstlisting}[language=Racket,numbers=none]
(leftDo MD
  (prize) <- uniformDoor
  (door) <- (host prize)
  () <- (observeLeft door)
  return (prize))
\end{lstlisting}
\end{minipage}
\quad ≠ \quad
\begin{minipage}{0.40\textwidth}
\begin{lstlisting}[language=Racket,numbers=none]
(rightDo MD
  (prize) <- uniformDoor
  (door) <- (host prize)
  () <- (observeLeft door)
  return (prize))
\end{lstlisting}
\end{minipage}
\end{equation*}
\end{example}

\paragraph{Four-statement case} When the property depends on an extra stochastic parameter, there are instead four steps to consider: (\emph{i}) we sample the parameter, (\emph{ii}) we have a prior, (\emph{iii}) we apply a stochastic channel, (\emph{iv}) and we observe some property.
The association of these four steps matters, and there are 5 ways of associating four statements. However, stochastic statements do associate, reducing the number of possibilities to four: of these, only two perform a Bayesian update with the three last statements.

\begin{proposition}[Pearl's and Jeffrey's updates]
  \label{prop:pearl-jeffrey}
  \AP Let $\mathsf{parameter} ፡ 1 → 𝐃Z$, and $\mathsf{prior} ፡ 1 → 𝐃X$, and $\mathsf{channel} ፡ X → 𝐃Y$ be stochastic kernels; let $\mathsf{property} ፡ Y × Z → 𝐌1$ be a partial function. There are two modes of update extending Bayes' update: (i) Pearl's update and (ii) Jeffrey's update.
  Pearl's and Jeffrey's updates do not coincide in general.
\begin{equation*}
\begin{minipage}{0.35\textwidth}
\begin{lstlisting}[language=Racket,numbers=none,mathescape]
(do (z) <- parameter
    (x) <- prior
    (y) <- (channel x)
    () <- (property y z)
   return (x))
\end{lstlisting}
\end{minipage}
\ \ \ ≠ \ \ 
\begin{minipage}{0.45\textwidth}
\begin{lstlisting}[language=Racket,numbers=none,mathescape]
(do (z) <- parameter
    (x y) <- (do
       (x) <- prior
       (y) <- (channel x)
       () <- (property y z)
       return (x y))
    return (x))
\end{lstlisting}
\end{minipage}
\end{equation*}
\end{proposition}
\begin{proof}
  One of our examples (the \emph{Unclear Test problem} in \Cref{sec:unclear-test}) follows these two modes of update and obtains different results. We do not characterize these further.
\end{proof}

\begin{remark}[Interpreting Pearl and Jeffrey]
  In the \emph{Unclear Test problem} \Cref{sec:unclear-test}, Pearl's update uses a single normalization constant for both the positive-result and the negative-result branches; Jeffrey's update normalizes both separately. The greater relative likelihood of the positive branch makes Pearl's update yield a higher probability of illness.
  In the work of Jacobs \cite{jacobs2019mathematics}, Pearl's update has been characterized as \emph{increasing validity}, while Jeffrey's update as \emph{decreasing divergence}.
\end{remark}

\section{Reasoning with the Metalanguage\label{sec:reasoning-metalanguage}}

This section introduces some principles to reason syntactically about programs involving observations in the magmadic metalanguage. Each axiom is based on one of the semantic principles we have established until now, except for the \kl{disintegration axiom}, which we now introduce.

\Cref{fig:axioms-metalanguage} lists the six axioms we will use: ① commutativity of $\mathsf{observe}$ and the \kl{Frobenius} axiom, with an extra equation stating that $\mathsf{F}$ has \kl{full support} (justified in \Cref{sec:partialfrobenius}); ② the \kl{disintegration} axiom and its algebraic version (justified in \Cref{sec:disintegration}); ③ the exchange axiom (justified in \Cref{sec:commutative-and-left-ex}); ④ discarding whenever $\mathsf{F}$ is stochastic, copying whenever $\mathsf{F}$ is partial, and the \kl{normal magmad} equation for arbitrary morphisms (justified in \Cref{sec:affine-relevant}); ⑤ associativity, whenever $\mathsf{F}$ is partial, $\mathsf{G}$ is partial, or $\mathsf{H}$ is stochastic (justified in \Cref{lemma:relevance}).

From these, we highlight how \emph{commutativity} and \emph{monoidality}—unlike in the monadic case—do not necessarily coincide for \kl{magmads}: we introduce a notion of \kl{left-exchanging magmad} (\Cref{def:left-exchange}).

\newsavebox{\codeAxTotalL}%
\begin{lrbox}{\codeAxTotalL}%
\begin{minipage}{0.20\textwidth}
\begin{lstlisting}[language=Racket,numbers=none,mathescape]
(do (x ...) <- F
    return ())
\end{lstlisting}
\end{minipage}	
\end{lrbox}

\newsavebox{\codeAxDeterministicL}
\begin{lrbox}{\codeAxDeterministicL}
\begin{minipage}{0.3\textwidth}
\begin{lstlisting}[language=Racket,numbers=none,mathescape]
(do (x ...) <- F
    return (x ... x ...))
\end{lstlisting}
\end{minipage}	
\end{lrbox}

\newsavebox{\codeAxTotalR}
\begin{lrbox}{\codeAxTotalR}
\begin{minipage}{0.20\textwidth}
\begin{lstlisting}[language=Racket,numbers=none,mathescape]
(do return ())
\end{lstlisting}
\end{minipage}	
\end{lrbox}

\newsavebox{\codeAxDeterministicR}
\begin{lrbox}{\codeAxDeterministicR}
\begin{minipage}{0.3\textwidth}
\begin{lstlisting}[language=Racket,numbers=none,mathescape]
(do (x ...) <- F
    (y ...) <- F
    return (x ... y ...))
\end{lstlisting}
\end{minipage}	
\end{lrbox}

\newsavebox{\codeAxInterchangeL}
\begin{lrbox}{\codeAxInterchangeL}
\begin{minipage}{0.3\textwidth}
\begin{lstlisting}[language=Racket,numbers=none,mathescape]
(x ...) <- (F a ...)
(y ...) <- (G b ...)
\end{lstlisting}
\end{minipage}	
\end{lrbox}

\newsavebox{\codeAxInterchangeR}
\begin{lrbox}{\codeAxInterchangeR}
\begin{minipage}{0.3\textwidth}
\begin{lstlisting}[language=Racket,numbers=none,mathescape]
(y ...) <- (G b ...)
(x ...) <- (F a ...)
\end{lstlisting}
\end{minipage}	
\end{lrbox}

\newsavebox{\codeAxFrobeniusL}
\begin{lrbox}{\codeAxFrobeniusL}
\begin{minipage}{0.3\textwidth}
\begin{lstlisting}[language=Racket,numbers=none,mathescape]
(x ...) <- (F a)
() <- (observe a b)
\end{lstlisting}
\end{minipage}	
\end{lrbox}

\newsavebox{\codeAxFrobeniusR}
\begin{lrbox}{\codeAxFrobeniusR}
\begin{minipage}{0.3\textwidth}
\begin{lstlisting}[language=Racket,numbers=none,mathescape]
() <- (observe a b)
(x ...) <- (F b)
\end{lstlisting}
\end{minipage}	
\end{lrbox}

\newsavebox{\codeAxFcommL}
\begin{lrbox}{\codeAxFcommL}
\begin{minipage}{0.3\textwidth}
\begin{lstlisting}[language=Racket,numbers=none,mathescape]
() <- (observe a b)
\end{lstlisting}
\end{minipage}	
\end{lrbox}

\newsavebox{\codeAxFcommR}
\begin{lrbox}{\codeAxFcommR}
\begin{minipage}{0.3\textwidth}
\begin{lstlisting}[language=Racket,numbers=none,mathescape]
() <- (observe b a)
\end{lstlisting}
\end{minipage}	
\end{lrbox}

\newsavebox{\codeAxDisintegrationL}
\begin{lrbox}{\codeAxDisintegrationL}
\begin{minipage}{0.35\textwidth}
\begin{lstlisting}[language=Racket,numbers=none,mathescape]
(do (x ... y ...) <- F
    return (x ... y ...))
\end{lstlisting}
\end{minipage}	
\end{lrbox}

\newsavebox{\codeAxDisintegrationR}
\begin{lrbox}{\codeAxDisintegrationR}
\begin{minipage}{0.35\textwidth}
\begin{lstlisting}[language=Racket,numbers=none,mathescape]
(do (x ... y ...) <- F
    (y ...) <- ((d F) x ...)
    return (x ... y ...))
\end{lstlisting}
\end{minipage}	
\end{lrbox}

\newsavebox{\codeAxDisintegrateA}
\begin{lrbox}{\codeAxDisintegrateA}
\begin{minipage}{0.2\textwidth}
\begin{lstlisting}[language=Racket,numbers=none,mathescape]
((d F) x ...)
\end{lstlisting}
\end{minipage}	
\end{lrbox}

\newsavebox{\codeAxDisintegrateB}
\begin{lrbox}{\codeAxDisintegrateB}
\begin{minipage}{0.42\textwidth}
\begin{lstlisting}[language=Racket,numbers=none,mathescape]
(do (a ... y ...) <- F
    () <- (observe (x ...) (a ...))
    return (y ...))
\end{lstlisting}
\end{minipage}	
\end{lrbox}

\newsavebox{\codeAxAssocL}
\begin{lrbox}{\codeAxAssocL}
\begin{minipage}{0.38\textwidth}
\begin{lstlisting}[language=Racket,numbers=none,mathescape]
(do (x ...) <- F
    (y ...) <- G
    (z ...) <- H
    return (x ... y ... z ...))
\end{lstlisting}
\end{minipage}	
\end{lrbox}

\newsavebox{\codeAxAssocR}
\begin{lrbox}{\codeAxAssocR}
\begin{minipage}{0.38\textwidth}
\begin{lstlisting}[language=Racket,numbers=none,mathescape]
(do (x ...) <- F
    (y ... z ...) <- (do
       (y ...) <- G
       (z ...) <- H
       return (y ... z ...))
    return (x ... y ... z ...))
\end{lstlisting}
\end{minipage}	
\end{lrbox}

\newsavebox{\codeAxQtR}
\begin{lrbox}{\codeAxQtR}
\begin{minipage}{0.35\textwidth}
\begin{lstlisting}[language=Racket, numbers=none]
(do (x1 ...) <- F
    (y1 ...) <- F
    return (x1 ...))
\end{lstlisting}
\end{minipage}
\end{lrbox}

\newsavebox{\codeAxQtL}
\begin{lrbox}{\codeAxQtL}
\begin{minipage}{0.35\textwidth}
\begin{lstlisting}[language=Racket, numbers=none]
(do (x1 ...) <- F
    return (x1 ...))
\end{lstlisting}
\end{minipage}
\end{lrbox}

\newsavebox{\codeAxSupportL}
\begin{lrbox}{\codeAxSupportL}
\begin{minipage}{0.4\textwidth}
\begin{lstlisting}[language=Racket,numbers=none,mathescape]
(do (y ...) <- F
    () <- (observe (x ...) (y ...))
    return (x ...))
\end{lstlisting}
\end{minipage}	
\end{lrbox}

\newsavebox{\codeAxSupportR}
\begin{lrbox}{\codeAxSupportR}
\begin{minipage}{0.4\textwidth}
\begin{lstlisting}[language=Racket,numbers=none,mathescape]
(do return (x ...))
\end{lstlisting}
\end{minipage}	
\end{lrbox}

\begin{figure}[ht!]
\begin{framed}
\centering
\begin{align*}
  \usebox{\codeAxFcommL} \quad 
  & \overset{①}{=} \quad \usebox{\codeAxFcommR} \\
  \usebox{\codeAxFrobeniusL} \quad 
  & \overset{①}{=} \quad \usebox{\codeAxFrobeniusR} \\
  \usebox{\codeAxSupportR} \quad 
  & \overset{①}{=} \quad \usebox{\codeAxSupportL} \\
  \usebox{\codeAxDisintegrationL} \quad 
  & \overset{②}{=} \quad \usebox{\codeAxDisintegrationR} \\
  \usebox{\codeAxDisintegrateA} \quad 
  & \overset{②}{=} \quad \usebox{\codeAxDisintegrateB} \\
  \usebox{\codeAxInterchangeL} \quad 
  & \overset{③}{=} \quad \usebox{\codeAxInterchangeR} \\
  \usebox{\codeAxTotalL} \quad 
  & \overset{④}{=} \quad \usebox{\codeAxTotalR} \\
  \usebox{\codeAxDeterministicL} \quad 
  & \overset{④}{=} \quad \usebox{\codeAxDeterministicR} \\
  \usebox{\codeAxQtL} \quad 
  & \overset{④}{=} \quad \usebox{\codeAxQtR} \\
  \usebox{\codeAxAssocL} \quad 
  & \overset{⑤}{=} \quad \usebox{\codeAxAssocR}
\end{align*}
\caption{Principles of the normalized stochastic magmadic metalanguage. Justified in ① \Cref{sec:partialfrobenius}, ② \Cref{sec:disintegration}, ③ \Cref{sec:commutative-and-left-ex}, ④ \Cref{sec:affine-relevant}, and ⑤ \Cref{lemma:relevance}.\label{fig:axioms-metalanguage}}
\end{framed}
\end{figure}

\subsection{Axiom—Partial Frobenius\label{sec:partialfrobenius}}

We ask for the existence of a partial equality check, $\mathsf{observe} ፡ X × X → 𝗡1$, which forces its two arguments to coincide: for the \kl{normalized distribution magmad}, this is defined by $\mathsf{observe}(x,y) = 0$ when $x ≠ y$ and by $\mathsf{observe}(x,x) = 1\ket{}$ otherwise. Apart from commutativity, the only axiom that we ask of this morphism is that it locally forces the equality of its two arguments.

\begin{definition}
  \label{def:partialFrobenius}
  \AP A \kl{magmad}, $(𝗡,μ^{𝗡},η^{𝗡})$, is \intro{partial Frobenius} whenever there exists, for each object $X$, a morphism, $\mathsf{observe} ፡ X × X → 𝗡1$, satisfying \emph{(i)} commutativity and \emph{(ii)} the Frobenius axiom.
  \begin{align*}
  \usebox{\codeAxFcommL}
  \quad & \overset{\emph{(ii)}}{=} \quad 
  \usebox{\codeAxFcommR}
  \\  
  \usebox{\codeAxFrobeniusL}
  \quad & \overset{\emph{(ii)}}{=} \quad 
  \usebox{\codeAxFrobeniusR}
  \end{align*}
\end{definition}

\begin{remark}
  \kl{Observe} operators correspond to \emph{partial Frobenius structures} in the sense of restriction categories \cite{diliberti21} and \emph{partial Markov categories} \cite{dilavore23:evidential,partialMarkov25}. The \emph{comparator} is obtained from the $\mathsf{observe}$ morphisms and variable-copying. We use it to derive a notion of \kl{full support}.
\end{remark}

\begin{definition}[Full support]
  \AP A Kleisli morphism, $f ፡ X → 𝗡 Y$, has \intro{full
  support} whenever the following equation holds.
  \[\usebox{\codeAxSupportL} \quad = \quad \usebox{\codeAxSupportR}\]
\end{definition}

\subsection{Axiom—Disintegration}
\label{sec:disintegration}

Disintegration was introduced by Cho and Jacobs \cite{cho2019disintegration} as a synthetic axiom for extracting a conditional channel, $d(f) ፡ X → 𝗡 Y$, from a joint distribution, $f ፡ 1 → 𝗡(X × Y)$ (\Cref{def:disintegration-existential}). For the \kl{normalized distribution magmad}, it is witnessed by the existence of conditional distributions.
\[
\usebox{\codeAxDisintegrationL}\quad \overset{\emph{(i)}}{=} \quad
\usebox{\codeAxDisintegrationR}
\]

A limitation of its formulation is that it is a non-equational axiom: it postulates that there exists such disintegration, but not how to construct it; in fact, it is unique if and only if the two projections coincide \cite[Proposition 11.15]{fritz2020synthetic}.
The situation is more fortunate in the \kl{normalized} case. A disintegration can always be constructed, and we may substitute the existential axiom by an algebraic one (\Cref{def:disintegration}).\footnote{Both the existential and the algebraic formulations of disintegration can be generalized to arbitrary monoidal magmoids. For the purposes of this text, we work exclusively with Kleisli magmoids in the category of sets.}

\begin{definition}[Disintegration, {{(c.f.,~\cite[Definition 3.5]{cho2019disintegration})}}]
  \label[definition]{def:disintegration-existential}
  \AP A \kl{magmad}, $(𝗡,μ^{𝗡},η^{𝗡})$ satisfies the \intro{disintegration axiom} whenever, for every function $f ፡ A → 𝗡(X × Y)$, there exists some function $d(f) ፡ X → 𝗡(Y)$ such that the following diagram commutes.
\[\begin{tikzcd}[ampersand replacement=\&]
	A \& {𝗡(X × Y)} \& {𝗡(X × X)} \\
	{𝗡(X × Y)} \& {𝗡𝗡(X× Y)} \& {𝗡(X × 𝗡(Y))}
	\arrow["f", from=1-1, to=1-2]
	\arrow["f"', from=1-1, to=2-1]
	\arrow["{𝗡(\delta × \varepsilon)}", from=1-2, to=1-3]
	\arrow["{𝗡(X × d(f))}", from=1-3, to=2-3]
	\arrow["{\mu^{𝗡}}"', from=2-2, to=2-1]
	\arrow["{\psi_r}"', from=2-3, to=2-2]
\end{tikzcd}\]
In that particular case, we say that $d(f)$ is the disintegration of $f$.
\end{definition}

\begin{definition}[Algebraic disintegration]
  \label[definition]{def:disintegration}
  \AP A \kl{magmad}, $(𝗡,μ^{𝗡},η^{𝗡})$ satisfies the \intro{algebraic disintegration} axiom whenever, for each $f ፡ A → 𝗡(X × Y)$, the following composite, $d(f) ፡ X → 𝗡Y$, makes the previous diagram commute.
  \[\begin{tikzcd}[ampersand replacement=\&,column sep=large]
    X \& {𝗡(X × X × Y)} \& {𝗡(𝗡(1) × Y)} \& {𝗡(Y)}
    \arrow["{(X \times f) 𑊩 \psi_{\ell} }", from=1-1, to=1-2]
    \arrow["{𝗡(\mathsf{obs} × Y)}", from=1-2, to=1-3]
    \arrow["{\psi_{\er} 𑊩 \mu^{𝗡}}", from=1-3, to=1-4]
  \end{tikzcd}\]
  Equivalently, the following definition makes the previous disintegration equation hold.
  \[
  \usebox{\codeAxDisintegrateA}\quad \overset{\emph{(ii)}}{=} \quad
  \usebox{\codeAxDisintegrateB}
  \]
\end{definition}

\subsection{Commutative and Left-Exchanging magmads}
\label{sec:commutative-and-left-ex}

Programming with \kl{commutative magmads} raises subtleties. For instance, we would expect that, for probabilistic programs, exchanging the lines of a program must not change its meaning, as long as the lines do not depend on each other \cite{staton2017commutative,culpepper18,zamdzhiev21}.\footnote{This property is usually called \emph{commutativity} itself.}  That is, whenever $x$ is free in $F_2$ and $y$ is free in $F_1$, the following \emph{exchange rule} holds.
\[
\begin{minipage}{0.35\textwidth}
\begin{lstlisting}[language=Racket, numbers=none]
(do ...
    (x1 ...) <- F1
    (y1 ...) <- F2
    ...)
\end{lstlisting}
\end{minipage}
\quad =\quad
\begin{minipage}{0.35\textwidth}
\begin{lstlisting}[language=Racket, numbers=none]
(do ...
    (y1 ...) <- F2
    (x1 ...) <- F1
    ...)
\end{lstlisting}
\end{minipage}
\]

As it happens for monads, every $\Set$-\kl{magmad} is \emph{strong}, meaning that there exist canonical maps $ψ_{ℓ} ፡ X ⊗ 𝐓(Y) → 𝐓(X ⊗ Y)$ and $ψ_{\er} ፡ 𝐓(X) ⊗ Y → 𝐓(X ⊗ Y)$. In some cases, the two ways in which they induce a map $ψ ፡ 𝐓(X) ⊗ 𝐓(Y) → 𝐓(X ⊗ Y)$ coincides: there exists a canonical way of interpreting two effects in parallel; we say that the \kl{magmad} is \kl[commutative magmad]{commutative} (\Cref{def:commutative-magmad}).

Perhaps surprisingly, the \emph{exchange rule} in \kl{left-do notation} does not follow from \kl[commutative monad]{commutativity} but is a stronger property of its own: \kl[left-exchanging magmad]{left-exchange} (\Cref{def:left-exchange}). \kl{Commutativity} only implies \kl{left-exchange} in the presence of the associativity axiom (\Cref{prop:monad-left-exchange,rem:counterexample-left-exchange}); but for \kl{magmads}, \kl{left-exchange} becomes an independent axiom. We prove that \kl{sesquilaws} induce \kl{left-exchanging magmads} (\Cref{prop:left-exchanging}).

\begin{definition}[Commutative magmad]
  \label[definition]{def:commutative-magmad}
  \AP A \intro{commutative magmad} is a \kl{magmad}, $(𝐓,μ^{𝐓},η^{𝐓})$, such that the following diagram commutes.
\[\begin{tikzcd}[cramped]
	{𝐓 X ⊗ 𝐓Y} & {𝐓(𝐓X ⊗ Y)} & {𝐓𝐓(X ⊗ Y)} \\
	{𝐓(X ⊗ 𝐓Y)} & {𝐓𝐓(X ⊗ Y)} & {𝐓(X ⊗ Y)}
	\arrow["{𝐓\psi_{\ell}}", from=1-1, to=1-2]
	\arrow["{𝐓\psi_{\er}}"', from=1-1, to=2-1]
	\arrow["{𝐓\psi_{\er}}", from=1-2, to=1-3]
	\arrow["{\mu^{𝐓}}", from=1-3, to=2-3]
	\arrow["{𝐓\psi_{\ell}}", from=2-1, to=2-2]
	\arrow["{\mu^{𝐓}}", from=2-2, to=2-3]
\end{tikzcd}\]
\end{definition}

\begin{remark}[Monoidal magmad]
  The general definition of \kl{monoidal magmad}—the generalization outside the category of sets—is not required for our development but is provided on the Appendix. All of the \kl{magmads} and \kl{monads} of this text are monoidal in this general sense.
\end{remark}

\begin{definition}[Left-exchanging magmad]
	\label[definition]{def:left-exchange}
  \AP A \intro{left-exchanging magmad} is a \kl{commutative magmad}, $(𝐓, μ^{𝐓}, η^{𝐓})$, such that, moreover, the following two diagrams commute.
\begin{equation*}
\begin{tikzcd}[column sep=small]
	{𝐓(𝐓 X ⊗ 𝐓Y)} & {𝐓𝐓(X ⊗ Y)} & {𝐓(X ⊗ Y)} \\
	{𝐓𝐓(X ⊗ 𝐓Y)} & {𝐓(X ⊗ 𝐓Y)} & {𝐓𝐓(X ⊗ Y)}
	\arrow["{𝐓\psi}", from=1-1, to=1-2]
	\arrow[""{name=0, anchor=center, inner sep=0}, "{𝐓\psi_{\ell}}"', from=1-1, to=2-1]
	\arrow["{\mu^{𝐓}}", from=1-2, to=1-3]
	\arrow["{\mu^{𝐓}}"', from=2-1, to=2-2]
	\arrow["{𝐓\sigma_{\er}}"', from=2-2, to=2-3]
	\arrow[""{name=1, anchor=center, inner sep=0}, "{\mu^{𝐓}}"', from=2-3, to=1-3]
	\arrow["{①}"{description}, draw=none, from=0, to=1]
\end{tikzcd} 
\begin{tikzcd}[column sep=small]
	{𝐓(𝐓 X ⊗ 𝐓Y)} & {𝐓𝐓(X ⊗ Y)} & {𝐓(X ⊗ Y)} \\
	{𝐓𝐓(𝐓X ⊗ Y)} & {𝐓(𝐓X ⊗ Y)} & {𝐓𝐓(X ⊗ Y)}
	\arrow["{𝐓\psi}", from=1-1, to=1-2]
	\arrow[""{name=0, anchor=center, inner sep=0}, "{𝐓\psi_{\er}}"', from=1-1, to=2-1]
	\arrow["{\mu^{𝐓}}", from=1-2, to=1-3]
	\arrow["{\mu^{𝐓}}"', from=2-1, to=2-2]
	\arrow["{𝐓\sigma_{\ell}}"', from=2-2, to=2-3]
	\arrow[""{name=1, anchor=center, inner sep=0}, "{\mu^{𝐓}}"', from=2-3, to=1-3]
	\arrow["{②}"{description}, draw=none, from=0, to=1]
\end{tikzcd}
\end{equation*}
\end{definition}  

\begin{propositionrep}
  \label[proposition]{prop:monad-left-exchange}
  A \kl{commutative monad} is a \kl{left-exchanging monad}. \label[proposition]{rem:counterexample-left-exchange}
  Not every \kl{commutative magmad} is a \kl{left-exchanging magmad}. 
\end{propositionrep}
\begin{proofsketch}
  For the proof, we employ naturality, the associativity of the \kl{monad}, and the definition of the joint strength. See the Appendix for details.

To prove that associativity is required, consider the \emph{average-writer magmad}, $(𝗤,μ^{𝗤},η^{𝗤})$, given by the functor $𝗤X = X × 𝐌 ℚ^{+}$, with the multiplication defined by $μ^{𝗤}(x,q_1,q_2) = (x, (q_1 + q_2)/ 2)$—with $μ^{𝗤}(x,⊥,q_2) = (x, q_2)$ and  $μ^{𝗤}(x,q_1,⊥) = (x, q_1)$—and the unit defined by $η^{𝗤}(x) = (x, ⊥)$. This is a \kl{commutative magmad} but not a \kl{left-exchanging magmad}.
\end{proofsketch}
\begin{proof}
  We employ ⓪ naturality, ① the associativity of the \kl{monad}, and ② the definition of the joint strength.
\[%
\begin{tikzcd}[ampersand replacement=\&,column sep=scriptsize]
	{𝐓(𝐓 X ⊗ 𝐓Y)} \& {𝐓𝐓( X ⊗ 𝐓Y)} \& {𝐓(X ⊗ 𝐓Y)} \& {𝐓𝐓(X ⊗ Y)} \& \\
	\& {𝐓𝐓( X ⊗ 𝐓Y)} \& {𝐓𝐓𝐓(X \otimes Y)} \& {𝐓𝐓(X ⊗ Y)} \& {𝐓(X \otimes Y)} \\
	{𝐓(𝐓 X ⊗ 𝐓Y)} \& {𝐓𝐓( X ⊗ 𝐓Y)} \& {𝐓𝐓𝐓(X \otimes Y)} \& {𝐓𝐓(X ⊗ Y)} \& {𝐓(X \otimes Y)} \\
	\& {𝐓(𝐓 X ⊗ 𝐓Y)} \& {𝐓𝐓(X ⊗ Y)} \& {𝐓(X \otimes Y)}
	\arrow["{𝐓\psi_{\ell}}", from=1-1, to=1-2]
	\arrow[no head, from=1-1, to=3-1]
	\arrow["{\mu^{𝐓}}", from=1-2, to=1-3]
	\arrow[""{name=0, anchor=center, inner sep=0}, no head, from=1-2, to=2-2]
	\arrow["{𝐓\psi_{\er}}", from=1-3, to=1-4]
	\arrow[""{name=1, anchor=center, inner sep=0}, no head, from=1-4, to=2-4]
	\arrow["{𝐓\psi_{\er}}", from=2-2, to=2-3]
	\arrow[no head, from=2-2, to=3-2]
	\arrow["{\mu^{𝐓}𝐓}", from=2-3, to=2-4]
	\arrow[""{name=2, anchor=center, inner sep=0}, no head, from=2-3, to=3-3]
	\arrow["{\mu^{𝐓}}", from=2-4, to=2-5]
	\arrow[""{name=3, anchor=center, inner sep=0}, no head, from=2-5, to=3-5]
	\arrow["{𝐓\psi_{\ell}}", from=3-1, to=3-2]
	\arrow[""{name=4, anchor=center, inner sep=0}, no head, from=3-1, to=4-2]
	\arrow["{𝐓\psi_{\er}}", from=3-2, to=3-3]
	\arrow["{𝐓\mu^{𝐓}}", from=3-3, to=3-4]
	\arrow["{\mu^{𝐓}}", from=3-4, to=3-5]
	\arrow[""{name=5, anchor=center, inner sep=0}, no head, from=3-4, to=4-3]
	\arrow[no head, from=3-5, to=4-4]
	\arrow["\psi", from=4-2, to=4-3]
	\arrow["{\mu^{𝐓}}", from=4-3, to=4-4]
	\arrow["{⓪}"{description}, draw=none, from=0, to=1]
	\arrow["{①}"{description}, draw=none, from=2, to=3]
	\arrow["{②}"{description}, draw=none, from=4, to=5]
\end{tikzcd}\]
Note how the proof employs the associativity of the monad.
To prove that associativity is required, consider the \emph{average-writer magmad}, $(𝗤,μ^{𝗤},η^{𝗤})$, given by the functor $𝗤X = X × 𝐌 ℚ^{+}$, with the multiplication defined by $μ^{𝗤}(x,q_1,q_2) = (x, (q_1 + q_2)/ 2)$—with $μ^{𝗤}(x,⊥,q_2) = (x, q_2)$ and  $μ^{𝗤}(x,q_1,⊥) = (x, q_1)$—and the unit defined by $η^{𝗤}(x) = (x, ⊥)$. This is a \kl{commutative magmad} but not a \kl{left-exchanging magmad}.
\end{proof}

\begin{proposition}
  \label[proposition]{prop:left-exchanging}
  \AP \kl{Sesquilaws} induce \kl{left-exchanging magmads}.
\end{proposition}

\subsection{Affine, Relevant, and Normal Magmads}
\label{sec:affine-relevant}

The category of sets and functions has the structure to copy and discard information; explicitly, there exist transformations $δ ፡ X → X × X$ and $ε ፡ X → 1$ that are natural. \kl{Monads} and \kl{magmads} do not necessarily preserve this structure: some effects cannot be copied or discarded without changing their meaning; the following two equations (①, ②) do not always hold. 

For instance, stochasticity cannot be copied (② fails)—throwing two coins is different from throwing a single coin and copying the result; and partiality cannot be discarded (① fails)—even if we discard the result of a partial function, its domain remains unaltered.

\[
\begin{minipage}{0.35\textwidth}
\begin{lstlisting}[language=Racket, numbers=none]
(do return ())
\end{lstlisting}
\end{minipage}
\quad \overset{\raisebox{0.2em}{①}}{=}\quad
\begin{minipage}{0.35\textwidth}
\begin{lstlisting}[language=Racket, numbers=none]
(do (x1 ...) <- F
    return ())
\end{lstlisting}
\end{minipage}
\]
\[
\begin{minipage}{0.35\textwidth}
\begin{lstlisting}[language=Racket, numbers=none]
(do (x1 ...) <- F
    (y1 ...) <- F
    return (x1 ... y1 ...))
\end{lstlisting}
\end{minipage}
\quad \overset{\raisebox{0.2em}{②}}{=}\quad
\begin{minipage}{0.35\textwidth}
\begin{lstlisting}[language=Racket, numbers=none]
(do (x1 ...) <- F
    return (x1 ... x1 ...))
\end{lstlisting}
\end{minipage}
\]

Instead, stochasticity can be discarded ①, and partiality can be copied ② (\Cref{distribution-affine-maybe-relevant}): monads that can be discarded are called \emph{\kl{affine}}; monads that can be copied are called \emph{\kl{relevant}} (\Cref{def:affine-relevant-magmad}).

\begin{definition}[Affine \& relevant magmads, c.f.~{{\cite{jacobs94}}}]
  \label[definition]{def:affine-relevant-magmad}
  \AP An \intro{affine magmad} and a \intro{relevant magmad} are \kl{commutative magmads}, $(𝐓, \smash{μ^{𝐓}}, \smash{η^{𝐓}}, \smash{ψ^{𝐓}})$, such that—respectively—the first or the second of the following diagrams commute.
  \begin{align*}
\begin{tikzcd}[ampersand replacement=\&]
	{𝐓X} \&\& {𝐓1} \\
	1
	\arrow[""{name=0, anchor=center, inner sep=0}, "{𝐓\varepsilon}", from=1-1, to=1-3]
	\arrow["\varepsilon"', from=1-1, to=2-1]
	\arrow["{\eta^{𝐓}}"', from=2-1, to=1-3]
	\arrow["{①}"{description}, draw=none, from=0, to=2-1]
\end{tikzcd}
	\quad
  \begin{tikzcd}[ampersand replacement=\&]
	{𝐓X} \& {𝐓(X \otimes X)} \\
	{𝐓X \otimes 𝐓X}
	\arrow[""{name=0, anchor=center, inner sep=0}, "{𝐓\delta}", from=1-1, to=1-2]
	\arrow["\delta"', from=1-1, to=2-1]
	\arrow["\psi"', from=2-1, to=1-2]
	\arrow["{②}"{description}, draw=none, from=0, to=2-1]
  \end{tikzcd}
\end{align*}
\end{definition}

\begin{proposition}
  \label[proposition]{distribution-affine-maybe-relevant}
  The \kl{distribution monad} $(𝐃,μ^{𝐃},η^{𝐃})$ is \kl{affine}; the \kl{maybe monad} $(𝐌,μ^{𝐌},η^{𝐌})$ is \kl{relevant}.
   \label[proposition]{prop:magmad-normalized-distributions}
  \AP The \kl{magmad} of \kl{normalized distributions}, $(𝐌𝐃,\smash{μ^{𝐌𝐃}_{∘}},\smash{η^{𝐌𝐃}_{∘}})$ is a \kl{commutative magmad} that is not \kl{affine} nor \kl{relevant}.
\end{proposition}

\kl{Normalized distributions} form a \kl{magmad} that is not \kl{affine} nor \kl{relevant}. Totality can be refined to account for updates: the resulting notion is that of a \kl{normal magmad}—or, in the terminology of \emph{partial Markov categories}, a magmad whose effects are \kl{quasitotal}~\cite{dilavore23:evidential,partialMarkov25,mohammed2025partializations}.

\[
\begin{minipage}{0.35\textwidth}
\begin{lstlisting}[language=Racket, numbers=none]
(do (x1 ...) <- F
    (y1 ...) <- F
    return (x1 ...))
\end{lstlisting}
	\end{minipage}
	\quad \overset{③}{=}\quad
	\begin{minipage}{0.35\textwidth}
\begin{lstlisting}[language=Racket, numbers=none]
(do (x1 ...) <- F
    return (x1 ...))
\end{lstlisting}
\end{minipage}
\]

\begin{definition}[Normal magmad]
\AP A \intro{normal magmad} is a \kl{commutative magmad}, $(𝐓, \smash{μ^{𝐓}}, \smash{η^{𝐓}}, \smash{ψ^{𝐓}})$, such that the following diagram commutes.
\[\begin{tikzcd}[ampersand replacement=\&]
	{𝐓X} \& {𝐓X} \\
	{𝐓X \otimes 𝐓X} \& {𝐓(X \otimes X)}
	\arrow[no head, from=1-1, to=1-2]
	\arrow[""{name=0, anchor=center, inner sep=0}, "\delta"', from=1-1, to=2-1]
	\arrow["\psi", from=2-1, to=2-2]
	\arrow[""{name=1, anchor=center, inner sep=0}, "{𝐓(X \otimes \varepsilon)}"', from=2-2, to=1-2]
	\arrow["{③}"{description}, draw=none, from=0, to=1]
\end{tikzcd}\]
\end{definition}

\begin{proposition}[Affine-relevant sesquilaw]
  \AP An \intro{affine-relevant sesquilaw} is a \kl{monoidal sesquilaw}, $(𝐒,𝐓,m,n)$, of an \kl{affine monad}, $(𝐓,\smash{μ^𝐓},\smash{η^𝐓})$ over a \kl{relevant monad}, $(𝐒,\smash{μ^𝐒},\smash{η^𝐒})$.
  \AP Any \kl{affine-relevant sesquilaw} $(𝐒,𝐓,m,n)$ induces a \kl{quasitotal magmad} $(𝐒𝐓,\smash{μ^{𝐒𝐓}_{n}},\smash{η^{𝐒𝐓}_{n}})$. 
\end{proposition}

\begin{corollary}
  \label[corollary]{cor:normal-magmad-examples}
  The \kl{magmad} of \kl{normalized distributions}, $(𝐌𝐃,\smash{μ^{𝐌𝐃}_{∘}},\smash{η^{𝐌𝐃}_{∘}})$ is a \kl{normal magmad}. The \kl{monad} of \kl{subdistributions}, $(𝐃𝐌,\smash{μ^{𝐌𝐃}_{•}},\smash{η^{𝐌𝐃}_{•}})$ is not a \kl[normal magmad]{normal monad}.
\end{corollary}

\subsection{Reasoning—Pearl's Front-door and Back-door criteria}
\label{sec:pearl-frontdoor-backdoor}

The axioms of the probabilistic magmadic metalanguage (\Cref{fig:axioms-metalanguage}) allow us to prove basic lemmas for synthetic causal inference (c.f.~\cite{pearl2009causality}). As a feature, they distinguish between interventions and observations by reparenthesizing. 

\begin{remark}
Jacobs, Széles, and Stein \cite{jacobs2025compositional} have recently shown how \kl{partial Markov categories} \cite{dilavore23:evidential}, extended with normalization boxes \cite{lorenz2023causal}, can be applied to synthetic causality (c.f.~\cite{yin22markovcausal,Pie23}). We prove similar results but, for the first time, using a normalized semantics and in a metalanguage with axioms: this seems to render the proofs slightly simpler and paves the way for generalization.
\end{remark}

\newsavebox{\codeBackdoorOne}
\begin{lrbox}{\codeBackdoorOne}
\begin{minipage}{0.25\textwidth}
\begin{lstlisting}[language=Racket,numbers=none,mathescape]
(do (u x y) <- P
    return (u x y))
\end{lstlisting}
\end{minipage}	
\end{lrbox}

\newsavebox{\codeBackdoorTwo}
\begin{lrbox}{\codeBackdoorTwo}
\begin{minipage}{0.25\textwidth}
\begin{lstlisting}[language=Racket,numbers=none,mathescape]
(do (u) <- F
    (x) <- (G x)
    (y) <- (H x u)
    return (u x y))
\end{lstlisting}
\end{minipage}	
\end{lrbox}

\newsavebox{\codeBackdoorThree}
\begin{lrbox}{\codeBackdoorThree}
\begin{minipage}{0.25\textwidth}
\begin{lstlisting}[language=Racket,numbers=none,mathescape]
(do (u) <- F
    (y) <- (H a u)
    return (y))
\end{lstlisting}
\end{minipage}	
\end{lrbox}

\newsavebox{\codeBackdoorFour}
\begin{lrbox}{\codeBackdoorFour}
\begin{minipage}{0.35\textwidth}
\begin{lstlisting}[language=Racket,numbers=none,mathescape]
(do
  (u0 x0 y0) <- P
  y <- (do
     (u x y) <- P
     () <- (observe u u0)
     () <- (observe x a)
     return (y))
  return (y))
\end{lstlisting}
\end{minipage}	
\end{lrbox}

\newsavebox{\codeBackproofOne}
\begin{lrbox}{\codeBackproofOne}
\begin{minipage}{0.25\textwidth}
\begin{lstlisting}[language=Racket,numbers=none,mathescape]
(do (u0 x0 y0) <- P
    return (u0))
\end{lstlisting}
\end{minipage}	
\end{lrbox}

\newsavebox{\codeBackproofTwo}
\begin{lrbox}{\codeBackproofTwo}
\begin{minipage}{0.25\textwidth}
\begin{lstlisting}[language=Racket,numbers=none,mathescape]
(do (u0) <- F
    (x) <- (G x)
    (y) <- (H x u)
    return (u0 x y))
\end{lstlisting}
\end{minipage}	
\end{lrbox}

\newsavebox{\codeBackproofThree}
\begin{lrbox}{\codeBackproofThree}
\begin{minipage}{0.25\textwidth}
\begin{lstlisting}[language=Racket,numbers=none,mathescape]
(do (u0) <- F
    return (u0 x y))
\end{lstlisting}
\end{minipage}	
\end{lrbox}

\newsavebox{\codeBackproofFour}
\begin{lrbox}{\codeBackproofFour}
\begin{minipage}{0.3\textwidth}
\begin{lstlisting}[language=Racket,numbers=none,mathescape]
(do (u x y) <- P
    () <- (observe u u0)
    () <- (observe x a)
    return (y))
\end{lstlisting}
\end{minipage}	
\end{lrbox}

\newsavebox{\codeBackproofFive}
\begin{lrbox}{\codeBackproofFive}
\begin{minipage}{0.3\textwidth}
\begin{lstlisting}[language=Racket,numbers=none,mathescape]
(do (u) <- F
    (x) <- (G x)
    (y) <- (H x u)
    () <- (observe u u0)
    () <- (observe x a)
    return (y))
\end{lstlisting}
\end{minipage}	
\end{lrbox}

\newsavebox{\codeBackproofSix}
\begin{lrbox}{\codeBackproofSix}
\begin{minipage}{0.3\textwidth}
\begin{lstlisting}[language=Racket,numbers=none,mathescape]
(do
  (u) <- F
  () <- (observe u u0)
  (x) <- (G x)
  (y) <- (H x u0)
  () <- (observe x a)
  return (y))
\end{lstlisting}
\end{minipage}	
\end{lrbox}

\newsavebox{\codeBackproofSeven}
\begin{lrbox}{\codeBackproofSeven}
\begin{minipage}{0.3\textwidth}
\begin{lstlisting}[language=Racket,numbers=none,mathescape]
(do
  (x) <- (G x)
  (y) <- (H x u0)
  () <- (observe x a)
  return (y))
\end{lstlisting}
\end{minipage}	
\end{lrbox}

\newsavebox{\codeBackproofEight}
\begin{lrbox}{\codeBackproofEight}
\begin{minipage}{0.35\textwidth}
\begin{lstlisting}[language=Racket,numbers=none,mathescape]
(do (u0 x0 y0) <- P
    y <- (do
       (u x y) <- P
       () <- (observe u u0)
       () <- (observe x a)
       return (y))
    return (y))
\end{lstlisting}
\end{minipage}	
\end{lrbox}

\newsavebox{\codeBackproofNine}
\begin{lrbox}{\codeBackproofNine}
\begin{minipage}{0.35\textwidth}
\begin{lstlisting}[language=Racket,numbers=none,mathescape]
(do (u0) <- F
    y <- (do
       (x) <- (G x)
       () <- (observe x a)
       (y) <- (H x u0)
       return (y))
    return (y))
\end{lstlisting}
\end{minipage}	
\end{lrbox}

\newsavebox{\codeBackproofTen}
\begin{lrbox}{\codeBackproofTen}
\begin{minipage}{0.35\textwidth}
\begin{lstlisting}[language=Racket,numbers=none,mathescape]
(do (u0) <- F
    (x) <- (do
       (x) <- (G x)
       () <- (observe x a)
       return (x))
    (y) <- (H x u0)
    return (y))
\end{lstlisting}
\end{minipage}	
\end{lrbox}

\newsavebox{\codeBackproofEleven}
\begin{lrbox}{\codeBackproofEleven}
\begin{minipage}{0.35\textwidth}
\begin{lstlisting}[language=Racket,numbers=none,mathescape]
(do (u0) <- F
    (y) <- (H a u0)
    return (y))
\end{lstlisting}
\end{minipage}	
\end{lrbox}

\newsavebox{\codeFrontOne}
\begin{lrbox}{\codeFrontOne}
\begin{minipage}{0.25\textwidth}
\begin{lstlisting}[language=Racket,numbers=none,mathescape]
(do (s t c) <- P
    return (s t c))
\end{lstlisting}
\end{minipage}	
\end{lrbox}

\newsavebox{\codeFrontTwo}
\begin{lrbox}{\codeFrontTwo}
\begin{minipage}{0.25\textwidth}
\begin{lstlisting}[language=Racket,numbers=none,mathescape]
(do (g) <- G
    (s) <- (S g)
    (t) <- (T s)
    (c) <- (C t g)
    return (s t c))
\end{lstlisting}
\end{minipage}	
\end{lrbox}

\newsavebox{\codeFrontThree}
\begin{lrbox}{\codeFrontThree}
\begin{minipage}{0.25\textwidth}
\begin{lstlisting}[language=Racket,numbers=none,mathescape]
(do (g) <- G  
    (t) <- (T a)
    (c) <- (C t g)
    return (c))
\end{lstlisting}
\end{minipage}	
\end{lrbox}

\newsavebox{\codeFrontFour}
\begin{lrbox}{\codeFrontFour}
\begin{minipage}{0.44\textwidth}
\begin{lstlisting}[language=Racket,numbers=none,mathescape]
(do (s t0 c0) <- P
    (t) <- (do (s1 t c1) <- P
               () <- (observe s1 a)
               return (t))
    (c) <- (do (s2 t2 c) <- P
               () <- (observe s2 s)
               () <- (observe t2 t)
               return (c))
    return (c))
\end{lstlisting}
\end{minipage}	
\end{lrbox}

\newsavebox{\codeFrontFive}
\begin{lrbox}{\codeFrontFive}
\begin{minipage}{0.25\textwidth}
\begin{lstlisting}[language=Racket,numbers=none,mathescape]
(do (s t0 c0) <- P
    return (s))
\end{lstlisting}
\end{minipage}	
\end{lrbox}

\newsavebox{\codeFrontSix}
\begin{lrbox}{\codeFrontSix}
\begin{minipage}{0.25\textwidth}
\begin{lstlisting}[language=Racket,numbers=none,mathescape]
(do (g) <- G
    (s) <- (S g)
    (t) <- (T s)
    (c) <- (C t g)
    return (s))
\end{lstlisting}
\end{minipage}	
\end{lrbox}

\newsavebox{\codeFrontSeven}
\begin{lrbox}{\codeFrontSeven}
\begin{minipage}{0.25\textwidth}
\begin{lstlisting}[language=Racket,numbers=none,mathescape]
(do (g) <- G
    (s) <- (S g)
    return (s))
\end{lstlisting}
\end{minipage}	
\end{lrbox}

\newsavebox{\codeFrontEight}
\begin{lrbox}{\codeFrontEight}
\begin{minipage}{0.35\textwidth}
\begin{lstlisting}[language=Racket,numbers=none,mathescape]
(do (s1 t c1) <- P
    () <- (observe s1 a)
    return (t))
\end{lstlisting}
\end{minipage}	
\end{lrbox}

\newsavebox{\codeFrontNine}
\begin{lrbox}{\codeFrontNine}
\begin{minipage}{0.35\textwidth}
\begin{lstlisting}[language=Racket,numbers=none,mathescape]
(do (g) <- G
    (s1) <- (S g)
    (t) <- (T s1)
    (c1) <- (C t g)
    () <- (observe s1 a)
    return (t))
\end{lstlisting}
\end{minipage}	
\end{lrbox}

\newsavebox{\codeFrontTen}
\begin{lrbox}{\codeFrontTen}
\begin{minipage}{0.35\textwidth}
\begin{lstlisting}[language=Racket,numbers=none,mathescape]
(do (g) <- G
    (s1) <- (S g)
    () <- (observe s1 a)
    (t) <- (T a)
    return (t))
\end{lstlisting}
\end{minipage}	
\end{lrbox}

\newsavebox{\codeFrontEleven}
\begin{lrbox}{\codeFrontEleven}
\begin{minipage}{0.35\textwidth}
\begin{lstlisting}[language=Racket,numbers=none,mathescape]
(do (t) <- (T a)
    return (t))
\end{lstlisting}
\end{minipage}	
\end{lrbox}

\newsavebox{\codeFrontTwelve}
\begin{lrbox}{\codeFrontTwelve}
\begin{minipage}{0.35\textwidth}
\begin{lstlisting}[language=Racket,numbers=none,mathescape]
(do (s2 t2 c) <- P
    () <- (observe s2 s)
    () <- (observe t2 t)
    return (c))
\end{lstlisting}
\end{minipage}	
\end{lrbox}

\newsavebox{\codeFrontThirteen}
\begin{lrbox}{\codeFrontThirteen}
\begin{minipage}{0.35\textwidth}
\begin{lstlisting}[language=Racket,numbers=none,mathescape]
(do (g) <- G
    (s2) <- (S g)
    (t2) <- (T s2)
    (c) <- (C t2 g)
    () <- (observe s2 s)
    () <- (observe t2 t)
    return (c))
\end{lstlisting}
\end{minipage}	
\end{lrbox}

\newsavebox{\codeFrontFourteen}
\begin{lrbox}{\codeFrontFourteen}
\begin{minipage}{0.35\textwidth}
\begin{lstlisting}[language=Racket,numbers=none,mathescape]
(do (g) <- G
    (s2) <- (S g)
    (g) <- ((d S) s2)
    (t2) <- (T s2)
    (c) <- (C t2 g)
    () <- (observe s2 s)
    () <- (observe t2 t)
    return (c))
\end{lstlisting}
\end{minipage}	
\end{lrbox}

\newsavebox{\codeFrontFifteen}
\begin{lrbox}{\codeFrontFifteen}
\begin{minipage}{0.35\textwidth}
\begin{lstlisting}[language=Racket,numbers=none,mathescape]
(do (g) <- G
    (s2) <- (S g)
    () <- (observe s2 s)
    (g) <- ((d S) s)
    (t2) <- (T s)
    (c) <- (C t2 g)
    () <- (observe t2 t)
    return (c))
\end{lstlisting}
\end{minipage}	
\end{lrbox}

\newsavebox{\codeFrontSixteen}
\begin{lrbox}{\codeFrontSixteen}
\begin{minipage}{0.35\textwidth}
\begin{lstlisting}[language=Racket,numbers=none,mathescape]
(do (g) <- ((d S) s)
    (t2) <- (T s)
    (c) <- (C t2 g)
    () <- (observe t2 t)
    return (c))
\end{lstlisting}
\end{minipage}	
\end{lrbox}

\newsavebox{\codeFrontSeventeen}
\begin{lrbox}{\codeFrontSeventeen}
\begin{minipage}{0.35\textwidth}
\begin{lstlisting}[language=Racket,numbers=none,mathescape]
(do (t2) <- (T s)
    () <- (observe t2 t)
    (g) <- ((d S) s)
    (c) <- (C t g)
    return (c))
\end{lstlisting}
\end{minipage}	
\end{lrbox}

\newsavebox{\codeFrontEighteen}
\begin{lrbox}{\codeFrontEighteen}
\begin{minipage}{0.35\textwidth}
\begin{lstlisting}[language=Racket,numbers=none,mathescape]
(do (g) <- ((d S) s)
    (c) <- (C t g)
    return (c))
\end{lstlisting}
\end{minipage}	
\end{lrbox}

\newsavebox{\codeFrontNineteen}
\begin{lrbox}{\codeFrontNineteen}
\begin{minipage}{0.35\textwidth}
\begin{lstlisting}[language=Racket,numbers=none,mathescape]
(do
  (t) <- (T a)  
  (g) <- G
  (s) <- (S g)
  (c) <- (do
    (g) <- ((d S) s)
    (c) <- (C t g)
    return (c))
  return (c))
\end{lstlisting}
\end{minipage}	
\end{lrbox}

\newsavebox{\codeFrontTwenty}
\begin{lrbox}{\codeFrontTwenty}
\begin{minipage}{0.35\textwidth}
\begin{lstlisting}[language=Racket,numbers=none,mathescape]
(do (t) <- (T a)  
    (g) <- G
    (s) <- (S g)
    (g) <- ((d S) s)
    (c) <- (C t g)
    return (c))
\end{lstlisting}
\end{minipage}	
\end{lrbox}

\newsavebox{\codeFrontTwentyone}
\begin{lrbox}{\codeFrontTwentyone}
\begin{minipage}{0.35\textwidth}
\begin{lstlisting}[language=Racket,numbers=none,mathescape]
(do (g) <- G
    (t) <- (T a)  
    (c) <- (C t g)
    return (c))
\end{lstlisting}
\end{minipage}	
\end{lrbox}

\begin{theoremrep}[Front-door adjustment formula {{\cite{pearl2009causality}}}]
  \label{thm:frontdoor}
  \AP Given a magmad, $(𝗡, μ^{𝗡}, η^{𝗡})$, satisfying the principles in \Cref{fig:axioms-metalanguage}, let a joint state, $p ፡ 1 → X ⊗
  Z ⊗ Y$, admit the following factorization into \kl{total} morphisms where, moreover, both $t$ and $(g ⨾ s)$ below have \kl{full support}.
  \begin{align*}
    \usebox{\codeFrontOne}
    \quad \overset{(i)}{=} \quad
    \usebox{\codeFrontTwo}
  \end{align*}
  Then, the following equation holds. We define \lstinline[language=Racket,  
basicstyle={\ttfamily\color{nord0}}]|(frontDoor P a)| as its right-hand side.
  \begin{align*}
    \usebox{\codeFrontThree}
    \quad \overset{(i)}{=} \quad
    \usebox{\codeFrontFour}
  \end{align*}
  In other words, an intervention on the variable\, $X$ can be rewritten as a
  composition, in the \kl{Markov magmoid}, of the observational data.
\end{theoremrep}
\begin{proofsketch}
  We simplify each one of the three parts of the equation. For the first part, \emph{(i)} we use the hypothesis, and \kl{totality} of $\mathsf{C}$. For the second part, \emph{(ii)} we also use \kl{Frobenius} and the \kl{full-support assumption}. For the third part, \emph{(iii)} we employ the \kl{disintegration} axiom. See the details in the Appendix.
  \begin{gather*}
    \usebox{\codeFrontFive}
    \quad \overset{(i)}{=} \quad
    \usebox{\codeFrontSeven} \\
    \usebox{\codeFrontEight}
    \quad \overset{(ii)}{=} \quad 
    \usebox{\codeFrontEleven} \\
    \usebox{\codeFrontTwelve}
    \quad \overset{(iii)}{=} \quad
    \usebox{\codeFrontEighteen}
  \end{gather*}

  With these equations, let us address the main claim. We \emph{(i)} substitute according to the previous three equations, \emph{(ii)} use the totality of $g$ and $c$, and \emph{(iii)} apply the \kl{disintegration} axiom.
  \begin{align*}
    &\usebox{\codeFrontFour}
    \quad \overset{(i)}{=} \quad 
    \usebox{\codeFrontNineteen} 
    \\ & \quad \overset{(ii)}{=} \quad  
    \usebox{\codeFrontTwenty}
    \quad \overset{(iii)}{=} \quad 
    \usebox{\codeFrontTwentyone}
  \end{align*}
\end{proofsketch}
\begin{proof}
  Let us simplify each one of the three factors. For the first factor, we use \emph{(i)} the hypothesis, and \emph{(ii)} \kl{totality} of $\mathsf{C}$.
  \begin{gather*}
    \usebox{\codeFrontFive}
    \quad \overset{(i)}{=} \quad
    \usebox{\codeFrontSix}
    \quad \overset{(ii)}{=} \quad
    \usebox{\codeFrontSeven} 
  \end{gather*}
  
  For the second factor computes the marginal of a conditional distribution. We use \emph{(i)} the hypothesis; \emph{(ii)} \kl[total]{totality} of $\mathsf{c}$ and \kl{Frobenius}; and \emph{(iii)} that $(g ⨾ s)$ has \kl{full support}.
  \begin{align*}
    & \usebox{\codeFrontEight}
    \quad \overset{(i)}{=} \quad 
    \usebox{\codeFrontNine}    
    \quad \overset{(ii)}{=} \quad \\
    & \usebox{\codeFrontTen}
    \quad \overset{(iii)}{=} \quad 
    \usebox{\codeFrontEleven}
  \end{align*}

  The third factor can be simplified as follows. We use \emph{(i)} the hypothesis, \emph{(ii)} the \kl{disintegration axiom}, \emph{(iii)} the \kl{Frobenius} axiom and \kl{exchange}, \emph{(iv)} that $(g ⨾ s)$ has \kl{full support}, \emph{(v)} \kl{exchange}, and \emph{(vi)} that $(t)$ has \kl{full support}.
  \begin{align*}
    & \usebox{\codeFrontTwelve}
    \quad \overset{(i)}{=} \quad 
    \usebox{\codeFrontThirteen}
    \quad \overset{(ii)}{=} \quad \\
    & \usebox{\codeFrontFourteen}
    \quad \overset{(iii)}{=} \quad 
    \usebox{\codeFrontFifteen}
    \quad \overset{(iv)}{=} \quad \\
    & \usebox{\codeFrontSixteen}
    \quad \overset{(v)}{=} \quad 
    \usebox{\codeFrontSeventeen} 
    \quad \overset{(vi)}{=} \quad \\
    & \usebox{\codeFrontEighteen}
  \end{align*}

  Finally, let us address the main claim. We \emph{(i)} substitute according to the previous three equations, \emph{(ii)} use the totality of $g$ and $c$, and \emph{(iii)} apply the \kl{disintegration} axiom.
  \begin{align*}
    &\usebox{\codeFrontFour}
    \quad \overset{(i)}{=} \quad 
    \usebox{\codeFrontNineteen} 
    \quad \overset{(ii)}{=} \quad \\ 
    &\usebox{\codeFrontTwenty}
    \quad \overset{(iii)}{=} \quad 
    \usebox{\codeFrontTwentyone}
  \end{align*}
  This concludes the proof: the intervention can be written exclusively in terms of $\mathsf{P}$.
\end{proof}

\begin{theoremrep}[Back-door adjustment formula {{\cite{pearl2009causality}}}]
  \label{prop:backdoorAdjustment}
  Given a magmad, $(𝗡, μ^{𝗡}, η^{𝗡})$, satisfying the principles in \Cref{fig:axioms-metalanguage}, let a joint state, $\mathsf{P} ∈ 𝗡(U × X × Y)$, admit the following factorization into \kl{total} morphisms—where, moreover, $\mathsf{F}$ has \kl{full support}.
  \begin{align*}
\usebox{\codeBackdoorOne}
\quad &\overset{(i)}{=} \quad
\usebox{\codeBackdoorTwo}
  \end{align*}
  Then, the following equation holds. We define \lstinline[language=Racket,  
basicstyle={\ttfamily\color{nord0}}]|(backDoor P a)| as the right-hand side.
  \begin{align*}
\usebox{\codeBackdoorThree}
\quad &\overset{(ii)}{=} \quad
\usebox{\codeBackdoorFour} \\
  \end{align*}
\end{theoremrep}
\begin{proofsketch}
  We apply the same reasoning principles as for the \kl{front-door formula}. See details in the Appendix.
\end{proofsketch}
\begin{proof}
  Let us first simplify the uppermost part of the diagram. We use \emph{(i)} the hypothesis, and \emph{(ii)} the \kl[total]{totality} of $\mathsf{G}$ and $\mathsf{H}$.
  \begin{align*}
    \usebox{\codeBackproofOne}
    \quad \overset{(i)}{=} \quad
    \usebox{\codeBackproofTwo}
    \quad \overset{(ii)}{=} \quad
    \usebox{\codeBackproofThree}
  \end{align*}

  Let us simplify the second part of the diagram.  We now use \emph{(i)} the hypothesis, \emph{(ii)} the \kl{Frobenius equation} and \kl{exchange}, and \emph{(iii)} the assumption that $\mathsf{F}$ has \kl{full support}.
  \begin{align*}
    & \usebox{\codeBackproofFour}
    \quad \overset{(i)}{=} \quad
    \usebox{\codeBackproofFive} 
    \quad \overset{(ii)}{=} \quad \\
    & \usebox{\codeBackproofSix} 
    \quad \overset{(iii)}{=} \quad
    \usebox{\codeBackproofSeven}
  \end{align*}

  Let us conclude the proof. We \emph{(i)} substitute both simplifications, \emph{(ii)} we use the totality of $\mathsf{H}$, and \emph{(iii)} the \kl{Frobenius equation} and the totality of $\mathsf{G}$.
  \begin{align*}
    & \usebox{\codeBackproofEight}
    \quad \overset{(i)}{=} \quad
    \usebox{\codeBackproofNine}
    \quad \overset{(ii)}{=} \quad \\
    & \usebox{\codeBackproofTen} 
    \quad \overset{(ii)}{=} \quad
    \usebox{\codeBackproofEleven}
  \end{align*}

  This concludes the proof: we have shown that the interventional formula can be written in terms of $\mathsf{P}$, exclusively.
\end{proof}

\newsavebox{\codeFrontdoorLeft}
\begin{lrbox}{\codeFrontdoorLeft}
\begin{minipage}{0.35\textwidth}
\begin{lstlisting}[language=Racket,numbers=none,mathescape]
(frontDoor P a)
\end{lstlisting}
\end{minipage}	
\end{lrbox}

From this perspective, causal inference is the resolution of equations in the magmoidal metalanguage of the \kl{normalized distribution monad}. Alternatively to most probabilistic programming literature, we do not use a normalization operator, but an ability to modulate associativity.

%

\section{Conclusions}

\kl{Normalization} has been under-explored in denotational probabilistic semantics and categorical probability
theory; most inference semantics either ignore \kl{normalization}—preferring \kl{substochastic kernels} or unnormalized kernels instead—or treat it as an ad-hoc operator or a semantic feature.
As a result, most denotational semantic universes for probabilistic programming are not \emph{normalized by construction}. We have introduced the \kl{magmad} of \kl[normalized distributions]{normalized distributions} as a first example of \emph{normalized-by-construction} denotational semantics for probabilistic inference.

Our proposed explanation for this gap is that \kl{normalization} arises from a failure of associativity, and failures of associativity are counterintuitive. For instance, we could naïvely say that, to solve an inference problem, one must (1) set up a prior distribution, (2) compute a stochastic kernel, and (3) update the prior with the observation. This description misses an essential point: how to associate these instructions or, in other words, \emph{when to normalize}. Any normalized probablistic semantics needs to clearly address this point, and a revision of existing categorical probabilitic semantics to understand their potential non-associative behaviour is warranted.

While non-associativity could seem too high of a price to pay, being explicit about associativity has conceptual advantages: it makes it possible to derive multiple properties of \kl{normalization} from a few axioms and it allows us to distinguish between observations and interventions as two different parenthesizations of the same expression.
That \kl{normalization} is not a \kl{distributive law} may be folklore, but it seems absent from the literature. We went one step further and classified this particular failure of normalization by introducing \kl{sesquilaws}: \kl{almost-distributive laws} that are multiplicative only up to idempotent.

\kl{Sesquilaws} that are monoidal induce not only \kl{monoidal magmads} but \kl{left-exchanging magmads}, another newly introduced notion refining the non-associative case. For \kl{left-exchanging magmads}, we provided both a variant of the monadic metalanguage that is not associative but preserves the line-exchange property. Remarkably, these seem to mostly coincide—although not exactly—with the existing string diagrammatic calculi for \kl{normalization} over \emph{partial Markov categories} \cite{dilavore23:evidential,lorenz2023causal,jacobs2025compositional}.

We claimed that there exist multiple ways of updating, that this is a non-associative phenomenon, and that an expressive probabilistic programming language must allow them in some form; in this, we follow Jacobs' description of Pearl's and Jeffrey's modes of update \cite{jacobs2019mathematics}.
Finally, we proposed a different approach to categorical causality: we saw how solving causality problems corresponds to solving equations on a magmadic metalanguage. While Pearl's \emph{do-calculus}, more than a synthetic axiomatization, is a collection of rules on top of probabilistic reasoning, we propose an algebraic presentation of the rules needed to discuss causality problems. Further work on automating these solutions—e.g. via rewriting, showing completeness for \emph{identifiability}—is warranted and not covered by this manuscript; it is a promising avenue for a categorical semantics of causality and causal programming.

More generally, we reassert—as weak distributive laws, duploids, or mass-chance
interpretations do—that failing distributive laws contain, in many cases,
mathematical structure worth studying on its own.

\begin{toappendix}

\end{toappendix}

\begin{toappendix}
  \subsection*{Single-point extensions}

\begin{definition}[Extension]
  Given two monads, $(𝐒, \smash{μ^{𝐒}}, \smash{η^{𝐒}})$ and $(𝐓, \smash{μ^{𝐓}}, \smash{η^{𝐓}})$, an $𝐒$-\intro{extension} of the $𝐓$-algebra $(X,α)$ is a $𝐓$-algebra of the form \((𝐒 X, α^{\#})\) such that \(α ⨾ η^{𝐒} = 𝐓 (η^{𝐒}) ⨾ α^{\#}\).
  \[\begin{tikzcd}[ampersand replacement=\&]
    𝐓 X \rar{𝐓 η^{𝐒}} \dar[swap]{α} \& 𝐓𝐒 X \dar{α^{\#}} \\
    X  \rar{η^{𝐒}} \& 𝐒 X
  \end{tikzcd}\]
\end{definition}

\begin{remark}
  A \emph{convex algebra} is an algebra for the finitely-supported distribution monad $(𝐃)$. A \emph{one-point extension} of a convex algebra—in the sense of Sokolova and Woracek \cite{sokolova18:termination}—is precisely an $𝐌$-\kl{extension}, an extension over the \kl{maybe monad}.
\end{remark}

\begin{theorem}%
  [Black-hole is unique {{\cite[Theorem~5.3]{sokolova18:termination}}}]%
  \label[proposition]{prop:one-point-extensions}
  The black-hole one-point extension, defined by
  \[\alpha^{\#}(d) = \begin{cases} d, & \text{if }d(\bot) = 0, \\ \bot, & \text{if } d(\bot)>0, \end{cases}\]
  is the only functorial one-point extension of convex algebras.
\end{theorem}

It is well-known that a \kl{distributive law}, \(ψ ፡ 𝐓𝐒 → 𝐒𝐓\), gives a lifting of the first monad, \(𝐒\), to the category of algebras of second, \(𝐓\)~\cite{beck1969distributive}. The lifted monad acts on algebras by extending them. Let us recall this construction here.

\begin{proposition}[Extensions from distributive laws]%
  \label[proposition]{prop:extensions-from-distributivity}
  Any \kl{distributive law}, \(ψ ፡ 𝐓𝐒 → 𝐒𝐓\), determines a functorial \(𝐒\)-\kl{extension} of\, \(𝐓\)-algebras.
\end{proposition}
\begin{proof}
  Given a \kl{distributive law}, $ψ ፡ 𝐓𝐒 → 𝐒𝐓$ and a $T$-algebra $(X,α)$, consider the morphism $α^{\#} = ψ_{X} ⨾ 𝐒(α)$. The lifted monad, $ℓ(𝐒)$, following Beck's work~\cite{beck1969distributive}, acts on objects $(X, α)$ by $ℓ(𝐒)(X, α) = (𝐒(X), α^{\#})$. 
  Let us now check that \((𝐒 X, α^{\#})\) is an \(𝐒\)-\kl{extension} of $(X, α)$. We use \emph{(i)} Beck's construction, \emph{(ii)} $𝐒$-unitality of the \kl{distributive law}, and \emph{(iii)} naturality.
  \begin{align*}
     𝐓(η^{𝐒}) ⨾ α^{\#}
    \overset{\emph{(i)}}{=} 𝐓(η^{𝐒}) ⨾ ψ ⨾ 𝐒(α)
     \overset{\emph{(ii)}}{=} η^{S} ⨾ 𝐒(α) 
     \overset{\emph{(iii)}}{=} α ⨾ η^{𝐒}.
  \end{align*}
  This concludes the proof.
\end{proof}
\end{toappendix}

\bibliographystyle{halpha} 
\bibliography{main}

\end{document}